\documentclass[10pt]{amsart}
\usepackage[cp1251]{inputenc}
\usepackage[english]{babel}
\usepackage{amsmath}
\usepackage{amssymb}
\usepackage{amsfonts}

\usepackage[linktocpage=true, colorlinks=true, linkcolor=blue, citecolor=blue, urlcolor=blue]{hyperref}

\begin{document}
\renewcommand{\refname}{References}

\thispagestyle{empty}

\title[Development and 
Application of the Fourier Method]
{Development and Application of the Fourier Method 
to the Mean-Square Approximation
of Iterated Ito and Stratonovich Stochastic Integrals}
\author[D.F. Kuznetsov]{Dmitriy F. Kuznetsov}
\address{Dmitriy Feliksovich Kuznetsov
\newline\hphantom{iii} Peter the Great Saint-Petersburg Polytechnic University,
\newline\hphantom{iii} Polytechnicheskaya ul., 29,
\newline\hphantom{iii} 195251, Saint-Petersburg, Russia}%
\email{sde\_kuznetsov@inbox.ru}

\maketitle {\small
\begin{quote}
\noindent{\sc Abstract.} 
The article is devoted to the mean-square approximation
of iterated Ito and Stratonovich stochastic integrals in the context
of the numerical integration of Ito stochastic differential equations.
The expansion of iterated Ito stochastic integrals of 
arbitrary mul\-tip\-li\-ci\-ty
$k$ $(k\in\mathbb{N})$ and expansions of iterated Stratonovich 
stochastic integrals of multipli\-ci\-ti\-es 1 to 6 have been obtained.
Considerable attention is paid to expansions based 
on multiple Fourier--Legendre series.
The exact and approximate expressions for the mean-square error of 
approximation of iterated Ito stochastic integrals are derived. 
The results of the article will be useful for
numerical integration of Ito stochastic differential equations
with non-commutative noise.

\medskip
\end{quote}
}


\setlength{\baselineskip}{2.0em}

\tableofcontents

\setlength{\baselineskip}{1.2em}


\vspace{5mm}

\section{Introduction}

\vspace{5mm}
 
Let $(\Omega,$ ${\rm F},$ ${\sf P})$ be a complete probability space, let 
$\{{\rm F}_t, t\in[0,T]\}$ be a nondecreasing right-continous family 
of $\sigma$-algebras of ${\rm F},$
and let ${\bf f}_t$ be a standard $m$-dimensional Wiener stochastic 
process, which is
${\rm F}_t$-measurable for any $t\in[0, T].$ We assume that the components
${\bf f}_{t}^{(i)}$ $(i=1,\ldots,m)$ of this process are independent. 
Consider
an Ito stochastic differential equation (SDE) in the integral form

\begin{equation}
\label{1.5.2}
{\bf x}_t={\bf x}_0+\int\limits_0^t {\bf a}({\bf x}_{\tau},\tau)d\tau+
\int\limits_0^t B({\bf x}_{\tau},\tau)d{\bf f}_{\tau},\ \ \
{\bf x}_0={\bf x}(0,\omega).
\end{equation}

\vspace{3mm}
\noindent
Here ${\bf x}_t$ is some $n$-dimensional stochastic process 
satisfying to (\ref{1.5.2}). 
The nonrandom functions ${\bf a}: \mathbb{R}^n\times[0, T]\to\mathbb{R}^n$,
$B: \mathbb{R}^n\times[0, T]\to\mathbb{R}^{n\times m}$
guarantee the existence and uniqueness up to stochastic 
equivalence of a solution
of the equation (\ref{1.5.2}) \cite{1}. The second integral on 
the right-hand side of (\ref{1.5.2}) is 
interpreted as an Ito stochastic integral.
Let ${\bf x}_0$ be an $n$-dimensional random variable, which is 
${\rm F}_0$-measurable and 
${\sf M}\{\left|{\bf x}_0\right|^2\}<\infty$ 
(${\sf M}$ denotes a mathematical expectation).
We assume that
${\bf x}_0$ and ${\bf f}_t-{\bf f}_0$ are independent when $t>0.$

It is well known \cite{KlPl2}-\cite{Mi3}
that Ito SDEs are 
adequate mathematical models of dynamic systems of 
different physical nature that are affected by random
perturbations.
For example, Ito SDEs are used as 
mathematical models in stochastic mathematical finance,
hydrology, seismology, geophysics, chemical kinetics, 
population dynamics,
electrodynamics, medicine and other fields 
\cite{KlPl2}-\cite{Mi3}.
Also these equations arise in optimal stochastic control, signal
filtering against the background of random noises, parameter estimation 
for stochastic systems as
well as in stochastic stability and bifurcations analysis 
\cite{KlPl2}, \cite{KPS}.

One of the effective approaches 
to the numerical integration of 
Ito SDEs is an approach based on the Taylor--Ito and 
Taylor--Stratonovich expansions
\cite{KlPl2}-\cite{20a-new-x}. The most important feature of such 
expan\-si\-ons is a presence in them of the so-called iterated
Ito and Stratonovich stochastic integrals, which play the key 
role for solving the 
problem of numerical integration of Ito SDEs and have the 
following form

\vspace{-3mm}
\begin{equation}
\label{ito}
J[\psi^{(k)}]_{T,t}=\int\limits_t^T\psi_k(t_k) \ldots \int\limits_t^{t_{2}}
\psi_1(t_1) d{\bf w}_{t_1}^{(i_1)}\ldots
d{\bf w}_{t_k}^{(i_k)},
\end{equation}

\vspace{1mm}
\begin{equation}
\label{str}
J^{*}[\psi^{(k)}]_{T,t}=
{\int\limits_t^{*}}^T \psi_k(t_k) \ldots {\int\limits_t^{*}}^{t_2}
\psi_1(t_1) d{\bf w}_{t_1}^{(i_1)}\ldots
d{\bf w}_{t_k}^{(i_k)},
\end{equation}

\vspace{4mm}
\noindent
where every $\psi_l(\tau)$ $(l=1,\ldots,k)$ is a nonrandom
function 
on $[t,T],$ ${\bf w}_{\tau}^{(i)}={\bf f}_{\tau}^{(i)}$
for $i=1,\ldots,m$ and
${\bf w}_{\tau}^{(0)}=\tau,$\ \
$i_1,\ldots,i_k = 0, 1,\ldots,m,$

\vspace{-1mm}
$$
\int\limits\ \ \hbox{and}\ \  \int\limits^{*}
$$ 

\vspace{3mm}
\noindent
denote Ito and 
Stratonovich stochastic integrals,
respectively (in this paper, 
we use the definition of the Stratonovich stochastic integral from \cite{KlPl2}).

Note that $\psi_l(\tau)\equiv 1$ $(l=1,\ldots,k)$ and
$i_1,\ldots,i_k = 0, 1,\ldots,m$ in  
\cite{KlPl2}-\cite{Mi3}. At the same time 
$\psi_l(\tau)\equiv (t-\tau)^{q_l}$ ($l=1,\ldots,k$; 
$q_1,\ldots,q_k=0, 1, 2,\ldots $) and $i_1,\ldots,i_k = 1,\ldots,m$ in
\cite{4}-\cite{20a-new-x}.

Effective solution 
of the problem of
combined mean-square approximation for collections 
of iterated Ito and Stratonovich stochastic integrals
(\ref{ito}) and (\ref{str})
composes the subject of this article (also see author's
publications \cite{7}-\cite{new-art-1-xxy}).

We want to mention in short that there are 
two main criteria of the numerical methods convergence 
for Ito SDEs \cite{KlPl2}-\cite{Mi3}:  
a strong or mean-square
criterion and a 
weak criterion, where the subject of approximation is not the solution 
of Ito SDE, simply stated, but the 
distribution of Ito SDE solution.

Using the strong numerical methods, we can build 
sample pathes
of Ito SDEs numerically. 
That is why
strong numerical methods are using when constructing new mathematical 
models on the basis of Ito SDEs. Moreover, these methods
are the tool for the numerical solution of different 
mathematical problems connected with Ito SDEs (see above)
\cite{KlPl2}-\cite{Mi3}. 

Strong numerical methods for Ito SDEs
require the combined mean-square approximation of collections 
of iterated Ito and Stratonovich stochastic integrals
(\ref{ito}) and (\ref{str}). 
The problem of effective jointly numerical modeling 
(with respect to the mean-square criterion of convergence) of iterated
Ito and Stratonovich stochastic integrals 
(\ref{ito}) and (\ref{str}) is 
complex from both
theoretical and computational points of view \cite{KlPl2}-\cite{rr}.

The only exception is connected with a narrow particular case, when 
$i_1=\ldots=i_k\ne 0$ and
$\psi_1(\tau),\ldots,\psi_k(\tau)\equiv \psi(\tau)$.
This case allows 
the investigation with using of the Ito formula 
\cite{KlPl2}-\cite{Mi3}.

Note that even for the mentioned coincidence ($i_1=\ldots=i_k\ne 0$), 
but for different 
functions $\psi_1(\tau),\ldots,\psi_k(\tau)$ the mentioned 
difficulties persist. As a result, 
relatively simple families of 
iterated Ito and Stratonovich stochastic integrals, 
which can be often 
met in the applications, cannot be expressed effectively in a finite 
form (with respect to the mean-square criterion
of approximation) 
using the system of standard 
Gaussian random variables. 

Why the problem of the mean-square approximation of iterated 
Ito and Stratonovich stochastic 
integrals is so complex?

Firstly, the mentioned stochastic integrals (in the case of fixed limits 
of integration) are the random variables, whose density functions are 
unknown in the general case. 
The exception is connected with
the narrow particular case which is
the simplest iterated Ito stochastic integral (\ref{ito}) 
with multiplicity 2
and $\psi_1(\tau), \psi_2(\tau)\equiv 1,$\
$i_1, i_2=1,\ldots,m$.
Nevertheless, the knowledge of this density function
not gives a simple way for approximation of iterated
Ito stochastic integral (\ref{ito}) of multiplicity 2 \cite{Wik}.

Secondly, we need to approximate not only one stochastic integral, 
but several iterated stochastic integrals which are complexly dependent
in a probabilistic meaning.

Often, the problem of combined mean-square approximation of iterated 
Ito and Stratonovich stochastic integrals occurs even in cases when the 
exact solution of Ito SDE is known. It means
that even if we know the solution of Ito SDE, 
we cannot model it numerically without the combined
numerical modeling of 
iterated Ito and Stratonovich 
stochastic integrals.

Note that for a number of special types of Ito SDEs
the problem of the mean-square approximation of iterated 
Ito and Stratonovich 
stochastic integrals can be simplified but cannot be solved. The equations
with additive vector noise, with scalar additive noise, 
with scalar non-additive  
noise, with a small parameter are related to such 
types of equations \cite{KlPl2}-\cite{Mi3}.

For the mentioned types of equations, 
simplifications are connected to the fact 
that some members
from stochastic Taylor expansions are equal to zero 
or we may neglect some members from these expansions
due to the presence of a small 
parameter \cite{KlPl2}-\cite{Mi3}.

Seems that iterated stochastic integrals may be approximated by multiple 
integral sums of different types \cite{Mi2}, \cite{Mi3}, \cite{Al}. 
However, this approach implies the partitioning of the interval 
of integration $[t, T]$ for iterated stochastic integrals. The length
$T-t$ of this interval is already fairly
small (because it is a step of integration of numerical methods for 
Ito SDEs) and does not need to be partitioned.
Computational experiments show that the application
of numerical simulation for iterated stochastic integrals 
(in which the interval of integration $[t, T]$ is
partitioned) leads to unacceptably high com\-pu\-ta\-ti\-o\-nal cost and 
accumulation of computation errors \cite{7}.

In \cite{Mi2} (also see \cite{KlPl2}, \cite{KPS}, \cite{Mi3}, 
\cite{KPW}, \cite{Zapad-9}), 
Milstein G.N. proposed to expand (\ref{ito}) or (\ref{str})
into the iterated series of products
of standard Gaussian random variables by representing the Wiener
process as a trigonometric Fourier series with random coefficients 
(version of the so-called Karhunen--Loeve expansion of the
Brownian bridge process).
To obtain the Milstein expansion of (\ref{str}), the truncated Fourier
expansions of components of the multidimensional
Wiener process ${\bf f}_s$ must be
iteratively substituted in the single integrals, and the integrals
must be calculated, starting from the innermost integral.
This is a complicated procedure that does not lead to a general
expansion of (\ref{str}) valid for an arbitrary multiplicity $k.$
For this reason, only expansions of single, double and triple
stochastic
integrals (\ref{str}) were obtained
\cite{KlPl2}, \cite{KPS}, \cite{KPW}, \cite{Zapad-9} ($k=1, 2, 3$),
\cite{Mi2}, \cite{Mi3} ($k=1, 2$) 
for the case $\psi_1(\tau), \psi_2(\tau),$ $\psi_3(\tau)$ $\equiv 1;$ 
$i_1, i_2, i_3=0,1,\ldots,m.$

It should be noted that 
the authors of the publications \cite{KlPl2}
(Sect.~5.8, pp.~202--204), \cite{KPS} (pp.~82-84),
\cite{KPW} (pp.~438-439),  
\cite{Zapad-9} (pp.~263-264) 
use the Wong--Zakai approximation 
\cite{W-Z-1}-\cite{Watanabe} (without rigorous proof)
within the frames of the mentioned approach \cite{Mi2}
based on the approximation of the Wiener
process in the form of its series expansion
(see discussion in Sect.~7 for details).

Note that in \cite{rra}, \cite{rr} the truncated expansions of the Wiener 
processes based on the Haar functions \cite{rr}
and trigonometric functions \cite{rra}, \cite{rr} were
applied for the expansion of 
double \cite{rra}, \cite{rr} and triple \cite{rra}
Ito stochastic integrals (\ref{ito}).
The expansions from \cite{rra}, \cite{rr} also lead to 
iterated application of the
operation of limit transition as in the Milstein approach \cite{Mi2}.

It is necessary to note that the Milstein approach \cite{Mi2} 
excelled at least
in several times (or even in several orders) 
the methods of multiple integral sums \cite{Mi2}, \cite{Mi3}, \cite{Al}
considering computational costs in the sense 
of their diminishing \cite{Mi2}, \cite{Mi3}, \cite{7}.

An alternative and more general strong approximation method was 
proposed for (\ref{str}) in \cite{20a} (Sect.~2.4) (also see \cite{4} (1998),
\cite{11}-\cite{16}, \cite{19}, \cite{20}, \cite{20a-new}, \cite{20a-new-x}, \cite{3} (1997)). 
In these papers
$J^{*}[\psi^{(k)}]_{T,t}$ was represented as the multiple stochastic 
integral
of the certain discontinuous nonrandom function of $k$ ($k\in\mathbb{N}$)
variables, and the 
function
was then expanded into the generalized iterated Fourier series by
complete system of 
continuously
differentiable functions that are orthonormal in the space
$L_2([t, T])$. As a result,
the general iterated series expansion of products
of standard Gaussian random variables 
was obtained in 
\cite{20a} (Sect.~2.4) (also see \cite{4} (1998),
\cite{11}-\cite{16}, \cite{19}, \cite{20}, \cite{20a-new}, \cite{20a-new-x}, \cite{3} (1997))
for the iterated Stratonovich
stochastic integrals (\ref{str})
of arbitrary multiplicity $k$ ($k\in\mathbb{N}$).
Hereinafter, this method is referred to as the method of 
generalized iterated Fourier series.

Consider the formulation of
the method of 
generalized iterated Fourier series.
Let us introduce the following 
function $K(t_1,\ldots,t_k)$ defined on the $k$-dimensional
hypercube $[t, T]^k$

\begin{equation}
\label{ppp}
K(t_1,\ldots,t_k)=
\begin{cases}
\psi_1(t_1)\ldots \psi_k(t_k),\ &t_1<\ldots<t_k\\
~\\
~\\
0,\  &\hbox{\rm otherwise}
\end{cases}\ \ \
=\ \ \
\prod\limits_{l=1}^k
\psi_l(t_l)\ \prod\limits_{l=1}^{k-1}{\bf 1}_{\{t_l<t_{l+1}\}},\ 
\end{equation}

\vspace{5mm}
\noindent
where $t_1,\ldots,t_k\in [t, T]$ $(k\ge 2)$ and 
$K(t_1)\equiv\psi_1(t_1)$ for $t_1\in[t, T].$ Here 
$\psi_1(\tau),\ldots,\psi_k(\tau)\in L_2([t, T])$
and
${\bf 1}_A$ denotes the indicator of the set $A$.

\vspace{2mm}

{\bf Theorem 1}\ \cite{20a} (Sect.~2.4) (also see \cite{4} (1998),
\cite{11}-\cite{16}, \cite{19}, \cite{20}, \cite{20a-new}, \cite{20a-new-x}, \cite{3} (1997)).\
{\it Suppose that every $\psi_l(\tau)$ $(l=1,\ldots,k)$ is twice continuously
differentiable function at the interval
$[t, T]$ and
$\{\phi_j(x)\}_{j=0}^{\infty}$ is a complete
orthonormal system of 
trigonometric functions in the space $L_2([t, T])$. 
Then, the iterated Stratonovich stochastic integral 
{\rm(\ref{str})}
is expanded into the 
conver\-ging 
in the mean of degree $2n$ $(n\in \mathbb{N})$
iterated series

\vspace{-3mm}
\begin{equation}
\label{1500}
J^{*}[\psi^{(k)}]_{T,t}=
\sum_{j_1=0}^{\infty}\ldots\sum_{j_k=0}^{\infty}
C_{j_k\ldots j_1}
\prod_{l=1}^k
\zeta^{(i_l)}_{j_l},
\end{equation}

\vspace{2mm}
\noindent
which means the following

\vspace{-2mm}
\begin{equation}
\label{1500e}
\lim\limits_{p_1\to\infty}\varlimsup\limits_{p_2\to\infty}
\ldots\varlimsup\limits_{p_k\to\infty}
{\sf M}\left\{\left(J^{*}[\psi^{(k)}]_{T,t}-
\sum_{j_1=0}^{p_1}\ldots\sum_{j_k=0}^{p_k}
C_{j_k\ldots j_1}
\prod_{l=1}^k
\zeta^{(i_l)}_{j_l}\right)^{2n}\right\}=0,
\end{equation}

\vspace{2mm}
\noindent
where $\varlimsup$ means $\limsup,$
\begin{equation}
\label{z1z}
\zeta_{j}^{(i)}=
\int\limits_t^T \phi_{j}(\tau) d{\bf w}_{\tau}^{(i)}
\end{equation} 

\vspace{2mm}
\noindent
are independent standard Gaussian random variables
for 
various
$i$ or $j$ {\rm(}if $i\ne 0${\rm)} and 

\vspace{-2mm}
\begin{equation}
\label{333.40}
C_{j_k\ldots j_1}=\int\limits_{[t,T]^k}
K(t_1,\ldots,t_k)\prod_{l=1}^{k}\phi_{j_l}(t_l)dt_1\ldots dt_k
\end{equation}

\vspace{2mm}
\noindent
is the Fourier coefficient.}

\vspace{2mm}

Note that an anlogue of Theorem~1 
for the 
case of Legendre 
polynomials, $n=1$ (the case of mean-square convergence), and $k=2$ is obtained in 
\cite{20a} (Sect.~2.4.1), \cite{23} (Sect.~2).

The proof of Theorem 1 is based on the following
statement.

\vspace{2mm}

{\bf Lemma 1}\ \cite{3} (also see \cite{4},
\cite{11}-\cite{16}, \cite{19}, \cite{20}-\cite{20a-new-x}).\
{\it Under the conditions of Theorem {\rm 1} the function
$K^{*}(t_1,\ldots,t_k)$
is represented in any internal point of the hypercube   
$[t,T]^k$ by the generalized iterated Fourier series

\vspace{-2mm}
\begin{equation}
\label{30.18}
K^{*}(t_1,\ldots,t_k)=
\sum_{j_1=0}^{\infty}\ldots \sum_{j_k=0}^{\infty}
C_{j_k\ldots j_1}\prod_{l=1}^{k} \phi_{j_l}(t_l),\ \ \ 
(t_1,\ldots,t_k)\in (t, T)^k,
\end{equation}

\vspace{2mm}
\noindent
where 
$$
K^{*}(t_1,\ldots,t_k)=\prod\limits_{l=1}^k\psi_l(t_l)
\prod_{l=1}^{k-1}\Biggl({\bf 1}_{\{t_l<t_{l+1}\}}+
\frac{1}{2}{\bf 1}_{\{t_l=t_{l+1}\}}\Biggr)=
$$

\vspace{3mm}
$$
=\prod_{l=1}^k \psi_l(t_l)\left(\prod_{l=1}^{k-1}
{\bf 1}_{\{t_l<t_{l+1}\}}+
\sum_{r=1}^{k-1}\frac{1}{2^r}
\sum_{\stackrel{s_r,\ldots,s_1=1}{{}_{s_r>\ldots>s_1}}}^{k-1}\ 
\prod_{l=1}^r {\bf 1}_{\{t_{s_l}=t_{s_l+1}\}}
\prod_{\stackrel{l=1}{{}_{l\ne s_1,\ldots, s_r}}}^{k-1}
{\bf 1}_{\{t_{l}<t_{l+1}\}}\right)
$$

\vspace{5mm}
\noindent
for $t_1,\ldots,t_k\in[t, T]$\ $(k\ge 2)$ and 
$K^{*}(t_1)\equiv\psi_1(t_1)$ for $t_1\in[t, T],$ 
${\bf 1}_A$ is the indicator of the set $A,$ the Fourier coefficient
$C_{j_k\ldots j_1}$ has the form {\rm (\ref{333.40})}.
At that{\rm ,} the iterated series {\rm (\ref{30.18})} converges at the 
boundary of the hypercube $[t,T]^k$
{\rm (}not necessarily to the function $K^{*}(t_1,\ldots,t_k)${\rm )}.}

\vspace{2mm}

In \cite{4}, \cite{11}-\cite{16}, \cite{19}, \cite{20}-\cite{20a-new-x}, 
\cite{3} it was shown that the method of 
generalized iterated Fourier series leads for $k=2$ and $\psi_1(\tau),
\psi_2(\tau)\equiv 1$ (the case of trigonometric system 
of functions) to the
Milstein expansion of (\ref{str}) \cite{Mi2}.

As we noted above, the method of generalized iterated Fourier series 
as well as the method from \cite{Mi2}
lead to iterated application of the operation of 
limit transition. So, the convergence problem of the following
approximation

\vspace{-2mm}
\begin{equation}
\label{pp}
J^{*}[\psi^{(k)}]_{T,t}^{p}=
\sum_{j_1,\ldots,j_k=0}^{p}
C_{j_k\ldots j_1}
\prod_{l=1}^k
\zeta^{(i_l)}_{j_l}
\end{equation}

\vspace{3mm}
\noindent
to $J^{*}[\psi^{(k)}]_{T,t}$ if $p\to\infty$ 
in the mean-square sense must be considered separately
(see Sect.~2 and discussion in Sect.~7 for details).
The mentioned problem appears for triple stochastic integrals
or even for some double stochastic integrals 
in the case, when $\psi_1(\tau),$ $\psi_2(\tau)\not\equiv 1$
(see above).

\vspace{5mm}

\section{Method of the Mean-Square Approximation
of Iterated Ito 
Stochastic Integrals Based on Generalized Multiple Fourier Series}

\vspace{5mm}

In the previous section we paid attention on the fact
that the method from \cite{Mi2} 
and the method of generalized 
iterated Fourier series 
\cite{4}, \cite{11}-\cite{16}, \cite{19}, \cite{20}-\cite{20a-new-x},
\cite{3} lead
to iterated application of the operation of limit transition.
So these methods may not converge in 
the mean-square sense 
to the appropriate iterated stochastic integrals (\ref{str}) 
for some methods of series summation (see (\ref{pp})).

The difficulties noted above can be overcome by the another method.
The idea of this method is as follows: the iterated Ito stochastic 
integral (\ref{ito})
of multiplicity $k$ ($k\in\mathbb{N}$) 
is represented as the multiple stochastic 
integral from the nonrandom discontinuous function 
$K(t_1,\ldots,t_k)$ defined 
on the hypercube $[t, T]^k$ by the relation (\ref{ppp}), 
where $[t, T]$ is the interval of 
integration of the iterated Ito stochastic integral. Then, 
the function $K(t_1,\ldots,t_k)$ 
is expanded in the hypercube $[t, T]^k$ into the generalized 
multiple Fourier series converging 
in the sense of norm in Hilbert space
$L_2([t,T]^k)$. After a number of nontrivial transformations we come 
(see Theorems 2--4 below) to the 
mean-square convergening expansion of the 
iterated Ito stochastic integral (\ref{ito})
into the multiple 
series of products
of standard  Gaussian random 
variables. The coefficients of this 
series are the coefficients of 
the generalized multiple Fourier series for the function $K(t_1,\ldots,t_k)$, 
which can be calculated using the explicit formula 
regardless of multiplicity $k$ of the iterated Ito stochastic integral.
Hereinafter, this method is referred to as the method of 
generalized multiple Fourier series.

Suppose that $\{\phi_j(x)\}_{j=0}^{\infty}$
is a complete orthonormal system of functions in 
the space $L_2([t, T])$. 
The function $K(t_1,\ldots,t_k)$ (defined by (\ref{ppp}))
belongs to the space $L_2([t, T]^k).$
At this situation it is well known that the generalized 
multiple Fourier series 
of $K(t_1,\ldots,t_k)\in L_2([t, T]^k)$ is converging 
to $K(t_1,\ldots,t_k)$ in the hypercube $[t, T]^k$ in 
the mean-square sense, i.e.

\begin{equation}
\label{sos1z}
\lim\limits_{p_1,\ldots,p_k\to \infty}
\Biggl\Vert
K(t_1,\ldots,t_k)-
\sum_{j_1=0}^{p_1}\ldots \sum_{j_k=0}^{p_k}
C_{j_k\ldots j_1}\prod_{l=1}^{k} \phi_{j_l}(t_l)\Biggr\Vert_{L_2([t, T]^k)}
=0,
\end{equation}

\vspace{4mm}
\noindent
where
\begin{equation}
\label{ppppa}
C_{j_k\ldots j_1}=\int\limits_{[t,T]^k}
K(t_1,\ldots,t_k)\prod_{l=1}^{k}\phi_{j_l}(t_l)dt_1\ldots dt_k
\end{equation}

\vspace{4mm}
\noindent
is the Fourier coefficient, and

$$
\left\Vert f\right\Vert_{L_2([t, T]^k)}=\left(\int\limits_{[t,T]^k}
f^2(t_1,\ldots,t_k)dt_1\ldots dt_k\right)^{1/2}.
$$

\vspace{6mm}

Consider the partition $\{\tau_j\}_{j=0}^N$ of $[t,T]$ such that

\begin{equation}
\label{1111}
t=\tau_0<\ldots <\tau_N=T,\ \ \
\Delta_N=
\hbox{\vtop{\offinterlineskip\halign{
\hfil#\hfil\cr
{\rm max}\cr
$\stackrel{}{{}_{0\le j\le N-1}}$\cr
}} }\Delta\tau_j\to 0\ \ \hbox{if}\ \ N\to \infty,\ \ \ 
\Delta\tau_j=\tau_{j+1}-\tau_j.
\end{equation}

\vspace{5mm}

{\bf Theorem 2} \cite{7} (2006), \cite{8}, \cite{19}, \cite{20}-\cite{31aaa}.\ 
{\it Suppose that
every $\psi_l(\tau)$ $(l=1,\ldots, k)$ is a continuous function on 
$[t, T]$ and
$\{\phi_j(x)\}_{j=0}^{\infty}$ is a complete orthonormal system  
of continuous functions in the space $L_2([t,T]).$ 
Then

\vspace{-1mm}
$$
J[\psi^{(k)}]_{T,t}\  =\ 
\hbox{\vtop{\offinterlineskip\halign{
\hfil#\hfil\cr
{\rm l.i.m.}\cr
$\stackrel{}{{}_{p_1,\ldots,p_k\to \infty}}$\cr
}} }\sum_{j_1=0}^{p_1}\ldots\sum_{j_k=0}^{p_k}
C_{j_k\ldots j_1}\Biggl(
\prod_{l=1}^k\zeta_{j_l}^{(i_l)}\ -
\Biggr.
$$

\vspace{2mm}
\begin{equation}
\label{tyyy}
-\ \ \Biggl.
\hbox{\vtop{\offinterlineskip\halign{
\hfil#\hfil\cr
{\rm l.i.m.}\cr
$\stackrel{}{{}_{N\to \infty}}$\cr
}} }\sum_{(l_1,\ldots,l_k)\in {\rm G}_k}
\phi_{j_{1}}(\tau_{l_1})
\Delta{\bf w}_{\tau_{l_1}}^{(i_1)}\ldots
\phi_{j_{k}}(\tau_{l_k})
\Delta{\bf w}_{\tau_{l_k}}^{(i_k)}\Biggr),
\end{equation}

\vspace{5mm}
\noindent
where

$$
{\rm G}_k={\rm H}_k\backslash{\rm L}_k,\ \ \
{\rm H}_k=\{(l_1,\ldots,l_k):\ l_1,\ldots,l_k=0,\ 1,\ldots,N-1\},
$$

$$
{\rm L}_k=\{(l_1,\ldots,l_k):\ l_1,\ldots,l_k=0,\ 1,\ldots,N-1;\
l_g\ne l_r\ (g\ne r);\ g, r=1,\ldots,k\},
$$

\vspace{4mm}
\noindent
${\rm l.i.m.}$ is a limit in the mean-square sense,
$i_1,\ldots,i_k=0,1,\ldots,m,$ 

\begin{equation}
\label{rr23}
\zeta_{j}^{(i)}=
\int\limits_t^T \phi_{j}(\tau) d{\bf w}_{\tau}^{(i)}
\end{equation} 

\vspace{3mm}
\noindent
are independent standard Gaussian random variables
for various
$i$ or $j$ {\rm(}if $i\ne 0${\rm),}
$C_{j_k\ldots j_1}$ is the Fourier coefficient {\rm(\ref{ppppa}),}
$\Delta{\bf w}_{\tau_{j}}^{(i)}=
{\bf w}_{\tau_{j+1}}^{(i)}-{\bf w}_{\tau_{j}}^{(i)}$
$(i=0,\ 1,\ldots,m),$\
$\left\{\tau_{j}\right\}_{j=0}^{N}$ is a partition of
$[t,T],$ which satisfies the condition {\rm (\ref{1111})}.}

\vspace{3mm}

{\bf Proof.}\ The proof
of Theorem 2 is based on Lemmas 5.1--5.3 \cite{20}
(P. A.253--A.259), Lemmas 1.1--1.3 \cite{20a}-\cite{20a-new-x} or
Lemmas 1--3 \cite{26a}.
According to Lemma 5.1 \cite{20}, Lemma 1.1 \cite{20a}-\cite{20a-new-x} or
Lemma 1 \cite{26a},
we have 

\vspace{1mm}
$$
J[\psi^{(k)}]_{T,t}=
\hbox{\vtop{\offinterlineskip\halign{
\hfil#\hfil\cr
{\rm l.i.m.}\cr
$\stackrel{}{{}_{N\to \infty}}$\cr
}} }\sum_{l_k=0}^{N-1}\ldots\sum_{l_1=0}^{l_2-1}
\psi_1(\tau_{l_1})\ldots\psi_k(\tau_{l_k})
\Delta{\bf w}_{\tau_{l_1}}^{(i_1)}
\ldots
\Delta{\bf w}_{\tau_{l_k}}^{(i_k)}=
$$

\vspace{3mm}
$$
=\hbox{\vtop{\offinterlineskip\halign{
\hfil#\hfil\cr
{\rm l.i.m.}\cr
$\stackrel{}{{}_{N\to \infty}}$\cr
}} }\sum_{l_k=0}^{N-1}\ldots\sum_{l_1=0}^{l_2-1}
K(\tau_{l_1},\ldots,\tau_{l_k})
\Delta{\bf w}_{\tau_{l_1}}^{(i_1)}
\ldots
\Delta{\bf w}_{\tau_{l_k}}^{(i_k)}=
$$

\vspace{3mm}
$$
=
\hbox{\vtop{\offinterlineskip\halign{
\hfil#\hfil\cr
{\rm l.i.m.}\cr
$\stackrel{}{{}_{N\to \infty}}$\cr
}} }\sum_{l_k=0}^{N-1}\ldots\sum_{l_1=0}^{N-1}
K(\tau_{l_1},\ldots,\tau_{l_k})
\Delta{\bf w}_{\tau_{l_1}}^{(i_1)}
\ldots
\Delta{\bf w}_{\tau_{l_k}}^{(i_k)}=
$$

\vspace{3mm}
$$
=\hbox{\vtop{\offinterlineskip\halign{
\hfil#\hfil\cr
{\rm l.i.m.}\cr
$\stackrel{}{{}_{N\to \infty}}$\cr
}} }
\sum\limits_{\stackrel{l_1,\ldots,l_k=0}{{}_{l_q\ne l_r;\ q\ne r;\ q, r=1,\ldots, k}}}^{N-1}
K(\tau_{l_1},\ldots,\tau_{l_k})
\Delta{\bf w}_{\tau_{l_1}}^{(i_1)}
\ldots
\Delta{\bf w}_{\tau_{l_k}}^{(i_k)}=
$$

\vspace{3mm}
\begin{equation}
\label{hehe100}
=
\int\limits_{t}^{T}
\ldots
\int\limits_{t}^{t_2}
\sum_{(t_1,\ldots,t_k)}\left(
K(t_1,\ldots,t_k)d{\bf w}_{t_1}^{(i_1)}
\ldots
d{\bf w}_{t_k}^{(i_k)}\right),
\end{equation}

\vspace{6mm}
\noindent
where permutations 
$(t_1,\ldots,t_k)$ when summing
are performed only 
in the expression, which is enclosed in pa\-ren\-the\-ses.

It is easy to see that (\ref{hehe100})
can be written in the form \cite{20}-\cite{20a-new-x}, \cite{26a}

$$
J[\psi^{(k)}]_{T,t}=\sum_{(t_1,\ldots,t_k)}
\int\limits_{t}^{T}
\ldots
\int\limits_{t}^{t_2}
K(t_1,\ldots,t_k)d{\bf w}_{t_1}^{(i_1)}
\ldots
d{\bf w}_{t_k}^{(i_k)},
$$

\vspace{3mm}
\noindent
where permutations
$(t_1,\ldots,t_k)$ when summing
are performed only in 
the values
$d{\bf w}_{t_1}^{(i_1)}
\ldots $
$d{\bf w}_{t_k}^{(i_k)}$. At the same time the indexes near upper 
limits of integration in the iterated stochastic integrals are changed 
correspondently and if $t_r$ swapped with $t_q$ in the  
permutation $(t_1,\ldots,t_k)$, then $i_r$ swapped with $i_q$ in 
the permutation $(i_1,\ldots,i_k)$.

Since the integration of bounded function with respect to the 
set of measure zero for Riemann or Lebesgue integrals gives zero result, then the 
following formula is correct for these integrals

\begin{equation}
\label{zx1}
\int\limits_{[t, T]^k}|G(t_1,\ldots,t_k)|dt_1\ldots dt_k=
\sum_{(t_1,\ldots,t_k)}
\int\limits_{t}^{T}
\ldots
\int\limits_{t}^{t_2}
|G(t_1,\ldots,t_k)|dt_1\ldots dt_k,
\end{equation}

\vspace{3mm}
\noindent
where permutations $(t_1,\ldots,t_k)$ when summing
are performed only 
in the 
va\-lues $dt_1,\ldots, dt_k$. At the same time the indexes near upper 
limits of integration are changed correspondently
and the function $|G(t_1,\ldots,t_k)|$ is supposed as integrated in 
the hypercube $[t, T]^k.$

According to 
Lemmas 5.1--5.3 \cite{20}
(P. A.253--A.259), Lemmas 1.1--1.3 \cite{20a}-\cite{20a-new-x} or
Lemmas 1--3 \cite{26a},
we get the following representation

\vspace{2mm}
$$
J[\psi^{(k)}]_{T,t}=
\sum_{j_1=0}^{p_1}\ldots
\sum_{j_k=0}^{p_k}
C_{j_k\ldots j_1}
\int\limits_{t}^{T}
\ldots
\int\limits_{t}^{t_2}
\sum_{(t_1,\ldots,t_k)}\left(
\phi_{j_1}(t_1)
\ldots
\phi_{j_k}(t_k)
d{\bf w}_{t_1}^{(i_1)}
\ldots
d{\bf w}_{t_k}^{(i_k)}\right)+
$$

\vspace{3mm}
$$
+R_{T,t}^{p_1,\ldots,p_k}=
$$

\vspace{5mm}
$$
=\sum_{j_1=0}^{p_1}\ldots
\sum_{j_k=0}^{p_k}
C_{j_k\ldots j_1}             
\hbox{\vtop{\offinterlineskip\halign{
\hfil#\hfil\cr
{\rm l.i.m.}\cr
$\stackrel{}{{}_{N\to \infty}}$\cr
}} }
\sum\limits_{\stackrel{l_1,\ldots,l_k=0}{{}_{l_q\ne l_r;\ 
q\ne r;\ q, r=1,\ldots, k}}}^{N-1}
\phi_{j_1}(\tau_{l_1})\ldots
\phi_{j_k}(\tau_{l_k})
\Delta{\bf w}_{\tau_{l_1}}^{(i_1)}
\ldots
\Delta{\bf w}_{\tau_{l_k}}^{(i_k)}+
$$

\vspace{3mm}
\begin{equation}
\label{e1}
+R_{T,t}^{p_1,\ldots,p_k}=
\end{equation}

\vspace{5mm}
$$
=\sum_{j_1=0}^{p_1}\ldots
\sum_{j_k=0}^{p_k}
C_{j_k\ldots j_1}\Biggl(
\hbox{\vtop{\offinterlineskip\halign{
\hfil#\hfil\cr
{\rm l.i.m.}\cr
$\stackrel{}{{}_{N\to \infty}}$\cr
}} }\sum_{l_1,\ldots,l_k=0}^{N-1}
\phi_{j_1}(\tau_{l_1})
\ldots
\phi_{j_k}(\tau_{l_k})
\Delta{\bf w}_{\tau_{l_1}}^{(i_1)}
\ldots
\Delta{\bf w}_{\tau_{l_k}}^{(i_k)}
-\Biggr.
$$

\vspace{4mm}
$$
-\Biggl.
\hbox{\vtop{\offinterlineskip\halign{
\hfil#\hfil\cr
{\rm l.i.m.}\cr
$\stackrel{}{{}_{N\to \infty}}$\cr
}} }\sum_{(l_1,\ldots,l_k)\in {\rm G}_k}
\phi_{j_{1}}(\tau_{l_1})
\Delta{\bf w}_{\tau_{l_1}}^{(i_1)}\ldots
\phi_{j_{k}}(\tau_{l_k})
\Delta{\bf w}_{\tau_{l_k}}^{(i_k)}\Biggr)+
$$

\vspace{3mm}
$$
+R_{T,t}^{p_1,\ldots,p_k}=
$$

\vspace{5mm}
$$
=\sum_{j_1=0}^{p_1}\ldots\sum_{j_k=0}^{p_k}
C_{j_k\ldots j_1}\left(
\prod_{l=1}^k\zeta_{j_l}^{(i_l)}-
\hbox{\vtop{\offinterlineskip\halign{
\hfil#\hfil\cr
{\rm l.i.m.}\cr
$\stackrel{}{{}_{N\to \infty}}$\cr
}} }\sum_{(l_1,\ldots,l_k)\in {\rm G}_k}
\phi_{j_{1}}(\tau_{l_1})
\Delta{\bf w}_{\tau_{l_1}}^{(i_1)}\ldots
\phi_{j_{k}}(\tau_{l_k})
\Delta{\bf w}_{\tau_{l_k}}^{(i_k)}\right)+
$$

\vspace{3mm}
\begin{equation}
\label{e2}
+R_{T,t}^{p_1,\ldots,p_k}\ \ \ \hbox{w.\ p.\ 1},
\end{equation}

\vspace{5mm}
\noindent
where

$$
R_{T,t}^{p_1,\ldots,p_k}=
$$

\begin{equation}
\label{y007}
=\sum_{(t_1,\ldots,t_k)}
\int\limits_{t}^{T}
\ldots
\int\limits_{t}^{t_2}
\left(K(t_1,\ldots,t_k)-
\sum_{j_1=0}^{p_1}\ldots
\sum_{j_k=0}^{p_k}
C_{j_k\ldots j_1}
\prod_{l=1}^k\phi_{j_l}(t_l)\right)
d{\bf w}_{t_1}^{(i_1)}
\ldots
d{\bf w}_{t_k}^{(i_k)},
\end{equation}

\vspace{5mm}
\noindent
where permutations $(t_1,\ldots,t_k)$ when summing are performed only 
in the values $d{\bf w}_{t_1}^{(i_1)}
\ldots $
$d{\bf w}_{t_k}^{(i_k)}$. At the same time the indexes near 
upper limits of integration in the iterated stochastic integrals 
are changed correspondently and if $t_r$ swapped with $t_q$ in the  
permutation $(t_1,\ldots,t_k)$, then $i_r$ swapped with $i_q$ in the 
permutation $(i_1,\ldots,i_k)$.

Let us estimate the remainder
$R_{T,t}^{p_1,\ldots,p_k}$ of the series.

By Lemma 5.2 \cite{20} (pp. A.257--A.258),
Lemma 1.2 \cite{20a}-\cite{20a-new-x}
or Lemma 2 \cite{26a} we have (see (\ref{zx1}))

\vspace{2mm}
$$
{\sf M}\left\{\left(R_{T,t}^{p_1,\ldots,p_k}\right)^2\right\}
\le 
$$

\vspace{2mm}
$$
\le C_k
\sum_{(t_1,\ldots,t_k)}
\int\limits_{t}^{T}
\ldots
\int\limits_{t}^{t_2}
\left(K(t_1,\ldots,t_k)-
\sum_{j_1=0}^{p_1}\ldots
\sum_{j_k=0}^{p_k}
C_{j_k\ldots j_1}
\prod_{l=1}^k\phi_{j_l}(t_l)\right)^2
dt_1
\ldots
dt_k=
$$

\vspace{2mm}
\begin{equation}
\label{obana1}
=C_k\int\limits_{[t,T]^k}
\left(K(t_1,\ldots,t_k)-
\sum_{j_1=0}^{p_1}\ldots
\sum_{j_k=0}^{p_k}
C_{j_k\ldots j_1}
\prod_{l=1}^k\phi_{j_l}(t_l)\right)^2
dt_1
\ldots
dt_k\to 0
\end{equation}

\vspace{5mm}
\noindent
if $p_1,\ldots,p_k\to\infty,$ where constant $C_k$ 
depends only
on the multiplicity $k$ of the iterated Ito stochastic integral. 
Theorem 2 is proved.

In order to evaluate the significance of Theorem 2 for practice we will
demonstrate its transformed particular cases for 
$k=1,\ldots,6$ \cite{7}-\cite{19}, \cite{20}-\cite{31aaa}

\begin{equation}
\label{a1}
J[\psi^{(1)}]_{T,t}
=\hbox{\vtop{\offinterlineskip\halign{
\hfil#\hfil\cr
{\rm l.i.m.}\cr
$\stackrel{}{{}_{p_1\to \infty}}$\cr
}} }\sum_{j_1=0}^{p_1}
C_{j_1}\zeta_{j_1}^{(i_1)},
\end{equation}

\vspace{2mm}
\begin{equation}
\label{a2}
J[\psi^{(2)}]_{T,t}
=\hbox{\vtop{\offinterlineskip\halign{
\hfil#\hfil\cr
{\rm l.i.m.}\cr
$\stackrel{}{{}_{p_1,p_2\to \infty}}$\cr
}} }\sum_{j_1=0}^{p_1}\sum_{j_2=0}^{p_2}
C_{j_2j_1}\Biggl(\zeta_{j_1}^{(i_1)}\zeta_{j_2}^{(i_2)}
-{\bf 1}_{\{i_1=i_2\ne 0\}}
{\bf 1}_{\{j_1=j_2\}}\Biggr),
\end{equation}

\vspace{5mm}
$$
J[\psi^{(3)}]_{T,t}=
\hbox{\vtop{\offinterlineskip\halign{
\hfil#\hfil\cr
{\rm l.i.m.}\cr
$\stackrel{}{{}_{p_1,\ldots,p_3\to \infty}}$\cr
}} }\sum_{j_1=0}^{p_1}\sum_{j_2=0}^{p_2}\sum_{j_3=0}^{p_3}
C_{j_3j_2j_1}\Biggl(
\zeta_{j_1}^{(i_1)}\zeta_{j_2}^{(i_2)}\zeta_{j_3}^{(i_3)}
-\Biggr.
$$
\begin{equation}
\label{leto5002}
\Biggl.-{\bf 1}_{\{i_1=i_2\ne 0\}}
{\bf 1}_{\{j_1=j_2\}}
\zeta_{j_3}^{(i_3)}
-{\bf 1}_{\{i_2=i_3\ne 0\}}
{\bf 1}_{\{j_2=j_3\}}
\zeta_{j_1}^{(i_1)}-
{\bf 1}_{\{i_1=i_3\ne 0\}}
{\bf 1}_{\{j_1=j_3\}}
\zeta_{j_2}^{(i_2)}\Biggr),
\end{equation}

\vspace{4mm}
$$
J[\psi^{(4)}]_{T,t}
=
\hbox{\vtop{\offinterlineskip\halign{
\hfil#\hfil\cr
{\rm l.i.m.}\cr
$\stackrel{}{{}_{p_1,\ldots,p_4\to \infty}}$\cr
}} }\sum_{j_1=0}^{p_1}\ldots\sum_{j_4=0}^{p_4}
C_{j_4\ldots j_1}\Biggl(
\prod_{l=1}^4\zeta_{j_l}^{(i_l)}
\Biggr.
-
$$
$$
-
{\bf 1}_{\{i_1=i_2\ne 0\}}
{\bf 1}_{\{j_1=j_2\}}
\zeta_{j_3}^{(i_3)}
\zeta_{j_4}^{(i_4)}
-
{\bf 1}_{\{i_1=i_3\ne 0\}}
{\bf 1}_{\{j_1=j_3\}}
\zeta_{j_2}^{(i_2)}
\zeta_{j_4}^{(i_4)}-
$$
$$
-
{\bf 1}_{\{i_1=i_4\ne 0\}}
{\bf 1}_{\{j_1=j_4\}}
\zeta_{j_2}^{(i_2)}
\zeta_{j_3}^{(i_3)}
-
{\bf 1}_{\{i_2=i_3\ne 0\}}
{\bf 1}_{\{j_2=j_3\}}
\zeta_{j_1}^{(i_1)}
\zeta_{j_4}^{(i_4)}-
$$
$$
-
{\bf 1}_{\{i_2=i_4\ne 0\}}
{\bf 1}_{\{j_2=j_4\}}
\zeta_{j_1}^{(i_1)}
\zeta_{j_3}^{(i_3)}
-
{\bf 1}_{\{i_3=i_4\ne 0\}}
{\bf 1}_{\{j_3=j_4\}}
\zeta_{j_1}^{(i_1)}
\zeta_{j_2}^{(i_2)}+
$$
$$
+
{\bf 1}_{\{i_1=i_2\ne 0\}}
{\bf 1}_{\{j_1=j_2\}}
{\bf 1}_{\{i_3=i_4\ne 0\}}
{\bf 1}_{\{j_3=j_4\}}
+
{\bf 1}_{\{i_1=i_3\ne 0\}}
{\bf 1}_{\{j_1=j_3\}}
{\bf 1}_{\{i_2=i_4\ne 0\}}
{\bf 1}_{\{j_2=j_4\}}+
$$
\begin{equation}
\label{leto5003}
+\Biggl.
{\bf 1}_{\{i_1=i_4\ne 0\}}
{\bf 1}_{\{j_1=j_4\}}
{\bf 1}_{\{i_2=i_3\ne 0\}}
{\bf 1}_{\{j_2=j_3\}}\Biggr),
\end{equation}

\vspace{6mm}
$$
J[\psi^{(5)}]_{T,t}
=\hbox{\vtop{\offinterlineskip\halign{
\hfil#\hfil\cr
{\rm l.i.m.}\cr
$\stackrel{}{{}_{p_1,\ldots,p_5\to \infty}}$\cr
}} }\sum_{j_1=0}^{p_1}\ldots\sum_{j_5=0}^{p_5}
C_{j_5\ldots j_1}\Biggl(
\prod_{l=1}^5\zeta_{j_l}^{(i_l)}
-\Biggr.
$$
$$
-
{\bf 1}_{\{i_1=i_2\ne 0\}}
{\bf 1}_{\{j_1=j_2\}}
\zeta_{j_3}^{(i_3)}
\zeta_{j_4}^{(i_4)}
\zeta_{j_5}^{(i_5)}-
{\bf 1}_{\{i_1=i_3\ne 0\}}
{\bf 1}_{\{j_1=j_3\}}
\zeta_{j_2}^{(i_2)}
\zeta_{j_4}^{(i_4)}
\zeta_{j_5}^{(i_5)}-
$$
$$
-
{\bf 1}_{\{i_1=i_4\ne 0\}}
{\bf 1}_{\{j_1=j_4\}}
\zeta_{j_2}^{(i_2)}
\zeta_{j_3}^{(i_3)}
\zeta_{j_5}^{(i_5)}-
{\bf 1}_{\{i_1=i_5\ne 0\}}
{\bf 1}_{\{j_1=j_5\}}
\zeta_{j_2}^{(i_2)}
\zeta_{j_3}^{(i_3)}
\zeta_{j_4}^{(i_4)}-
$$
$$
-
{\bf 1}_{\{i_2=i_3\ne 0\}}
{\bf 1}_{\{j_2=j_3\}}
\zeta_{j_1}^{(i_1)}
\zeta_{j_4}^{(i_4)}
\zeta_{j_5}^{(i_5)}-
{\bf 1}_{\{i_2=i_4\ne 0\}}
{\bf 1}_{\{j_2=j_4\}}
\zeta_{j_1}^{(i_1)}
\zeta_{j_3}^{(i_3)}
\zeta_{j_5}^{(i_5)}-
$$
$$
-
{\bf 1}_{\{i_2=i_5\ne 0\}}
{\bf 1}_{\{j_2=j_5\}}
\zeta_{j_1}^{(i_1)}
\zeta_{j_3}^{(i_3)}
\zeta_{j_4}^{(i_4)}
-{\bf 1}_{\{i_3=i_4\ne 0\}}
{\bf 1}_{\{j_3=j_4\}}
\zeta_{j_1}^{(i_1)}
\zeta_{j_2}^{(i_2)}
\zeta_{j_5}^{(i_5)}-
$$
$$
-
{\bf 1}_{\{i_3=i_5\ne 0\}}
{\bf 1}_{\{j_3=j_5\}}
\zeta_{j_1}^{(i_1)}
\zeta_{j_2}^{(i_2)}
\zeta_{j_4}^{(i_4)}
-{\bf 1}_{\{i_4=i_5\ne 0\}}
{\bf 1}_{\{j_4=j_5\}}
\zeta_{j_1}^{(i_1)}
\zeta_{j_2}^{(i_2)}
\zeta_{j_3}^{(i_3)}+
$$
$$
+
{\bf 1}_{\{i_1=i_2\ne 0\}}
{\bf 1}_{\{j_1=j_2\}}
{\bf 1}_{\{i_3=i_4\ne 0\}}
{\bf 1}_{\{j_3=j_4\}}\zeta_{j_5}^{(i_5)}+
{\bf 1}_{\{i_1=i_2\ne 0\}}
{\bf 1}_{\{j_1=j_2\}}
{\bf 1}_{\{i_3=i_5\ne 0\}}
{\bf 1}_{\{j_3=j_5\}}\zeta_{j_4}^{(i_4)}+
$$
$$
+
{\bf 1}_{\{i_1=i_2\ne 0\}}
{\bf 1}_{\{j_1=j_2\}}
{\bf 1}_{\{i_4=i_5\ne 0\}}
{\bf 1}_{\{j_4=j_5\}}\zeta_{j_3}^{(i_3)}+
{\bf 1}_{\{i_1=i_3\ne 0\}}
{\bf 1}_{\{j_1=j_3\}}
{\bf 1}_{\{i_2=i_4\ne 0\}}
{\bf 1}_{\{j_2=j_4\}}\zeta_{j_5}^{(i_5)}+
$$
$$
+
{\bf 1}_{\{i_1=i_3\ne 0\}}
{\bf 1}_{\{j_1=j_3\}}
{\bf 1}_{\{i_2=i_5\ne 0\}}
{\bf 1}_{\{j_2=j_5\}}\zeta_{j_4}^{(i_4)}+
{\bf 1}_{\{i_1=i_3\ne 0\}}
{\bf 1}_{\{j_1=j_3\}}
{\bf 1}_{\{i_4=i_5\ne 0\}}
{\bf 1}_{\{j_4=j_5\}}\zeta_{j_2}^{(i_2)}+
$$
$$
+
{\bf 1}_{\{i_1=i_4\ne 0\}}
{\bf 1}_{\{j_1=j_4\}}
{\bf 1}_{\{i_2=i_3\ne 0\}}
{\bf 1}_{\{j_2=j_3\}}\zeta_{j_5}^{(i_5)}+
{\bf 1}_{\{i_1=i_4\ne 0\}}
{\bf 1}_{\{j_1=j_4\}}
{\bf 1}_{\{i_2=i_5\ne 0\}}
{\bf 1}_{\{j_2=j_5\}}\zeta_{j_3}^{(i_3)}+
$$
$$
+
{\bf 1}_{\{i_1=i_4\ne 0\}}
{\bf 1}_{\{j_1=j_4\}}
{\bf 1}_{\{i_3=i_5\ne 0\}}
{\bf 1}_{\{j_3=j_5\}}\zeta_{j_2}^{(i_2)}+
{\bf 1}_{\{i_1=i_5\ne 0\}}
{\bf 1}_{\{j_1=j_5\}}
{\bf 1}_{\{i_2=i_3\ne 0\}}
{\bf 1}_{\{j_2=j_3\}}\zeta_{j_4}^{(i_4)}+
$$
$$
+
{\bf 1}_{\{i_1=i_5\ne 0\}}
{\bf 1}_{\{j_1=j_5\}}
{\bf 1}_{\{i_2=i_4\ne 0\}}
{\bf 1}_{\{j_2=j_4\}}\zeta_{j_3}^{(i_3)}+
{\bf 1}_{\{i_1=i_5\ne 0\}}
{\bf 1}_{\{j_1=j_5\}}
{\bf 1}_{\{i_3=i_4\ne 0\}}
{\bf 1}_{\{j_3=j_4\}}\zeta_{j_2}^{(i_2)}+
$$
$$
+
{\bf 1}_{\{i_2=i_3\ne 0\}}
{\bf 1}_{\{j_2=j_3\}}
{\bf 1}_{\{i_4=i_5\ne 0\}}
{\bf 1}_{\{j_4=j_5\}}\zeta_{j_1}^{(i_1)}+
{\bf 1}_{\{i_2=i_4\ne 0\}}
{\bf 1}_{\{j_2=j_4\}}
{\bf 1}_{\{i_3=i_5\ne 0\}}
{\bf 1}_{\{j_3=j_5\}}\zeta_{j_1}^{(i_1)}+
$$
\begin{equation}
\label{a5}
+\Biggl.
{\bf 1}_{\{i_2=i_5\ne 0\}}
{\bf 1}_{\{j_2=j_5\}}
{\bf 1}_{\{i_3=i_4\ne 0\}}
{\bf 1}_{\{j_3=j_4\}}\zeta_{j_1}^{(i_1)}\Biggr),
\end{equation}

\vspace{8mm}

$$
J[\psi^{(6)}]_{T,t}
=\hbox{\vtop{\offinterlineskip\halign{
\hfil#\hfil\cr
{\rm l.i.m.}\cr
$\stackrel{}{{}_{p_1,\ldots,p_6\to \infty}}$\cr
}} }\sum_{j_1=0}^{p_1}\ldots\sum_{j_6=0}^{p_6}
C_{j_6\ldots j_1}\Biggl(
\prod_{l=1}^6
\zeta_{j_l}^{(i_l)}
-\Biggr.
$$
$$
-
{\bf 1}_{\{i_1=i_6\ne 0\}}
{\bf 1}_{\{j_1=j_6\}}
\zeta_{j_2}^{(i_2)}
\zeta_{j_3}^{(i_3)}
\zeta_{j_4}^{(i_4)}
\zeta_{j_5}^{(i_5)}-
{\bf 1}_{\{i_2=i_6\ne 0\}}
{\bf 1}_{\{j_2=j_6\}}
\zeta_{j_1}^{(i_1)}
\zeta_{j_3}^{(i_3)}
\zeta_{j_4}^{(i_4)}
\zeta_{j_5}^{(i_5)}-
$$
$$
-
{\bf 1}_{\{i_3=i_6\ne 0\}}
{\bf 1}_{\{j_3=j_6\}}
\zeta_{j_1}^{(i_1)}
\zeta_{j_2}^{(i_2)}
\zeta_{j_4}^{(i_4)}
\zeta_{j_5}^{(i_5)}-
{\bf 1}_{\{i_4=i_6\ne 0\}}
{\bf 1}_{\{j_4=j_6\}}
\zeta_{j_1}^{(i_1)}
\zeta_{j_2}^{(i_2)}
\zeta_{j_3}^{(i_3)}
\zeta_{j_5}^{(i_5)}-
$$
$$
-
{\bf 1}_{\{i_5=i_6\ne 0\}}
{\bf 1}_{\{j_5=j_6\}}
\zeta_{j_1}^{(i_1)}
\zeta_{j_2}^{(i_2)}
\zeta_{j_3}^{(i_3)}
\zeta_{j_4}^{(i_4)}-
{\bf 1}_{\{i_1=i_2\ne 0\}}
{\bf 1}_{\{j_1=j_2\}}
\zeta_{j_3}^{(i_3)}
\zeta_{j_4}^{(i_4)}
\zeta_{j_5}^{(i_5)}
\zeta_{j_6}^{(i_6)}-
$$
$$
-
{\bf 1}_{\{i_1=i_3\ne 0\}}
{\bf 1}_{\{j_1=j_3\}}
\zeta_{j_2}^{(i_2)}
\zeta_{j_4}^{(i_4)}
\zeta_{j_5}^{(i_5)}
\zeta_{j_6}^{(i_6)}-
{\bf 1}_{\{i_1=i_4\ne 0\}}
{\bf 1}_{\{j_1=j_4\}}
\zeta_{j_2}^{(i_2)}
\zeta_{j_3}^{(i_3)}
\zeta_{j_5}^{(i_5)}
\zeta_{j_6}^{(i_6)}-
$$
$$
-
{\bf 1}_{\{i_1=i_5\ne 0\}}
{\bf 1}_{\{j_1=j_5\}}
\zeta_{j_2}^{(i_2)}
\zeta_{j_3}^{(i_3)}
\zeta_{j_4}^{(i_4)}
\zeta_{j_6}^{(i_6)}-
{\bf 1}_{\{i_2=i_3\ne 0\}}
{\bf 1}_{\{j_2=j_3\}}
\zeta_{j_1}^{(i_1)}
\zeta_{j_4}^{(i_4)}
\zeta_{j_5}^{(i_5)}
\zeta_{j_6}^{(i_6)}-
$$
$$
-
{\bf 1}_{\{i_2=i_4\ne 0\}}
{\bf 1}_{\{j_2=j_4\}}
\zeta_{j_1}^{(i_1)}
\zeta_{j_3}^{(i_3)}
\zeta_{j_5}^{(i_5)}
\zeta_{j_6}^{(i_6)}-
{\bf 1}_{\{i_2=i_5\ne 0\}}
{\bf 1}_{\{j_2=j_5\}}
\zeta_{j_1}^{(i_1)}
\zeta_{j_3}^{(i_3)}
\zeta_{j_4}^{(i_4)}
\zeta_{j_6}^{(i_6)}-
$$
$$
-
{\bf 1}_{\{i_3=i_4\ne 0\}}
{\bf 1}_{\{j_3=j_4\}}
\zeta_{j_1}^{(i_1)}
\zeta_{j_2}^{(i_2)}
\zeta_{j_5}^{(i_5)}
\zeta_{j_6}^{(i_6)}-
{\bf 1}_{\{i_3=i_5\ne 0\}}
{\bf 1}_{\{j_3=j_5\}}
\zeta_{j_1}^{(i_1)}
\zeta_{j_2}^{(i_2)}
\zeta_{j_4}^{(i_4)}
\zeta_{j_6}^{(i_6)}-
$$
$$
-
{\bf 1}_{\{i_4=i_5\ne 0\}}
{\bf 1}_{\{j_4=j_5\}}
\zeta_{j_1}^{(i_1)}
\zeta_{j_2}^{(i_2)}
\zeta_{j_3}^{(i_3)}
\zeta_{j_6}^{(i_6)}+
$$
$$
+
{\bf 1}_{\{i_1=i_2\ne 0\}}
{\bf 1}_{\{j_1=j_2\}}
{\bf 1}_{\{i_3=i_4\ne 0\}}
{\bf 1}_{\{j_3=j_4\}}
\zeta_{j_5}^{(i_5)}
\zeta_{j_6}^{(i_6)}+
{\bf 1}_{\{i_1=i_2\ne 0\}}
{\bf 1}_{\{j_1=j_2\}}
{\bf 1}_{\{i_3=i_5\ne 0\}}
{\bf 1}_{\{j_3=j_5\}}
\zeta_{j_4}^{(i_4)}
\zeta_{j_6}^{(i_6)}+
$$
$$
+
{\bf 1}_{\{i_1=i_2\ne 0\}}
{\bf 1}_{\{j_1=j_2\}}
{\bf 1}_{\{i_4=i_5\ne 0\}}
{\bf 1}_{\{j_4=j_5\}}
\zeta_{j_3}^{(i_3)}
\zeta_{j_6}^{(i_6)}
+
{\bf 1}_{\{i_1=i_3\ne 0\}}
{\bf 1}_{\{j_1=j_3\}}
{\bf 1}_{\{i_2=i_4\ne 0\}}
{\bf 1}_{\{j_2=j_4\}}
\zeta_{j_5}^{(i_5)}
\zeta_{j_6}^{(i_6)}+
$$
$$
+
{\bf 1}_{\{i_1=i_3\ne 0\}}
{\bf 1}_{\{j_1=j_3\}}
{\bf 1}_{\{i_2=i_5\ne 0\}}
{\bf 1}_{\{j_2=j_5\}}
\zeta_{j_4}^{(i_4)}
\zeta_{j_6}^{(i_6)}
+{\bf 1}_{\{i_1=i_3\ne 0\}}
{\bf 1}_{\{j_1=j_3\}}
{\bf 1}_{\{i_4=i_5\ne 0\}}
{\bf 1}_{\{j_4=j_5\}}
\zeta_{j_2}^{(i_2)}
\zeta_{j_6}^{(i_6)}+
$$
$$
+
{\bf 1}_{\{i_1=i_4\ne 0\}}
{\bf 1}_{\{j_1=j_4\}}
{\bf 1}_{\{i_2=i_3\ne 0\}}
{\bf 1}_{\{j_2=j_3\}}
\zeta_{j_5}^{(i_5)}
\zeta_{j_6}^{(i_6)}
+
{\bf 1}_{\{i_1=i_4\ne 0\}}
{\bf 1}_{\{j_1=j_4\}}
{\bf 1}_{\{i_2=i_5\ne 0\}}
{\bf 1}_{\{j_2=j_5\}}
\zeta_{j_3}^{(i_3)}
\zeta_{j_6}^{(i_6)}+
$$
$$
+
{\bf 1}_{\{i_1=i_4\ne 0\}}
{\bf 1}_{\{j_1=j_4\}}
{\bf 1}_{\{i_3=i_5\ne 0\}}
{\bf 1}_{\{j_3=j_5\}}
\zeta_{j_2}^{(i_2)}
\zeta_{j_6}^{(i_6)}
+
{\bf 1}_{\{i_1=i_5\ne 0\}}
{\bf 1}_{\{j_1=j_5\}}
{\bf 1}_{\{i_2=i_3\ne 0\}}
{\bf 1}_{\{j_2=j_3\}}
\zeta_{j_4}^{(i_4)}
\zeta_{j_6}^{(i_6)}+
$$
$$
+
{\bf 1}_{\{i_1=i_5\ne 0\}}
{\bf 1}_{\{j_1=j_5\}}
{\bf 1}_{\{i_2=i_4\ne 0\}}
{\bf 1}_{\{j_2=j_4\}}
\zeta_{j_3}^{(i_3)}
\zeta_{j_6}^{(i_6)}
+
{\bf 1}_{\{i_1=i_5\ne 0\}}
{\bf 1}_{\{j_1=j_5\}}
{\bf 1}_{\{i_3=i_4\ne 0\}}
{\bf 1}_{\{j_3=j_4\}}
\zeta_{j_2}^{(i_2)}
\zeta_{j_6}^{(i_6)}+
$$
$$
+
{\bf 1}_{\{i_2=i_3\ne 0\}}
{\bf 1}_{\{j_2=j_3\}}
{\bf 1}_{\{i_4=i_5\ne 0\}}
{\bf 1}_{\{j_4=j_5\}}
\zeta_{j_1}^{(i_1)}
\zeta_{j_6}^{(i_6)}
+
{\bf 1}_{\{i_2=i_4\ne 0\}}
{\bf 1}_{\{j_2=j_4\}}
{\bf 1}_{\{i_3=i_5\ne 0\}}
{\bf 1}_{\{j_3=j_5\}}
\zeta_{j_1}^{(i_1)}
\zeta_{j_6}^{(i_6)}+
$$
$$
+
{\bf 1}_{\{i_2=i_5\ne 0\}}
{\bf 1}_{\{j_2=j_5\}}
{\bf 1}_{\{i_3=i_4\ne 0\}}
{\bf 1}_{\{j_3=j_4\}}
\zeta_{j_1}^{(i_1)}
\zeta_{j_6}^{(i_6)}
+
{\bf 1}_{\{i_6=i_1\ne 0\}}
{\bf 1}_{\{j_6=j_1\}}
{\bf 1}_{\{i_3=i_4\ne 0\}}
{\bf 1}_{\{j_3=j_4\}}
\zeta_{j_2}^{(i_2)}
\zeta_{j_5}^{(i_5)}+
$$
$$
+
{\bf 1}_{\{i_6=i_1\ne 0\}}
{\bf 1}_{\{j_6=j_1\}}
{\bf 1}_{\{i_3=i_5\ne 0\}}
{\bf 1}_{\{j_3=j_5\}}
\zeta_{j_2}^{(i_2)}
\zeta_{j_4}^{(i_4)}
+
{\bf 1}_{\{i_6=i_1\ne 0\}}
{\bf 1}_{\{j_6=j_1\}}
{\bf 1}_{\{i_2=i_5\ne 0\}}
{\bf 1}_{\{j_2=j_5\}}
\zeta_{j_3}^{(i_3)}
\zeta_{j_4}^{(i_4)}+
$$
$$
+
{\bf 1}_{\{i_6=i_1\ne 0\}}
{\bf 1}_{\{j_6=j_1\}}
{\bf 1}_{\{i_2=i_4\ne 0\}}
{\bf 1}_{\{j_2=j_4\}}
\zeta_{j_3}^{(i_3)}
\zeta_{j_5}^{(i_5)}
+
{\bf 1}_{\{i_6=i_1\ne 0\}}
{\bf 1}_{\{j_6=j_1\}}
{\bf 1}_{\{i_4=i_5\ne 0\}}
{\bf 1}_{\{j_4=j_5\}}
\zeta_{j_2}^{(i_2)}
\zeta_{j_3}^{(i_3)}+
$$
$$
+
{\bf 1}_{\{i_6=i_1\ne 0\}}
{\bf 1}_{\{j_6=j_1\}}
{\bf 1}_{\{i_2=i_3\ne 0\}}
{\bf 1}_{\{j_2=j_3\}}
\zeta_{j_4}^{(i_4)}
\zeta_{j_5}^{(i_5)}
+
{\bf 1}_{\{i_6=i_2\ne 0\}}
{\bf 1}_{\{j_6=j_2\}}
{\bf 1}_{\{i_3=i_5\ne 0\}}
{\bf 1}_{\{j_3=j_5\}}
\zeta_{j_1}^{(i_1)}
\zeta_{j_4}^{(i_4)}+
$$
$$
+
{\bf 1}_{\{i_6=i_2\ne 0\}}
{\bf 1}_{\{j_6=j_2\}}
{\bf 1}_{\{i_4=i_5\ne 0\}}
{\bf 1}_{\{j_4=j_5\}}
\zeta_{j_1}^{(i_1)}
\zeta_{j_3}^{(i_3)}
+
{\bf 1}_{\{i_6=i_2\ne 0\}}
{\bf 1}_{\{j_6=j_2\}}
{\bf 1}_{\{i_3=i_4\ne 0\}}
{\bf 1}_{\{j_3=j_4\}}
\zeta_{j_1}^{(i_1)}
\zeta_{j_5}^{(i_5)}+
$$
$$
+
{\bf 1}_{\{i_6=i_2\ne 0\}}
{\bf 1}_{\{j_6=j_2\}}
{\bf 1}_{\{i_1=i_5\ne 0\}}
{\bf 1}_{\{j_1=j_5\}}
\zeta_{j_3}^{(i_3)}
\zeta_{j_4}^{(i_4)}
+
{\bf 1}_{\{i_6=i_2\ne 0\}}
{\bf 1}_{\{j_6=j_2\}}
{\bf 1}_{\{i_1=i_4\ne 0\}}
{\bf 1}_{\{j_1=j_4\}}
\zeta_{j_3}^{(i_3)}
\zeta_{j_5}^{(i_5)}+
$$
$$
+
{\bf 1}_{\{i_6=i_2\ne 0\}}
{\bf 1}_{\{j_6=j_2\}}
{\bf 1}_{\{i_1=i_3\ne 0\}}
{\bf 1}_{\{j_1=j_3\}}
\zeta_{j_4}^{(i_4)}
\zeta_{j_5}^{(i_5)}
+
{\bf 1}_{\{i_6=i_3\ne 0\}}
{\bf 1}_{\{j_6=j_3\}}
{\bf 1}_{\{i_2=i_5\ne 0\}}
{\bf 1}_{\{j_2=j_5\}}
\zeta_{j_1}^{(i_1)}
\zeta_{j_4}^{(i_4)}+
$$
$$
+
{\bf 1}_{\{i_6=i_3\ne 0\}}
{\bf 1}_{\{j_6=j_3\}}
{\bf 1}_{\{i_4=i_5\ne 0\}}
{\bf 1}_{\{j_4=j_5\}}
\zeta_{j_1}^{(i_1)}
\zeta_{j_2}^{(i_2)}
+
{\bf 1}_{\{i_6=i_3\ne 0\}}
{\bf 1}_{\{j_6=j_3\}}
{\bf 1}_{\{i_2=i_4\ne 0\}}
{\bf 1}_{\{j_2=j_4\}}
\zeta_{j_1}^{(i_1)}
\zeta_{j_5}^{(i_5)}+
$$
$$
+
{\bf 1}_{\{i_6=i_3\ne 0\}}
{\bf 1}_{\{j_6=j_3\}}
{\bf 1}_{\{i_1=i_5\ne 0\}}
{\bf 1}_{\{j_1=j_5\}}
\zeta_{j_2}^{(i_2)}
\zeta_{j_4}^{(i_4)}
+
{\bf 1}_{\{i_6=i_3\ne 0\}}
{\bf 1}_{\{j_6=j_3\}}
{\bf 1}_{\{i_1=i_4\ne 0\}}
{\bf 1}_{\{j_1=j_4\}}
\zeta_{j_2}^{(i_2)}
\zeta_{j_5}^{(i_5)}+
$$
$$
+
{\bf 1}_{\{i_6=i_3\ne 0\}}
{\bf 1}_{\{j_6=j_3\}}
{\bf 1}_{\{i_1=i_2\ne 0\}}
{\bf 1}_{\{j_1=j_2\}}
\zeta_{j_4}^{(i_4)}
\zeta_{j_5}^{(i_5)}
+
{\bf 1}_{\{i_6=i_4\ne 0\}}
{\bf 1}_{\{j_6=j_4\}}
{\bf 1}_{\{i_3=i_5\ne 0\}}
{\bf 1}_{\{j_3=j_5\}}
\zeta_{j_1}^{(i_1)}
\zeta_{j_2}^{(i_2)}+
$$
$$
+
{\bf 1}_{\{i_6=i_4\ne 0\}}
{\bf 1}_{\{j_6=j_4\}}
{\bf 1}_{\{i_2=i_5\ne 0\}}
{\bf 1}_{\{j_2=j_5\}}
\zeta_{j_1}^{(i_1)}
\zeta_{j_3}^{(i_3)}
+
{\bf 1}_{\{i_6=i_4\ne 0\}}
{\bf 1}_{\{j_6=j_4\}}
{\bf 1}_{\{i_2=i_3\ne 0\}}
{\bf 1}_{\{j_2=j_3\}}
\zeta_{j_1}^{(i_1)}
\zeta_{j_5}^{(i_5)}+
$$
$$
+
{\bf 1}_{\{i_6=i_4\ne 0\}}
{\bf 1}_{\{j_6=j_4\}}
{\bf 1}_{\{i_1=i_5\ne 0\}}
{\bf 1}_{\{j_1=j_5\}}
\zeta_{j_2}^{(i_2)}
\zeta_{j_3}^{(i_3)}
+
{\bf 1}_{\{i_6=i_4\ne 0\}}
{\bf 1}_{\{j_6=j_4\}}
{\bf 1}_{\{i_1=i_3\ne 0\}}
{\bf 1}_{\{j_1=j_3\}}
\zeta_{j_2}^{(i_2)}
\zeta_{j_5}^{(i_5)}+
$$
$$
+
{\bf 1}_{\{i_6=i_4\ne 0\}}
{\bf 1}_{\{j_6=j_4\}}
{\bf 1}_{\{i_1=i_2\ne 0\}}
{\bf 1}_{\{j_1=j_2\}}
\zeta_{j_3}^{(i_3)}
\zeta_{j_5}^{(i_5)}
+
{\bf 1}_{\{i_6=i_5\ne 0\}}
{\bf 1}_{\{j_6=j_5\}}
{\bf 1}_{\{i_3=i_4\ne 0\}}
{\bf 1}_{\{j_3=j_4\}}
\zeta_{j_1}^{(i_1)}
\zeta_{j_2}^{(i_2)}+
$$
$$
+
{\bf 1}_{\{i_6=i_5\ne 0\}}
{\bf 1}_{\{j_6=j_5\}}
{\bf 1}_{\{i_2=i_4\ne 0\}}
{\bf 1}_{\{j_2=j_4\}}
\zeta_{j_1}^{(i_1)}
\zeta_{j_3}^{(i_3)}
+
{\bf 1}_{\{i_6=i_5\ne 0\}}
{\bf 1}_{\{j_6=j_5\}}
{\bf 1}_{\{i_2=i_3\ne 0\}}
{\bf 1}_{\{j_2=j_3\}}
\zeta_{j_1}^{(i_1)}
\zeta_{j_4}^{(i_4)}+
$$
$$
+
{\bf 1}_{\{i_6=i_5\ne 0\}}
{\bf 1}_{\{j_6=j_5\}}
{\bf 1}_{\{i_1=i_4\ne 0\}}
{\bf 1}_{\{j_1=j_4\}}
\zeta_{j_2}^{(i_2)}
\zeta_{j_3}^{(i_3)}
+
{\bf 1}_{\{i_6=i_5\ne 0\}}
{\bf 1}_{\{j_6=j_5\}}
{\bf 1}_{\{i_1=i_3\ne 0\}}
{\bf 1}_{\{j_1=j_3\}}
\zeta_{j_2}^{(i_2)}
\zeta_{j_4}^{(i_4)}+
$$
$$
+
{\bf 1}_{\{i_6=i_5\ne 0\}}
{\bf 1}_{\{j_6=j_5\}}
{\bf 1}_{\{i_1=i_2\ne 0\}}
{\bf 1}_{\{j_1=j_2\}}
\zeta_{j_3}^{(i_3)}
\zeta_{j_4}^{(i_4)}-
$$
$$
-
{\bf 1}_{\{i_6=i_1\ne 0\}}
{\bf 1}_{\{j_6=j_1\}}
{\bf 1}_{\{i_2=i_5\ne 0\}}
{\bf 1}_{\{j_2=j_5\}}
{\bf 1}_{\{i_3=i_4\ne 0\}}
{\bf 1}_{\{j_3=j_4\}}-
$$
$$
-
{\bf 1}_{\{i_6=i_1\ne 0\}}
{\bf 1}_{\{j_6=j_1\}}
{\bf 1}_{\{i_2=i_4\ne 0\}}
{\bf 1}_{\{j_2=j_4\}}
{\bf 1}_{\{i_3=i_5\ne 0\}}
{\bf 1}_{\{j_3=j_5\}}-
$$
$$
-
{\bf 1}_{\{i_6=i_1\ne 0\}}
{\bf 1}_{\{j_6=j_1\}}
{\bf 1}_{\{i_2=i_3\ne 0\}}
{\bf 1}_{\{j_2=j_3\}}
{\bf 1}_{\{i_4=i_5\ne 0\}}
{\bf 1}_{\{j_4=j_5\}}-
$$
$$
-
{\bf 1}_{\{i_6=i_2\ne 0\}}
{\bf 1}_{\{j_6=j_2\}}
{\bf 1}_{\{i_1=i_5\ne 0\}}
{\bf 1}_{\{j_1=j_5\}}
{\bf 1}_{\{i_3=i_4\ne 0\}}
{\bf 1}_{\{j_3=j_4\}}-
$$
$$
-
{\bf 1}_{\{i_6=i_2\ne 0\}}
{\bf 1}_{\{j_6=j_2\}}
{\bf 1}_{\{i_1=i_4\ne 0\}}
{\bf 1}_{\{j_1=j_4\}}
{\bf 1}_{\{i_3=i_5\ne 0\}}
{\bf 1}_{\{j_3=j_5\}}-
$$
$$
-
{\bf 1}_{\{i_6=i_2\ne 0\}}
{\bf 1}_{\{j_6=j_2\}}
{\bf 1}_{\{i_1=i_3\ne 0\}}
{\bf 1}_{\{j_1=j_3\}}
{\bf 1}_{\{i_4=i_5\ne 0\}}
{\bf 1}_{\{j_4=j_5\}}-
$$
$$
-
{\bf 1}_{\{i_6=i_3\ne 0\}}
{\bf 1}_{\{j_6=j_3\}}
{\bf 1}_{\{i_1=i_5\ne 0\}}
{\bf 1}_{\{j_1=j_5\}}
{\bf 1}_{\{i_2=i_4\ne 0\}}
{\bf 1}_{\{j_2=j_4\}}-
$$
$$
-
{\bf 1}_{\{i_6=i_3\ne 0\}}
{\bf 1}_{\{j_6=j_3\}}
{\bf 1}_{\{i_1=i_4\ne 0\}}
{\bf 1}_{\{j_1=j_4\}}
{\bf 1}_{\{i_2=i_5\ne 0\}}
{\bf 1}_{\{j_2=j_5\}}-
$$
$$
-
{\bf 1}_{\{i_3=i_6\ne 0\}}
{\bf 1}_{\{j_3=j_6\}}
{\bf 1}_{\{i_1=i_2\ne 0\}}
{\bf 1}_{\{j_1=j_2\}}
{\bf 1}_{\{i_4=i_5\ne 0\}}
{\bf 1}_{\{j_4=j_5\}}-
$$
$$
-
{\bf 1}_{\{i_6=i_4\ne 0\}}
{\bf 1}_{\{j_6=j_4\}}
{\bf 1}_{\{i_1=i_5\ne 0\}}
{\bf 1}_{\{j_1=j_5\}}
{\bf 1}_{\{i_2=i_3\ne 0\}}
{\bf 1}_{\{j_2=j_3\}}-
$$
$$
-
{\bf 1}_{\{i_6=i_4\ne 0\}}
{\bf 1}_{\{j_6=j_4\}}
{\bf 1}_{\{i_1=i_3\ne 0\}}
{\bf 1}_{\{j_1=j_3\}}
{\bf 1}_{\{i_2=i_5\ne 0\}}
{\bf 1}_{\{j_2=j_5\}}-
$$
$$
-
{\bf 1}_{\{i_6=i_4\ne 0\}}
{\bf 1}_{\{j_6=j_4\}}
{\bf 1}_{\{i_1=i_2\ne 0\}}
{\bf 1}_{\{j_1=j_2\}}
{\bf 1}_{\{i_3=i_5\ne 0\}}
{\bf 1}_{\{j_3=j_5\}}-
$$
$$
-
{\bf 1}_{\{i_6=i_5\ne 0\}}
{\bf 1}_{\{j_6=j_5\}}
{\bf 1}_{\{i_1=i_4\ne 0\}}
{\bf 1}_{\{j_1=j_4\}}
{\bf 1}_{\{i_2=i_3\ne 0\}}
{\bf 1}_{\{j_2=j_3\}}-
$$
$$
-
{\bf 1}_{\{i_6=i_5\ne 0\}}
{\bf 1}_{\{j_6=j_5\}}
{\bf 1}_{\{i_1=i_2\ne 0\}}
{\bf 1}_{\{j_1=j_2\}}
{\bf 1}_{\{i_3=i_4\ne 0\}}
{\bf 1}_{\{j_3=j_4\}}-
$$
\begin{equation}
\label{a6}
\Biggl.-
{\bf 1}_{\{i_6=i_5\ne 0\}}
{\bf 1}_{\{j_6=j_5\}}
{\bf 1}_{\{i_1=i_3\ne 0\}}
{\bf 1}_{\{j_1=j_3\}}
{\bf 1}_{\{i_2=i_4\ne 0\}}
{\bf 1}_{\{j_2=j_4\}}\Biggr),
\end{equation}

\vspace{7mm}
\noindent
where ${\bf 1}_A$ is the indicator of the set $A$.

The convergence 
in the mean of degree $2n$ ($n\in \mathbb{N}$) 
in Theorem 2 is proved in
in \cite{20a} (Sect.~1.1.9, 1.11, 1.12), \cite{26a} (Sect.~6, 15, 16).
Moreover, the complete orthonormal systems of Haar and 
Rade\-ma\-cher--Walsh functions in the space $L_2([t,T])$ also
can be applied in Theorem 2
\cite{9}-\cite{16}, \cite{19}, \cite{20}-\cite{20a-new-x}, \cite{26a}.
The convergence w.~p.~1 in Theorem 2 
is proved in \cite{20a}-\cite{20a-new-x}, \cite{22}, \cite{26},
\cite{30b}
for complete orthonormal systems of Legendre polynomials 
and trigonometric functions
in the space $L_2([t,T])$.
The modifications of Theorem 2 
were obtained in \cite{20}-\cite{20a-new-x}, \cite{26b}
for complete 
orthonormal with weight  
$r(t_1)\ldots r(t_k)\ge 0$
systems of functions in the space $L_2([t,T]^k)$ ($k\in\mathbb{N}$) 
as well 
as for some other types of iterated stochastic 
integrals (iterated stochastic integrals 
with respect to martingale Poisson measures and 
iterated stochastic integrals with respect 
to martingales).
Application of Theorem 2 and Theorem 4 (see below) for the mean-square
approximation of iterated stochastic integrals 
with respect to the 
infinite-dimensional $Q$-Wiener process can be found
in the monographs \cite{20a}-\cite{20a-new-x} (Chapter 7) and in \cite{31}-\cite{31aaa}.

Note that the correctness of formulas (\ref{a1})--(\ref{a6}) 
can be 
verified 
by the fact that if 
$i_1=\ldots=i_6=i=1,\ldots,m$
and $\psi_1(\tau),\ldots,\psi_6(\tau)\equiv \psi(\tau)$,
then we can derive from (\ref{a1})--(\ref{a6}) the well known
equalities

$$
J[\psi^{(1)}]_{T,t}
=\frac{1}{1!}\delta_{T,t},
$$

$$
J[\psi^{(2)}]_{T,t}
=\frac{1}{2!}\left(\delta^2_{T,t}-\Delta_{T,t}\right),\
$$

$$
J[\psi^{(3)}]_{T,t}
=\frac{1}{3!}\left(\delta_{T,t}^3-3\delta_{T,t}\Delta_{T,t}\right),
$$

$$
J[\psi^{(4)}]_{T,t}
=\frac{1}{4!}\left(\delta^4_{T,t}-6\delta_{T,t}^2\Delta_{T,t}
+3\Delta^2_{T,t}\right),\
$$

$$
J[\psi^{(5)}]_{T,t}
=\frac{1}{5!}\left(\delta^5_{T,t}-10\delta_{T,t}^3\Delta_{T,t}
+15\delta_{T,t}\Delta^2_{T,t}\right),
$$

$$
J[\psi^{(6)}]_{T,t}
=\frac{1}{6!}\left(\delta^6_{T,t}-15\delta_{T,t}^4\Delta_{T,t}
+45\delta_{T,t}^2\Delta^2_{T,t}-15\Delta_{T,t}^3\right)
$$

\vspace{5mm}
\noindent
w.~p.~1 
\cite{7}-\cite{16}, \cite{19}-\cite{20a-new-x}, where 

$$
\delta_{T,t}=\int\limits_t^T\psi(\tau)d{\bf f}_{\tau}^{(i)},\ \ \
\Delta_{T,t}=\int\limits_t^T\psi^2(\tau)d\tau.
$$

\vspace{3mm}

Note that the mentioned equalities
can be independently  
obtained using the Ito formula and Hermite polynomials.

For further consideration, let us 
consider the generalization of formulas (\ref{a1})--(\ref{a6})                 
for the case of an arbitrary multiplicity $k$ $(k\in\mathbb{N})$ of 
the iterated Ito stochastic integral $J[\psi^{(k)}]_{T,t}$ defined by (\ref{ito}).
In order to do this, let us
introduce some notations. 
Consider the unordered
set $\{1, 2, \ldots, k\}$ 
and separate it into two parts:
the first part consists of $r$ unordered 
pairs (sequence order of these pairs is also unimportant) and the 
second one consists of the 
remaining $k-2r$ numbers.
So, we have

\begin{equation}
\label{leto5007}
(\{
\underbrace{\{g_1, g_2\}, \ldots, 
\{g_{2r-1}, g_{2r}\}}_{\small{\hbox{part 1}}}
\},
\{\underbrace{q_1, \ldots, q_{k-2r}}_{\small{\hbox{part 2}}}
\}),
\end{equation}

\vspace{2mm}
\noindent
where 
$$\{g_1, g_2, \ldots, 
g_{2r-1}, g_{2r}, q_1, \ldots, q_{k-2r}\}=\{1, 2, \ldots, k\},
$$

\vspace{4mm}
\noindent
braces   
mean an unordered 
set, and pa\-ren\-the\-ses mean an ordered set.

We will say that (\ref{leto5007}) is a partition 
and consider the sum with respect to all possible
partitions

\begin{equation}
\label{leto5008}
\sum_{\stackrel{(\{\{g_1, g_2\}, \ldots, 
\{g_{2r-1}, g_{2r}\}\}, \{q_1, \ldots, q_{k-2r}\})}
{{}_{\{g_1, g_2, \ldots, 
g_{2r-1}, g_{2r}, q_1, \ldots, q_{k-2r}\}=\{1, 2, \ldots, k\}}}}
a_{g_1 g_2, \ldots, 
g_{2r-1} g_{2r}, q_1 \ldots q_{k-2r}},
\end{equation}

\vspace{4mm}
\noindent
where $a_{g_1 g_2, \ldots, 
g_{2r-1} g_{2r}, q_1 \ldots q_{k-2r}}\in\mathbb{R}.$

Below there are several examples of sums in the form (\ref{leto5008})

\vspace{2mm}
$$
\sum_{\stackrel{(\{g_1, g_2\})}{{}_{\{g_1, g_2\}=\{1, 2\}}}}
a_{g_1 g_2}=a_{12},
$$

\vspace{3mm}
$$
\sum_{\stackrel{(\{\{g_1, g_2\}, \{g_3, g_4\}\})}
{{}_{\{g_1, g_2, g_3, g_4\}=\{1, 2, 3, 4\}}}}
a_{g_1 g_2, g_3 g_4}=a_{12,34} + a_{13,24} + a_{23,14},
$$

\vspace{3mm}
$$
\sum_{\stackrel{(\{g_1, g_2\}, \{q_1, q_{2}\})}
{{}_{\{g_1, g_2, q_1, q_{2}\}=\{1, 2, 3, 4\}}}}
a_{g_1 g_2, q_1 q_{2}}=
$$

$$
=a_{12,34}+a_{13,24}+a_{14,23}
+a_{23,14}+a_{24,13}+a_{34,12},
$$

\vspace{3mm}
$$
\sum_{\stackrel{(\{g_1, g_2\}, \{q_1, q_{2}, q_3\})}
{{}_{\{g_1, g_2, q_1, q_{2}, q_3\}=\{1, 2, 3, 4, 5\}}}}
a_{g_1 g_2, q_1 q_{2}q_3}
=
$$

$$
=a_{12,345}+a_{13,245}+a_{14,235}
+a_{15,234}+a_{23,145}+a_{24,135}+
$$
$$
+a_{25,134}+a_{34,125}+a_{35,124}+a_{45,123},
$$

\vspace{3mm}
$$
\sum_{\stackrel{(\{\{g_1, g_2\}, \{g_3, g_{4}\}\}, \{q_1\})}
{{}_{\{g_1, g_2, g_3, g_{4}, q_1\}=\{1, 2, 3, 4, 5\}}}}
a_{g_1 g_2, g_3 g_{4},q_1}
=
$$

$$
=
a_{12,34,5}+a_{13,24,5}+a_{14,23,5}+
a_{12,35,4}+a_{13,25,4}+a_{15,23,4}+
$$
$$
+a_{12,54,3}+a_{15,24,3}+a_{14,25,3}+a_{15,34,2}+a_{13,54,2}+a_{14,53,2}+
$$
$$
+
a_{52,34,1}+a_{53,24,1}+a_{54,23,1}.
$$

\vspace{6mm}

Now, we can formulate Theorem 2
(see (\ref{tyyy})) 
using the alternative form.

\vspace{2mm}

{\bf Theorem 3} \cite{10} (2009) (also see \cite{11}-\cite{16}, 
\cite{19}, \cite{20}-\cite{20a-new-x}.\
{\it Under the conditions of Theorem {\rm 2} 
the following expansion

$$
J[\psi^{(k)}]_{T,t}=
\hbox{\vtop{\offinterlineskip\halign{
\hfil#\hfil\cr
{\rm l.i.m.}\cr
$\stackrel{}{{}_{p_1,\ldots,p_k\to \infty}}$\cr
}} }
\sum\limits_{j_1=0}^{p_1}\ldots
\sum\limits_{j_k=0}^{p_k}
C_{j_k\ldots j_1}\Biggl(
\prod_{l=1}^k\zeta_{j_l}^{(i_l)}+\sum\limits_{r=1}^{[k/2]}
(-1)^r \times
\Biggr.
$$

\vspace{3mm}
\begin{equation}
\label{leto6000hh}
\times
\sum_{\stackrel{(\{\{g_1, g_2\}, \ldots, 
\{g_{2r-1}, g_{2r}\}\}, \{q_1, \ldots, q_{k-2r}\})}
{{}_{\{g_1, g_2, \ldots, 
g_{2r-1}, g_{2r}, q_1, \ldots, q_{k-2r}\}=\{1, 2, \ldots, k\}}}}
\prod\limits_{s=1}^r
{\bf 1}_{\{i_{g_{{}_{2s-1}}}=~i_{g_{{}_{2s}}}\ne 0\}}
\Biggl.{\bf 1}_{\{j_{g_{{}_{2s-1}}}=~j_{g_{{}_{2s}}}\}}
\prod_{l=1}^{k-2r}\zeta_{j_{q_l}}^{(i_{q_l})}\Biggr)
\end{equation}

\vspace{5mm}
\noindent
con\-verg\-ing in the mean-square sense is valid,
where $[x]$ is an integer part of a real number $x,$
$\prod\limits_{\emptyset}
\stackrel{\sf def}{=}1,$ $\sum\limits_{\emptyset}
\stackrel{\sf def}{=}0;$
another notations are the same as in Theorem {\rm 2.}}

\vspace{2mm}

In particular, from (\ref{leto6000hh}) for $k=5$ we obtain

\vspace{1mm}

$$
J[\psi^{(5)}]_{T,t}=
\hbox{\vtop{\offinterlineskip\halign{
\hfil#\hfil\cr
{\rm l.i.m.}\cr
$\stackrel{}{{}_{p_1,\ldots,p_5\to \infty}}$\cr
}} }\sum_{j_1=0}^{p_1}\ldots\sum_{j_5=0}^{p_5}
C_{j_5\ldots j_1}\Biggl(
\prod_{l=1}^5\zeta_{j_l}^{(i_l)}-\Biggr.
$$

\vspace{2mm}
$$
-
\sum\limits_{\stackrel{(\{g_1, g_2\}, \{q_1, q_{2}, q_3\})}
{{}_{\{g_1, g_2, q_{1}, q_{2}, q_3\}=\{1, 2, 3, 4, 5\}}}}
{\bf 1}_{\{i_{g_{{}_{1}}}=~i_{g_{{}_{2}}}\ne 0\}}
{\bf 1}_{\{j_{g_{{}_{1}}}=~j_{g_{{}_{2}}}\}}
\prod_{l=1}^{3}\zeta_{j_{q_l}}^{(i_{q_l})}+
$$

\vspace{2mm}
$$
+
\sum_{\stackrel{(\{\{g_1, g_2\}, 
\{g_{3}, g_{4}\}\}, \{q_1\})}
{{}_{\{g_1, g_2, g_{3}, g_{4}, q_1\}=\{1, 2, 3, 4, 5\}}}}
{\bf 1}_{\{i_{g_{{}_{1}}}=~i_{g_{{}_{2}}}\ne 0\}}
{\bf 1}_{\{j_{g_{{}_{1}}}=~j_{g_{{}_{2}}}\}}
\Biggl.{\bf 1}_{\{i_{g_{{}_{3}}}=~i_{g_{{}_{4}}}\ne 0\}}
{\bf 1}_{\{j_{g_{{}_{3}}}=~j_{g_{{}_{4}}}\}}
\zeta_{j_{q_1}}^{(i_{q_1})}\Biggr).
$$

\vspace{5mm}
\noindent
The last equality obviously agrees with
(\ref{a5}).

\vspace{5mm}

\section{A Generalization of Theorems 2, 3 to the Case of an Arbitrary 
Complete Ortho\-nor\-mal System of Functions in the Space $L_2([t, T])$
and $\psi_1(\tau),$ $\ldots,\psi_k(\tau)\in L_2([t, T])$}

\vspace{5mm}

In this section, we will use the definition of the multiple Wiener 
stochastic integral from \cite{ito1951}, \cite{Kuo} to generalize Theorems 
2, 3 to the case of an arbitrary 
complete orthonormal system of functions in the space $L_2([t, T])$
and $\psi_1(\tau),$ $\ldots,\psi_k(\tau)\in L_2([t, T]).$ 

Consider the following step function on the hypercube $[t, T]^k$

\vspace{-2mm}
\begin{equation}
\label{chain3}
\Phi_N(t_1,\ldots,t_k)=\sum\limits_{l_1,\ldots,l_k=0}^{N-1}
a_{l_1 \ldots l_k} {\bf 1}_{[\tau_{l_1},\tau_{l_1+1})}(t_1) \ldots
{\bf 1}_{[\tau_{l_k},\tau_{l_k+1})}(t_k),
\end{equation}

\vspace{3mm}
\noindent
where $a_{l_1 \ldots l_k}\in\mathbb{R}$ and such that 
$a_{l_1 \ldots l_k}=0$ if $l_p=l_q$ for some $p\ne q,$

\vspace{1mm}
$$
{\bf 1}_A (\tau)=\left\{
\begin{matrix}
1\ &{\rm if}\ \tau\in A \cr\cr
0\ &\hbox{\rm otherwise}
\end{matrix}\right.,
$$

\vspace{4mm}
\noindent
$N\in\mathbb{N},$ $\left\{\tau_{j}\right\}_{j=0}^{N}$ is a partition of
$[t,T],$ which satisfies the condition (\ref{1111}):     

\vspace{1mm}
\begin{equation}
\label{1111xxx1}
t=\tau_0<\ldots <\tau_N=T,\ \ \
\Delta_N=
\hbox{\vtop{\offinterlineskip\halign{
\hfil#\hfil\cr
{\rm max}\cr
$\stackrel{}{{}_{0\le j\le N-1}}$\cr
}} }\Delta\tau_j\to 0\ \ \hbox{if}\ \ N\to \infty,\ \ \ 
\Delta\tau_j=\tau_{j+1}-\tau_j.
\end{equation}

\vspace{4mm}

Let us define the multiple Wiener stochastic integral for $\Phi_N(t_1,\ldots,t_k)$ 
\cite{ito1951}, \cite{Kuo}

\vspace{1mm}
\begin{equation}
\label{chain9}
J'[\Phi_N]_{T,t}^{(i_1\ldots i_k)}\stackrel{\sf def}{=}
\sum\limits_{l_1,\ldots,l_k=0}^{N-1}
a_{l_1 \ldots l_k}
\Delta{\bf w}_{\tau_{l_1}}^{(i_1)}\ldots \Delta{\bf w}_{\tau_{l_k}}^{(i_k)},
\end{equation}

\vspace{4mm}
\noindent
where $\Delta{\bf w}_{\tau_{j}}^{(i)}=
{\bf w}_{\tau_{j+1}}^{(i)}-{\bf w}_{\tau_{j}}^{(i)},$\
$i=0, 1,\ldots,m,$\ ${\bf w}_{\tau}^{(0)}=\tau.$

It is known (see \cite{Kuo}, Lemma~9.6.4)
that for any $\Phi(t_1,\ldots,t_k)\in L_2([t, T]^k)$ 
there exists a sequence of step functions $\Phi_N(t_1,\ldots,t_k)$ of the form (\ref{chain3})
such that

\vspace{1mm}
\begin{equation}
\label{chain15}
\lim\limits_{N\to\infty} \int\limits_{[t,T]^k}
\left(\Phi(t_1,\ldots,t_k)-\Phi_N(t_1,\ldots,t_k)\right)^2 dt_1\ldots dt_k=0.
\end{equation}

\vspace{4mm}

We have

\vspace{-1mm}
$$
\Phi_N(t_1,\ldots,t_k)=\sum\limits_{l_1,\ldots,l_k=0}^{N-1}
a_{l_1 \ldots l_k} {\bf 1}_{[\tau_{l_1},\tau_{l_1+1})}(t_1) \ldots
{\bf 1}_{[\tau_{l_k},\tau_{l_k+1})}(t_k)=
$$

\vspace{1mm}
\begin{equation}
\label{chain5}
=\sum\limits_{(l_1,\ldots,l_k)}
\sum_{\stackrel{l_1,\ldots,l_k=0}{{}_{l_1<l_2<\ldots < l_k}}}^{N-1}
a_{l_1 \ldots l_k} {\bf 1}_{[\tau_{l_1},\tau_{l_1+1})}(t_1) \ldots
{\bf 1}_{[\tau_{l_k},\tau_{l_k+1})}(t_k),
\end{equation}

\vspace{4mm}
\noindent
where permutations $(l_1,\ldots,l_k)$ when summing are 
performed only in the expression $l_1<l_2<\ldots < l_k$
(recall that $a_{l_1 \ldots l_k}=0$ if $l_p=l_q$ for some $p\ne q$).

Using (\ref{chain5}), we get

\vspace{-1mm}
\begin{equation}
\label{chain30}
\sum_{(t_1,\ldots,t_k)}
\int\limits_{t}^{T}
\ldots
\int\limits_{t}^{t_2}
\Phi_N(t_1,\ldots,t_k)d{\bf w}_{t_1}^{(i_1)}
\ldots
d{\bf w}_{t_k}^{(i_k)}=
\end{equation}

\vspace{1mm}
$$
=\sum\limits_{(l_1,\ldots,l_k)}
\sum_{\stackrel{l_1,\ldots,l_k=0}{{}_{l_1<l_2<\ldots < l_k}}}^{N-1}
a_{l_1 \ldots l_k} 
\Delta{\bf w}_{\tau_{l_1}}^{(i_1)} \ldots \Delta{\bf w}_{\tau_{l_k}}^{(i_k)}=
$$

\vspace{1mm}
$$
=\sum\limits_{\stackrel{l_1,\ldots,l_k=0}{{}_{l_q\ne l_r;\ q\ne r;\ 
q, r=1,\ldots, k}}}^{N-1}
a_{l_1 \ldots l_k} 
\Delta{\bf w}_{\tau_{l_1}}^{(i_1)} \ldots \Delta{\bf w}_{\tau_{l_k}}^{(i_k)}=
$$

\vspace{1mm}
\begin{equation}
\label{chain10}
=J'[\Phi_N]_{T,t}^{(i_1\ldots i_k)}\ \ \ \hbox{w.\ p.\ 1},
\end{equation}

\vspace{3mm}
\noindent
where permutations $(t_1,\ldots,t_k)$ when summing are 
performed only in the values
$d{\bf w}_{t_1}^{(i_1)}
\ldots $
$d{\bf w}_{t_k}^{(i_k)}$ 
and permutations $(l_1,\ldots,l_k)$ when summing are 
performed only in the expression $l_1<l_2<\ldots < l_k.$
At the same time the indices near 
upper 
limits of integration in the iterated stochastic integrals in (\ref{chain30}) are changed 
correspondently and if $t_r$ swapped with $t_q$ in the  
permutation $(t_1,\ldots,t_k)$, then $i_r$ swapped with $i_q$ in 
the permutation $(i_1,\ldots,i_k)$ (see (\ref{chain30})).
In addition, the multiple Wiener stochastic integral 
$J'[\Phi_N]_{T,t}^{(i_1\ldots i_k)}$ is defined by (\ref{chain9})
and 

\vspace{1mm}
$$
\int\limits_{t}^{T}
\ldots
\int\limits_{t}^{t_2}
\Phi_N(t_1,\ldots,t_k)d{\bf w}_{t_1}^{(i_1)}
\ldots
d{\bf w}_{t_k}^{(i_k)}
$$

\vspace{4mm}
\noindent
is the iterated Ito stochastic integral.

Using (\ref{chain15}), (\ref{chain10}), Lemma 2, and (\ref{zx1}) for Lebesgue integrals, 
we have

\vspace{1mm}
$$
{\sf M}\left\{\left(J'[\Phi_N]_{T,t}^{(i_1\ldots i_k)}-
J'[\Phi_M]_{T,t}^{(i_1\ldots i_k)}\right)^2\right\}\le
$$

\vspace{2mm}
$$
\le C_k 
\sum_{(t_1,\ldots,t_k)}
\int\limits_{t}^{T}
\ldots
\int\limits_{t}^{t_2}
\left(\Phi_N(t_1,\ldots,t_k)-\Phi_M(t_1,\ldots,t_k)\right)^2 dt_1
\ldots dt_k=
$$

\vspace{2mm}
$$
=C_k 
\int\limits_{[t,T]^k}
\left(\Phi_N(t_1,\ldots,t_k)-\Phi_M(t_1,\ldots,t_k)\right)^2 dt_1
\ldots dt_k=
$$

\vspace{2mm}
$$
=C_k\left\Vert \Phi_N-\Phi_M\right\Vert_{L_2([t, T]^k)}^2\le
$$

\vspace{2mm}
$$
\le 2 C_k \left(\left\Vert \Phi_N-\Phi\right\Vert_{L_2([t, T]^k)}^2+
\left\Vert \Phi-\Phi_M\right\Vert_{L_2([t, T]^k)}^2\right)^2\ \to 0
$$

\vspace{5mm}
\noindent
if $N,M\to\infty,$ 
where constant $C_k$ 
depends only
on the multiplicity $k$ of the multiple Wiener stochastic integral.

Thus, there exists the limit 

$$
\hbox{\vtop{\offinterlineskip\halign{
\hfil#\hfil\cr
{\rm l.i.m.}\cr
$\stackrel{}{{}_{N\to \infty}}$\cr
}} }J'[\Phi_N]_{T,t}^{(i_1\ldots i_k)}.
$$

\vspace{4mm}

We will define the multiple Wiener stochastic integral for $\Phi(t_1,\ldots,t_k)\in L_2([t, T]^k)$ 
by the formula \cite{ito1951}, \cite{Kuo}
\begin{equation}
\label{WiI}
J'[\Phi]_{T,t}^{(i_1\ldots i_k)}\stackrel{\sf def}{=}
\hbox{\vtop{\offinterlineskip\halign{
\hfil#\hfil\cr
{\rm l.i.m.}\cr
$\stackrel{}{{}_{N\to \infty}}$\cr
}} }J'[\Phi_N]_{T,t}^{(i_1\ldots i_k)}=
\hbox{\vtop{\offinterlineskip\halign{
\hfil#\hfil\cr
{\rm l.i.m.}\cr
$\stackrel{}{{}_{N\to \infty}}$\cr
}} }
\sum\limits_{l_1,\ldots,l_k=0}^{N-1}
a_{l_1 \ldots l_k}
\Delta{\bf w}_{\tau_{l_1}}^{(i_1)}\ldots \Delta{\bf w}_{\tau_{l_k}}^{(i_k)},
\end{equation}

\vspace{4mm}
\noindent
where $\Phi_N(t_1,\ldots,t_k)$ is defined by 
(\ref{chain3}),
$\Delta{\bf w}_{\tau_{j}}^{(i)}=
{\bf w}_{\tau_{j+1}}^{(i)}-{\bf w}_{\tau_{j}}^{(i)},$\
$i=0, 1,\ldots,m,$\ ${\bf w}_{\tau}^{(0)}=\tau.$

Let us prove the following equality 

\begin{equation}
\label{Wi110}
J'[\Phi]_{T,t}^{(i_1\ldots i_k)}=\sum_{(t_1,\ldots,t_k)}
\int\limits_{t}^{T}
\ldots
\int\limits_{t}^{t_2}
\Phi(t_1,\ldots,t_k)d{\bf w}_{t_1}^{(i_1)}
\ldots
d{\bf w}_{t_k}^{(i_k)}\ \ \ \hbox{w.\ p.\ 1},
\end{equation}

\vspace{4mm}
\noindent
where permutations $(t_1,\ldots,t_k)$ when summing are 
performed only in the values
$d{\bf w}_{t_1}^{(i_1)}
\ldots $
$d{\bf w}_{t_k}^{(i_k)}.$ At the same time the indices near 
upper 
limits of integration in the iterated stochastic integrals are changed 
correspondently and if $t_r$ swapped with $t_q$ in the  
permutation $(t_1,\ldots,t_k)$, then $i_r$ swapped with $i_q$ in 
the permutation $(i_1,\ldots,i_k).$ 
In addition, the multiple Wiener stochastic integral 
$J'[\Phi]_{T,t}^{(i_1\ldots i_k)}$ is defined by (\ref{WiI})
and 

$$
\int\limits_{t}^{T}
\ldots
\int\limits_{t}^{t_2}
\Phi(t_1,\ldots,t_k)d{\bf w}_{t_1}^{(i_1)}
\ldots
d{\bf w}_{t_k}^{(i_k)}
$$

\vspace{4mm}
\noindent
is the iterated Ito stochastic integral.

The equality (\ref{Wi110}) has already been proved for the case 
$\Phi(t_1,\ldots,t_k)=\Phi_N(t_1,\ldots,t_k)$ (see (\ref{chain10})).
From (\ref{chain10}) we have
$$
J'[\Phi_N]_{T,t}^{(i_1\ldots i_k)}=
\sum_{(t_1,\ldots,t_k)}
\int\limits_{t}^{T}
\ldots
\int\limits_{t}^{t_2}
\Phi_N(t_1,\ldots,t_k)d{\bf w}_{t_1}^{(i_1)}
\ldots
d{\bf w}_{t_k}^{(i_k)}=
$$

\vspace{2mm}
$$
=\sum_{(t_1,\ldots,t_k)}
\int\limits_{t}^{T}
\ldots
\int\limits_{t}^{t_2}
\Phi(t_1,\ldots,t_k)d{\bf w}_{t_1}^{(i_1)}
\ldots
d{\bf w}_{t_k}^{(i_k)}+
$$

\vspace{2mm}
\begin{equation}
\label{chain11}
+\sum_{(t_1,\ldots,t_k)}
\int\limits_{t}^{T}
\ldots
\int\limits_{t}^{t_2}
\left(\Phi_N(t_1,\ldots,t_k)-\Phi(t_1,\ldots,t_k)\right)d{\bf w}_{t_1}^{(i_1)}
\ldots
d{\bf w}_{t_k}^{(i_k)}\ \ \ \hbox{w.~p.~1.}
\end{equation}

\vspace{5mm}

Passing to the limit $\hbox{\vtop{\offinterlineskip\halign{
\hfil#\hfil\cr
{\rm l.i.m.}\cr
$\stackrel{}{{}_{N\to \infty}}$\cr
}} }$ in the equality (\ref{chain11}), we obtain

\vspace{1mm}
$$
J'[\Phi]_{T,t}^{(i_1\ldots i_k)}=
\sum_{(t_1,\ldots,t_k)}
\int\limits_{t}^{T}
\ldots
\int\limits_{t}^{t_2}
\Phi(t_1,\ldots,t_k)d{\bf w}_{t_1}^{(i_1)}
\ldots
d{\bf w}_{t_k}^{(i_k)}+
$$

\vspace{2mm}
\begin{equation}
\label{chain12}
+\hbox{\vtop{\offinterlineskip\halign{
\hfil#\hfil\cr
{\rm l.i.m.}\cr
$\stackrel{}{{}_{N\to \infty}}$\cr
}} }\sum_{(t_1,\ldots,t_k)}
\int\limits_{t}^{T}
\ldots
\int\limits_{t}^{t_2}
\left(\Phi_N(t_1,\ldots,t_k)-\Phi(t_1,\ldots,t_k)\right)d{\bf w}_{t_1}^{(i_1)}
\ldots
d{\bf w}_{t_k}^{(i_k)}\ \ \ \hbox{w.~p.~1.}
\end{equation}

\vspace{5mm}
             
Using Lemma 1.2 \cite{20a}-\cite{20a-new-x} or
Lemma 2 \cite{26a} as well as (\ref{zx1}) for Lebesgue integrals and (\ref{chain15}), we get

\vspace{1mm}
$$
{\sf M}\left\{\left(
\sum_{(t_1,\ldots,t_k)}
\int\limits_{t}^{T}
\ldots
\int\limits_{t}^{t_2}
\left(\Phi_N(t_1,\ldots,t_k)-\Phi(t_1,\ldots,t_k)\right)d{\bf w}_{t_1}^{(i_1)}
\ldots
d{\bf w}_{t_k}^{(i_k)}\right)^2\right\}\le
$$

\vspace{2mm}
$$
\le C_k 
\sum_{(t_1,\ldots,t_k)}
\int\limits_{t}^{T}
\ldots
\int\limits_{t}^{t_2}
\left(\Phi_N(t_1,\ldots,t_k)-\Phi(t_1,\ldots,t_k)\right)^2 dt_1
\ldots dt_k=
$$

\vspace{2mm}
\begin{equation}
\label{chain20}
=C_k 
\int\limits_{[t,T]^k}
\left(\Phi_N(t_1,\ldots,t_k)-\Phi(t_1,\ldots,t_k)\right)^2 dt_1
\ldots dt_k\ \to 0
\end{equation}

\vspace{5mm}
\noindent
if $N\to\infty,$ 
where constant $C_k$ 
depends only
on the multiplicity $k$ of the multiple Wiener stochastic integral.

The relations (\ref{chain12}) and (\ref{chain20}) prove the equality 
(\ref{Wi110}).
From (\ref{Wi110}) we have

\begin{equation}
\label{wi1001}
J[\psi^{(k)}]_{T,t}^{(i_1\ldots i_k)}=\int\limits_t^T\psi_k(t_k) \ldots \int\limits_t^{t_{2}}
\psi_1(t_1) d{\bf w}_{t_1}^{(i_1)}\ldots
d{\bf w}_{t_k}^{(i_k)}=J'[K]_{T,t}^{(i_1\ldots i_k)}\ \ \ \hbox{w.\ p.\ 1},
\end{equation}

\vspace{3mm}
\noindent
where 
$K=K(t_1,\ldots,t_k)$ is defined by (\ref{ppp}).

Applying (\ref{wi1001}), we obtain

\vspace{1mm}
\begin{equation}
\label{chain102}
J[\psi^{(k)}]_{T,t}^{(i_1\ldots i_k)}=J'[K]_{T,t}^{(i_1\ldots i_k)}
=\sum_{j_1=0}^{p_1}\ldots
\sum_{j_k=0}^{p_k}
C_{j_k\ldots j_1}
J'[\phi_{j_1}\ldots \phi_{j_k}]_{T,t}^{(i_1\ldots i_k)}+
J'[R_{p_1\ldots p_k}]_{T,t}^{(i_1\ldots i_k)}
\end{equation}

\vspace{4mm}
\noindent
w.~p.~1, where
\begin{equation}
\label{chain30001}
R_{p_1\ldots p_k}(t_1,\ldots,t_k)\stackrel{{\rm def}}{=}
K(t_1,\ldots,t_k)-
\sum_{j_1=0}^{p_1}\ldots
\sum_{j_k=0}^{p_k}
C_{j_k\ldots j_1}
\prod_{l=1}^k\phi_{j_l}(t_l)
\end{equation}

\vspace{3mm}
\noindent
and
\begin{equation}
\label{chain300}
C_{j_k\ldots j_1}=\int\limits_{[t,T]^k}
K(t_1,\ldots,t_k)\prod_{l=1}^{k}\phi_{j_l}(t_l)dt_1\ldots dt_k
\end{equation}

\vspace{4mm}
\noindent
is the Fourier coefficient corresponding to $K(t_1,\ldots,t_k).$

Again applying (\ref{Wi110}), we have

\begin{equation}
\label{wi2005}
J'[R_{p_1\ldots p_k}]_{T,t}^{(i_1\ldots i_k)}
=
\sum_{(t_1,\ldots,t_k)}
\int\limits_{t}^{T}
\ldots
\int\limits_{t}^{t_2}
\Biggl(K(t_1,\ldots,t_k)-
\sum_{j_1=0}^{p_1}\ldots
\sum_{j_k=0}^{p_k}
C_{j_k\ldots j_1}
\prod_{l=1}^k\phi_{j_l}(t_l)\Biggr)
d{\bf w}_{t_1}^{(i_1)}
\ldots
d{\bf w}_{t_k}^{(i_k)},
\end{equation}

\vspace{4mm}
\noindent
where permutations $(t_1,\ldots,t_k)$ when summing are performed only 
in the values $d{\bf w}_{t_1}^{(i_1)}
\ldots $
$d{\bf w}_{t_k}^{(i_k)}$. At the same time the indices near 
upper limits of integration in the iterated stochastic integrals 
are changed correspondently and if $t_r$ swapped with $t_q$ in the  
permutation $(t_1,\ldots,t_k)$, then $i_r$ swapped with $i_q$ in the 
permutation $(i_1,\ldots,i_k).$
In addition, the multiple Wiener stochastic integral
$J'[R_{p_1\ldots p_k}]_{T,t}^{(i_1\ldots i_k)}$ is defined by 
(\ref{WiI}).

Using Lemma 1.2 \cite{20a}-\cite{20a-new-x} or
Lemma 2 \cite{26a}, (\ref{sos1z}) as well as 
(\ref{zx1}) for Lebesgue integrals, we have

\vspace{1mm}
$$
{\sf M}\left\{\left(J'[R_{p_1\ldots p_k}]_{T,t}^{(i_1\ldots i_k)}\right)^2\right\}
\le 
$$

\vspace{2mm}
$$
\le C_k
\sum_{(t_1,\ldots,t_k)}
\int\limits_{t}^{T}
\ldots
\int\limits_{t}^{t_2}
\left(K(t_1,\ldots,t_k)-
\sum_{j_1=0}^{p_1}\ldots
\sum_{j_k=0}^{p_k}
C_{j_k\ldots j_1}
\prod_{l=1}^k\phi_{j_l}(t_l)\right)^2
dt_1
\ldots
dt_k=
$$

\vspace{2mm}
\begin{equation}
\label{chain7771}
=C_k\int\limits_{[t,T]^k}
\left(K(t_1,\ldots,t_k)-
\sum_{j_1=0}^{p_1}\ldots
\sum_{j_k=0}^{p_k}
C_{j_k\ldots j_1}
\prod_{l=1}^k\phi_{j_l}(t_l)\right)^2
dt_1
\ldots
dt_k\to 0
\end{equation}

\vspace{4mm}
\noindent
if $p_1,\ldots,p_k\to\infty,$ where constant $C_k$ 
depends only
on the multiplicity $k$ of the 
iterated  Ito stochastic integral
$J[\psi^{(k)}]_{T,t}^{(i_1\ldots i_k)}$.

Thus, the following theorem is proved.

\vspace{2mm}

{\bf Theorem~4}\ \cite{20a} (Sect.~1.11), \cite{26a} (Sect.~15).
{\it Suppose that
$\psi_1(\tau),\ldots,\psi_k(\tau)\in L_2([t, T])$ and
$\{\phi_j(x)\}_{j=0}^{\infty}$ is an arbitrary complete orthonormal system  
of functions in the space $L_2([t,T]).$
Then the following expansion

\vspace{1mm}
$$
J[\psi^{(k)}]_{T,t}=
\hbox{\vtop{\offinterlineskip\halign{
\hfil#\hfil\cr
{\rm l.i.m.}\cr
$\stackrel{}{{}_{p_1,\ldots,p_k\to \infty}}$\cr
}} }
\sum\limits_{j_1=0}^{p_1}\ldots
\sum\limits_{j_k=0}^{p_k}
C_{j_k\ldots j_1}\Biggl(
\prod_{l=1}^k\zeta_{j_l}^{(i_l)}+\sum\limits_{r=1}^{[k/2]}
(-1)^r \times
\Biggr.
$$

\vspace{3mm}
$$
\times
\sum_{\stackrel{(\{\{g_1, g_2\}, \ldots, 
\{g_{2r-1}, g_{2r}\}\}, \{q_1, \ldots, q_{k-2r}\})}
{{}_{\{g_1, g_2, \ldots, 
g_{2r-1}, g_{2r}, q_1, \ldots, q_{k-2r}\}=\{1, 2, \ldots, k\}}}}
\prod\limits_{s=1}^r
{\bf 1}_{\{i_{g_{{}_{2s-1}}}=~i_{g_{{}_{2s}}}\ne 0\}}
\Biggl.{\bf 1}_{\{j_{g_{{}_{2s-1}}}=~j_{g_{{}_{2s}}}\}}
\prod_{l=1}^{k-2r}\zeta_{j_{q_l}}^{(i_{q_l})}\Biggr)
$$

\vspace{5mm}
\noindent
con\-verg\-ing in the mean-square sense is valid,
where $[x]$ is an integer part of a real number $x;$
another notations are the same as in Theorems~{\rm 2, 3}.}

\vspace{2mm}

It should be noted that an analogue of Theorem 4 was considered 
in \cite{Rybakov1000}. 
Note that we use another notations 
\cite{20a} (Sect.~1.11), \cite{26a} (Sect.~15)
in comparison with \cite{Rybakov1000}.
Moreover, the proof of an analogue of Theorem 4
from \cite{Rybakov1000} is different from the proof given in 
\cite{20a} (Sect.~1.11), \cite{26a} (Sect.~15).

\vspace{5mm}

\section{Expansions of Iterated Stratonovich Stochastic Integrals of Multiplicities 1 to 6
Based on Multiple Fourier--Legendre Series and Multiple Trigonometric Fourier Series}

\vspace{5mm}

In a number of works of the author \cite{12}-\cite{16}, \cite{19},
\cite{20}-\cite{21}, \cite{25}
Theorems 2, 4 have been adapted for the iterated Stratonovich stochastic integrals
(\ref{str}) of multiplicities 2 to 6. Let us collect some old results
in the following statement.

\vspace{2mm}

{\bf Theorem 5} \cite{12}-\cite{16}, \cite{19}, \cite{20}-\cite{21}, 
\cite{25}.\ 
{\it Suppose that 
$\{\phi_j(x)\}_{j=0}^{\infty}$ is a complete orthonormal system of 
Legendre polynomials or trigonometric functions in the space $L_2([t, T]).$
At the same time $\psi_2(\tau)$ is a continuously differentiable 
function on $[t, T]$ and $\psi_1(\tau), \psi_3(\tau)$ are twice 
continuously differentiable functions on $[t, T]$. Then

\vspace{-1mm}
\begin{equation}
\label{a}
J^{*}[\psi^{(2)}]_{T,t}=
\hbox{\vtop{\offinterlineskip\halign{
\hfil#\hfil\cr
{\rm l.i.m.}\cr
$\stackrel{}{{}_{p_1,p_2\to \infty}}$\cr
}} }\sum_{j_1=0}^{p_1}\sum_{j_2=0}^{p_2}
C_{j_2j_1}\zeta_{j_1}^{(i_1)}\zeta_{j_2}^{(i_2)}\ \ \ (i_1,i_2=1,\ldots,m),
\end{equation}

\vspace{1mm}
\begin{equation}
\label{feto19000ab}
J^{*}[\psi^{(3)}]_{T,t}=
\hbox{\vtop{\offinterlineskip\halign{
\hfil#\hfil\cr
{\rm l.i.m.}\cr
$\stackrel{}{{}_{p_1,p_2,p_3\to \infty}}$\cr
}} }\sum_{j_1=0}^{p_1}\sum_{j_2=0}^{p_2}\sum_{j_3=0}^{p_3}
C_{j_3 j_2 j_1}\zeta_{j_1}^{(i_1)}\zeta_{j_2}^{(i_2)}\zeta_{j_3}^{(i_3)}\ \ \
(i_1,i_2,i_3=0, 1,\ldots,m),
\end{equation}

\vspace{1mm}
\begin{equation}
\label{feto19000a}
J^{*}[\psi^{(3)}]_{T,t}=
\hbox{\vtop{\offinterlineskip\halign{
\hfil#\hfil\cr
{\rm l.i.m.}\cr
$\stackrel{}{{}_{p\to \infty}}$\cr
}} }
\sum\limits_{j_1,j_2,j_3=0}^{p}
C_{j_3 j_2 j_1}\zeta_{j_1}^{(i_1)}\zeta_{j_2}^{(i_2)}\zeta_{j_3}^{(i_3)}\ \ \
(i_1,i_2,i_3=1,\ldots,m),
\end{equation}

\vspace{1mm}
\begin{equation}
\label{uu}
J^{*}[\psi^{(4)}]_{T,t}=
\hbox{\vtop{\offinterlineskip\halign{
\hfil#\hfil\cr
{\rm l.i.m.}\cr
$\stackrel{}{{}_{p\to \infty}}$\cr
}} }
\sum\limits_{j_1, \ldots, j_4=0}^{p}
C_{j_4 j_3 j_2 j_1}\zeta_{j_1}^{(i_1)}
\zeta_{j_2}^{(i_2)}\zeta_{j_3}^{(i_3)}\zeta_{j_4}^{(i_4)}\ \ \
(i_1,\ldots,i_4=0, 1,\ldots,m),
\end{equation}

\vspace{3mm}
\noindent
where $J^{*}[\psi^{(k)}]_{T,t}$ is defined by {\rm (\ref{str})} and
$\psi_l(\tau)\equiv 1$ $(l=1,\ldots,4)$ in {\rm (\ref{feto19000ab})}, 
{\rm (\ref{uu});} another notations are the same as in Theorems {\rm 2--4.}
}

Note that the formula (\ref{a}) is generalized to the case
of continuous functions $\psi_1(\tau), \psi_2(\tau)$
in \cite{20a} (Sect.~2.1.4).

Recently, a new approach to the expansion and mean-square 
approximation of iterated Stratonovich stochastic integrals has been obtained
\cite{20a} (Sect.~2.10--2.16), \cite{21} (Sect.~13--19), 
\cite{25} (Sect.~5--11), \cite{arxiv-11} (Sect.~7--13),
\cite{new-art-1-xxy} (Sect.~4--9).
Let us formulate four theorems that were obtained using this approach.

\vspace{2mm}

{\bf Theorem 6}\ \cite{20a}, \cite{21}, \cite{25}, \cite{arxiv-11}, \cite{new-art-1-xxy}.\
{\it Suppose 
that $\{\phi_j(x)\}_{j=0}^{\infty}$ is a complete orthonormal system of 
Legendre polynomials or trigonometric functions in the space $L_2([t, T]).$
Furthermore, let $\psi_1(\tau), \psi_2(\tau),$ $\psi_3(\tau)$ are continuously dif\-ferentiable 
nonrandom functions on $[t, T].$ 
Then, for the 
iterated Stra\-to\-no\-vich stochastic integral of third multiplicity

\vspace{-1mm}
$$
J^{*}[\psi^{(3)}]_{T,t}={\int\limits_t^{*}}^T\psi_3(t_3)
{\int\limits_t^{*}}^{t_3}\psi_2(t_2)
{\int\limits_t^{*}}^{t_2}\psi_1(t_1)
d{\bf w}_{t_1}^{(i_1)}
d{\bf w}_{t_2}^{(i_2)}d{\bf w}_{t_3}^{(i_3)}\ \ \ (i_1,i_2,i_3=0,1,\ldots,m)
$$

\vspace{3mm}
\noindent
the following 
relations

\vspace{-2mm}
\begin{equation}
\label{fin1}
J^{*}[\psi^{(3)}]_{T,t}
=\hbox{\vtop{\offinterlineskip\halign{
\hfil#\hfil\cr
{\rm l.i.m.}\cr
$\stackrel{}{{}_{p\to \infty}}$\cr
}} }
\sum\limits_{j_1, j_2, j_3=0}^{p}
C_{j_3 j_2 j_1}\zeta_{j_1}^{(i_1)}\zeta_{j_2}^{(i_2)}\zeta_{j_3}^{(i_3)},
\end{equation}

\vspace{1mm}
\begin{equation}
\label{fin2}
{\sf M}\left\{\left(
J^{*}[\psi^{(3)}]_{T,t}-
\sum\limits_{j_1, j_2, j_3=0}^{p}
C_{j_3 j_2 j_1}\zeta_{j_1}^{(i_1)}\zeta_{j_2}^{(i_2)}\zeta_{j_3}^{(i_3)}\right)^2\right\}
\le \frac{C}{p}
\end{equation}

\vspace{3mm}
\noindent
are fulfilled, where $i_1, i_2, i_3=0,1,\ldots,m$ in {\rm (\ref{fin1})} and 
$i_1, i_2, i_3=1,\ldots,m$ in {\rm (\ref{fin2})},
constant $C$ is independent of $p,$

\vspace{-1mm}
$$
C_{j_3 j_2 j_1}=\int\limits_t^T\psi_3(t_3)\phi_{j_3}(t_3)
\int\limits_t^{t_3}\psi_2(t_2)\phi_{j_2}(t_2)
\int\limits_t^{t_2}\psi_1(t_1)\phi_{j_1}(t_1)dt_1dt_2dt_3
$$

\vspace{3mm}
\noindent
and
$$
\zeta_{j}^{(i)}=
\int\limits_t^T \phi_{j}(\tau) d{\bf f}_{\tau}^{(i)}
$$ 

\vspace{2mm}
\noindent
are independent standard Gaussian random variables for various 
$i$ or $j$ {\rm (}in the case when $i\ne 0${\rm );} 
another notations are the same as in Theorems~{\rm 2--4}.}

\vspace{2mm}

{\bf Theorem 7}\ \cite{20a}, \cite{21}, \cite{25}, \cite{arxiv-11}, \cite{new-art-1-xxy}.\ 
{\it Let
$\{\phi_j(x)\}_{j=0}^{\infty}$ be a complete orthonormal system of 
Legendre polynomials or trigonometric functions in the space $L_2([t, T]).$
Furthermore, let $\psi_1(\tau), \ldots,$ $\psi_4(\tau)$ be continuously dif\-ferentiable 
nonrandom functions on $[t, T].$ 
Then, for the 
iterated Stra\-to\-no\-vich stochastic integral of fourth multiplicity

\vspace{-1mm}
\begin{equation}
\label{fin0}
J^{*}[\psi^{(4)}]_{T,t}={\int\limits_t^{*}}^T\psi_4(t_4)
{\int\limits_t^{*}}^{t_4}\psi_3(t_3)
{\int\limits_t^{*}}^{t_3}\psi_2(t_2)
{\int\limits_t^{*}}^{t_2}\psi_1(t_1)
d{\bf w}_{t_1}^{(i_1)}
d{\bf w}_{t_2}^{(i_2)}d{\bf w}_{t_3}^{(i_3)}d{\bf w}_{t_4}^{(i_4)}
\end{equation}

\vspace{3mm}
\noindent
the following 
relations

\vspace{-1mm}
\begin{equation}
\label{fin3}
J^{*}[\psi^{(4)}]_{T,t}
=\hbox{\vtop{\offinterlineskip\halign{
\hfil#\hfil\cr
{\rm l.i.m.}\cr
$\stackrel{}{{}_{p\to \infty}}$\cr
}} }
\sum\limits_{j_1, j_2, j_3,j_4=0}^{p}
C_{j_4j_3 j_2 j_1}\zeta_{j_1}^{(i_1)}\zeta_{j_2}^{(i_2)}\zeta_{j_3}^{(i_3)}\zeta_{j_4}^{(i_4)},
\end{equation}

\vspace{1mm}

\begin{equation}
\label{fin4}
{\sf M}\left\{\left(
J^{*}[\psi^{(4)}]_{T,t}-
\sum\limits_{j_1, j_2, j_3, j_4=0}^{p}
C_{j_4 j_3 j_2 j_1}\zeta_{j_1}^{(i_1)}\zeta_{j_2}^{(i_2)}\zeta_{j_3}^{(i_3)}
\zeta_{j_4}^{(i_4)}
\right)^2\right\}
\le \frac{C}{p^{1-\varepsilon}}
\end{equation}

\vspace{3mm}
\noindent
are fulfilled, where $i_1, \ldots , i_4=0,1,\ldots,m$ in {\rm (\ref{fin0}),} {\rm (\ref{fin3})} 
and $i_1, \ldots, i_4=1,\ldots,m$ in {\rm (\ref{fin4}),}
constant $C$ does not depend on $p,$
$\varepsilon$ is an arbitrary
small positive real number 
for the case of complete orthonormal system of 
Legendre polynomials in the space $L_2([t, T])$
and $\varepsilon=0$ for the case of
complete orthonormal system of 
trigonometric functions in the space $L_2([t, T]),$

\vspace{-1mm}
$$
C_{j_4 j_3 j_2 j_1}=
$$

$$
=
\int\limits_t^T\psi_4(t_4)\phi_{j_4}(t_4)
\int\limits_t^{t_4}\psi_3(t_3)\phi_{j_3}(t_3)
\int\limits_t^{t_3}\psi_2(t_2)\phi_{j_2}(t_2)
\int\limits_t^{t_2}\psi_1(t_1)\phi_{j_1}(t_1)dt_1dt_2dt_3dt_4;
$$

\vspace{3mm}
\noindent
another notations are the same as in Theorem~{\rm 6}.}

\vspace{2mm}

{\bf Theorem 8}\ \cite{20a}, \cite{21}, \cite{25}, \cite{arxiv-11}, \cite{new-art-1-xxy}.\
{\it Assume 
that $\{\phi_j(x)\}_{j=0}^{\infty}$ is a complete orthonormal system of 
Legendre polynomials or trigonometric functions in the space $L_2([t, T])$
and $\psi_1(\tau), \ldots,$ $\psi_5(\tau)$ are continuously dif\-ferentiable 
nonrandom functions on $[t, T].$ 
Then, for the 
iterated Stra\-to\-no\-vich stochastic integral of fifth multiplicity

\vspace{-1mm}
\begin{equation}
\label{fin7}
J^{*}[\psi^{(5)}]_{T,t}={\int\limits_t^{*}}^T\psi_5(t_5)
\ldots
{\int\limits_t^{*}}^{t_2}\psi_1(t_1)
d{\bf w}_{t_1}^{(i_1)}
\ldots d{\bf w}_{t_5}^{(i_5)}
\end{equation}

\vspace{3mm}
\noindent
the following 
relations

\vspace{-1mm}
\begin{equation}
\label{fin8}
J^{*}[\psi^{(5)}]_{T,t}
=\hbox{\vtop{\offinterlineskip\halign{
\hfil#\hfil\cr
{\rm l.i.m.}\cr
$\stackrel{}{{}_{p\to \infty}}$\cr
}} }
\sum\limits_{j_1,\ldots,j_5=0}^{p}
C_{j_5 \ldots j_1}\zeta_{j_1}^{(i_1)}\ldots \zeta_{j_5}^{(i_5)},
\end{equation}

\vspace{1mm}

\begin{equation}
\label{fin9}
{\sf M}\left\{\left(
J^{*}[\psi^{(5)}]_{T,t}-
\sum\limits_{j_1, \ldots, j_5=0}^{p}
C_{j_5 \ldots j_1}\zeta_{j_1}^{(i_1)}\ldots
\zeta_{j_5}^{(i_5)}
\right)^2\right\}
\le \frac{C}{p^{1-\varepsilon}}
\end{equation}

\vspace{3mm}
\noindent
are fulfilled, where $i_1, \ldots , i_5=0,1,\ldots,m$ in {\rm (\ref{fin7}),} {\rm (\ref{fin8})} 
and $i_1, \ldots, i_5=1,\ldots,m$ in {\rm (\ref{fin9}),}
constant $C$ is independent of $p,$
$\varepsilon$ is an arbitrary
small positive real number 
for the case of complete orthonormal system of 
Legendre polynomials in the space $L_2([t, T])$
and $\varepsilon=0$ for the case of
complete orthonormal system of 
trigonometric functions in the space $L_2([t, T]),$

\vspace{-1mm}
$$
C_{j_5 \ldots j_1}=
\int\limits_t^T\psi_5(t_5)\phi_{j_5}(t_5)\ldots
\int\limits_t^{t_2}\psi_1(t_1)\phi_{j_1}(t_1)dt_1\ldots dt_5;
$$

\vspace{3mm}
\noindent
another notations are the same as in Theorems~{\rm 6, 7}.}

\vspace{2mm}

{\bf Theorem 9}\ \cite{20a}, \cite{21}, \cite{25}, \cite{arxiv-11}.\
{\it Suppose that 
$\{\phi_j(x)\}_{j=0}^{\infty}$ is a complete orthonormal system of 
Legendre polynomials or trigonometric functions in the space $L_2([t, T]).$
Then, for the 
iterated Stratonovich stochastic integral of sixth multiplicity

$$
J_{T,t}^{*(i_1\ldots i_6)}={\int\limits_t^{*}}^T
\ldots
{\int\limits_t^{*}}^{t_2}
d{\bf w}_{t_1}^{(i_1)}
\ldots d{\bf w}_{t_6}^{(i_6)}
$$

\vspace{3mm}
\noindent
the following 
expansion 

\vspace{-3mm}
$$
J_{T,t}^{*(i_1\ldots i_6)}
=\hbox{\vtop{\offinterlineskip\halign{
\hfil#\hfil\cr
{\rm l.i.m.}\cr
$\stackrel{}{{}_{p\to \infty}}$\cr
}} }
\sum\limits_{j_1, \ldots, j_6=0}^{p}
C_{j_6 \ldots j_1}\zeta_{j_1}^{(i_1)}\ldots
\zeta_{j_6}^{(i_6)}
$$

\vspace{4mm}
\noindent
that converges in the mean-square sense is valid, where
$i_1, \ldots, i_6=0, 1,\ldots,m,$

\vspace{-1mm}
$$
C_{j_6 \ldots j_1}=
\int\limits_t^T\phi_{j_6}(t_6)\ldots
\int\limits_t^{t_2}\phi_{j_1}(t_1)dt_1\ldots dt_6;
$$

\vspace{3mm}
\noindent
another notations are the same as in Theorems~{\rm 6--8}.}

\vspace{2mm}

The results of this section were developed in \cite{20a} (Chapter~2), \cite{2024xx}-\cite{2025xxxaaa}.
In particular, analogues of Theorem~9 for iterated Stratonovich stochastic
integrals of multiplicities 7 and 8 were obtained in \cite{20a} (Sect.~2.36, 2.37).
In addition, the variants of Thorems 5--9 
were obtained
for the case when $\{\phi_j(x)\}_{j=0}^{\infty}$ is an arbitrary complete orthonormal system
of functions in $L_2([t, T])$ \cite{20a} (Sect.~2.1.4, 2.23, 2.24, 2.31--2.34),
\cite{2024xx}-\cite{2025xxxaaa}.

\vspace{5mm}

\section{Exact and Approximate Expressions
for the Mean-Square Error of Approximation of Iterated Ito
Stochastic Integrals in Theorems 2, 4}

\vspace{5mm}

{\bf Theorem 10}\ \cite{17}-\cite{19}, \cite{20}-\cite{20a-new-x}, \cite{26}.\ 
{\it Suppose that
every $\psi_l(\tau)$ $(l=1,\ldots, k)$ is a continuous function on 
$[t, T]$ and
$\{\phi_j(x)\}_{j=0}^{\infty}$ is a complete orthonormal system  
of continuous functions in the space $L_2([t,T]).$ Then

$$
{\sf M}\left\{\left(J[\psi^{(k)}]_{T,t}-
J[\psi^{(k)}]_{T,t}^p\right)^2\right\}
= \int\limits_{[t,T]^k} K^2(t_1,\ldots,t_k)
dt_1\ldots dt_k - 
$$

\vspace{1mm}
$$
- \sum_{j_1=0}^{p}\ldots\sum_{j_k=0}^{p}
C_{j_k\ldots j_1}
{\sf M}\left\{J[\psi^{(k)}]_{T,t}
\sum\limits_{(j_1,\ldots,j_k)}
\int\limits_t^T \phi_{j_k}(t_k)
\ldots
\int\limits_t^{t_{2}}\phi_{j_{1}}(t_{1})
d{\bf f}_{t_1}^{(i_1)}\ldots
d{\bf f}_{t_k}^{(i_k)}\right\},
$$

\vspace{6mm}
\noindent
where
$$
J[\psi^{(k)}]_{T,t}=\int\limits_t^T\psi_k(t_k) \ldots \int\limits_t^{t_{2}}
\psi_1(t_1) d{\bf f}_{t_1}^{(i_1)}\ldots
d{\bf f}_{t_k}^{(i_k)}\ \ \ (i_1,\ldots,i_k=1,\ldots,m),
$$

\vspace{2mm}
\begin{equation}
\label{yeee2}
J[\psi^{(k)}]_{T,t}^p=
\sum_{j_1=0}^{p}\ldots\sum_{j_k=0}^{p}
C_{j_k\ldots j_1}\left(
\prod_{l=1}^k\zeta_{j_l}^{(i_l)}-S_{j_1,\ldots,j_k}^{(i_1\ldots i_k)}
\right),
\end{equation}

\vspace{4mm}
\begin{equation}
\label{ppp1}
S_{j_1,\ldots,j_k}^{(i_1\ldots i_k)}=
\hbox{\vtop{\offinterlineskip\halign{
\hfil#\hfil\cr
{\rm l.i.m.}\cr
$\stackrel{}{{}_{N\to \infty}}$\cr
}} }\sum_{(l_1,\ldots,l_k)\in {\rm G}_k}
\phi_{j_{1}}(\tau_{l_1})
\Delta{\bf f}_{\tau_{l_1}}^{(i_1)}\ldots
\phi_{j_{k}}(\tau_{l_k})
\Delta{\bf f}_{\tau_{l_k}}^{(i_k)},
\end{equation}

\vspace{8mm}
\noindent
the Fourier coefficient $C_{j_k\ldots j_1}$ has the form {\rm (\ref{ppppa})},

\begin{equation}
\label{rr232}
\zeta_{j}^{(i)}=
\int\limits_t^T \phi_{j}(\tau) d{\bf f}_{\tau}^{(i)}
\end{equation}

\vspace{3mm}
\noindent
are independent standard Gaussian random variables
for various
$i$ or $j$ $(i=1,\ldots,m),$

$$
\sum\limits_{(j_1,\ldots,j_k)}
$$ 

\vspace{3mm}
\noindent
means the sum with respect to all
possible permutations 
$(j_1,\ldots,j_k)$ {\rm (}at that, if 
$j_r$ swapped with $j_q$ in the permutation $(j_1,\ldots,j_k)$,
then $i_r$ swapped with $i_q$ in the permutation
$(i_1,\ldots,i_k));$
another notations are the same as in Theorem {\rm 2.}}

\vspace{2mm}

{\bf Proof.}
Using Theorem 2 for the case $p_1=\ldots=p_k=p$ and
$i_1,\ldots,i_k=1,\ldots,m$, we obtain

\vspace{1mm}
\begin{equation}
\label{yyye1}
J[\psi^{(k)}]_{T,t}=\
\hbox{\vtop{\offinterlineskip\halign{
\hfil#\hfil\cr
{\rm l.i.m.}\cr
$\stackrel{}{{}_{p\to \infty}}$\cr
}} }\sum_{j_1=0}^{p}\ldots\sum_{j_k=0}^{p}
C_{j_k\ldots j_1}\left(
\prod_{l=1}^k\zeta_{j_l}^{(i_l)}-S_{j_1,\ldots,j_k}^{(i_1\ldots i_k)}
\right).
\end{equation}

\vspace{5mm}

For $n>p$ we can write 

\vspace{2mm}
$$
J[\psi^{(k)}]_{T,t}^n=
\left(\sum_{j_1=0}^{p}+\sum_{j_1=p+1}^n\right)\ldots
\left(\sum_{j_k=0}^{p}+\sum_{j_k=p+1}^n\right)
C_{j_k\ldots j_1}\left(
\prod_{l=1}^k\zeta_{j_l}^{(i_l)}-S_{j_1,\ldots,j_k}^{(i_1\ldots i_k)}
\right)=
$$

\vspace{4mm}
\begin{equation}
\label{yyye}
=J[\psi^{(k)}]_{T,t}^p + \xi[\psi^{(k)}]_{T,t}^{p+1,n}.
\end{equation}

\vspace{5mm}

Let us prove 
that 
due to the special structure of random variables 
$S_{j_1,\ldots,j_k}^{(i_1\ldots i_k)}$ (also see (\ref{a2})--(\ref{a6}))
the following relations 
are correct

\vspace{1mm}
\begin{equation}
\label{tyty}
{\sf M}\left\{
\prod_{l=1}^k\zeta_{j_l}^{(i_l)}-S_{j_1,\ldots,j_k}^{(i_1\ldots i_k)}
\right\}=0,
\end{equation}

\vspace{2mm}
\begin{equation}
\label{tyty1}
{\sf M}\left\{
\left(\prod_{l=1}^k\zeta_{j_l}^{(i_l)}-S_{j_1,\ldots,j_k}^{(i_1\ldots i_k)}
\right)
\left(\prod_{l=1}^k\zeta_{j_l'}^{(i_l)}-S_{j_1',\ldots,j_k'}^{(i_1\ldots i_k)}
\right)\right\}=0,
\end{equation}

\vspace{3mm}
\noindent
where
$$
(j_1,\ldots,j_k)\in{\rm K}_p,\ \ \ (j_1',\ldots,j_k')
\in{\rm K}_n\backslash {\rm K}_{p}
$$ 

\vspace{2mm}
and
$$
{\rm K}_n=\left\{(j_1,\ldots,j_k):\ 0\le j_1,\ldots,j_k\le n\right\},
$$

$$
{\rm K}_p=\left\{(j_1,\ldots,j_k):\ 0\le j_1,\ldots,j_k\le p\right\}.
$$

\vspace{6mm}

Let us prove (\ref{tyty}).
For the case $i_1,\ldots,i_k=1,\ldots,m$ 
and $p_1=\ldots=p_k=p$ from (\ref{e1}) and (\ref{e2}) we obtain

\vspace{2mm}
$$
\prod_{l=1}^k\zeta_{j_l}^{(i_l)}-
S_{j_1,\ldots,j_k}^{(i_1\ldots i_k)}\ \ =\ \ 
\hbox{\vtop{\offinterlineskip\halign{
\hfil#\hfil\cr
{\rm l.i.m.}\cr
$\stackrel{}{{}_{N\to \infty}}$\cr
}} }
\sum_{\stackrel{l_1,\ldots,l_k=0}
{{}_{l_q\ne l_r;\ q\ne r;\
q,r=1,\ldots,k}}}^{N-1}
\phi_{j_1}(\tau_{l_1})\ldots
\phi_{j_k}(\tau_{l_k})
\Delta{\bf f}_{\tau_{l_1}}^{(i_1)}
\ldots
\Delta{\bf f}_{\tau_{l_k}}^{(i_k)}=
$$

\vspace{2mm}
\begin{equation}
\label{ttt211}
=\sum\limits_{(j_1,\ldots,j_k)}
\int\limits_t^T \phi_{j_k}(t_k)
\ldots
\int\limits_t^{t_{2}}\phi_{j_{1}}(t_{1})
d{\bf f}_{t_1}^{(i_1)}\ldots
d{\bf f}_{t_k}^{(i_k)}\ \ \ {\rm w.\ p.\ 1,}
\end{equation}

\vspace{3mm}
\noindent
where 
$$
\sum\limits_{(j_1,\ldots,j_k)}
$$ 

\vspace{3mm}
\noindent
means the sum with respect to all
possible permutations
$(j_1,\ldots,j_k).$ At the same time if 
$j_r$ swapped  with $j_q$ in the permutation $(j_1,\ldots,j_k)$,
then $i_r$ swapped  with $i_q$ in the permutation
$(i_1,\ldots,i_k).$

From (\ref{ttt211}) due to the moment property of Ito
stochastic integral we obtain (\ref{tyty}).

Let us prove (\ref{tyty1}).
From (\ref{ttt211}) we have

\vspace{2mm}
$$
0\le \left\vert{\sf M}\left\{
\left(\prod_{l=1}^k\zeta_{j_l}^{(i_l)}-S_{j_1,\ldots,j_k}^{(i_1\ldots i_k)}
\right)
\left(\prod_{l=1}^k\zeta_{j_l'}^{(i_l)}-S_{j_1',\ldots,j_k'}^{(i_1\ldots i_k)}
\right)\right\}\right\vert=
$$

\vspace{5mm}
$$
=\left\vert
{\sf M}\left\{\sum\limits_{(j_1,\ldots,j_k)}\sum\limits_{(j_1',\ldots,j_k')}
\int\limits_t^T \phi_{j_k}(t_k)
\ldots
\int\limits_t^{t_{2}}\phi_{j_{1}}(t_{1})
d{\bf f}_{t_1}^{(i_1)}\ldots
d{\bf f}_{t_k}^{(i_k)}\times\right.\right.
$$

\vspace{3mm}
$$
\left.\left.\times
\int\limits_t^T \phi_{j_k'}(t_k)
\ldots
\int\limits_t^{t_{2}}\phi_{j_1'}(t_{1})
d{\bf f}_{t_1}^{(i_1)}\ldots
d{\bf f}_{t_k}^{(i_k)}\right\}\right\vert\le
$$

\vspace{3mm}

$$
\le\sum\limits_{(j_1',\ldots,j_k')}\ 
\int\limits_t^T \phi_{j_k}(t_k)\phi_{j_k'}(t_k)dt_k
\ldots
\int\limits_t^{T}\phi_{j_1}(t_{1})\phi_{j_1'}(t_{1})
dt_1=
$$

\vspace{3mm}
\begin{equation}
\label{hq11}
=
\sum\limits_{(j_1',\ldots,j_k')}
{\bf 1}_{\{j_1=j_1'\}}\ldots {\bf 1}_{\{j_k=j_k'\}},
\end{equation}

\vspace{5mm}
\noindent
where 
${\bf 1}_A$ is the indicator of the set $A$.
From (\ref{hq11}) we obtain (\ref{tyty1}).

Consider in detail the case $k=3$ in (\ref{hq11}).
We have

\vspace{1mm}
$$
\left\vert{\sf M}\left\{\sum\limits_{(j_1,j_2,j_3)}\sum\limits_{(j_1',j_2',j_3')}
\int\limits_t^T \phi_{j_3}(t_3)
\int\limits_t^{t_{3}}\phi_{j_{2}}(t_{2})
\int\limits_t^{t_{2}}\phi_{j_{1}}(t_{1})
d{\bf f}_{t_1}^{(i_1)}d{\bf f}_{t_2}^{(i_2)}
d{\bf f}_{t_3}^{(i_3)}\times\right.\right.
$$

\vspace{3mm}
$$
\left.\left.\times
\int\limits_t^T \phi_{j_3'}(t_3)
\int\limits_t^{t_{3}}\phi_{j_{2}'}(t_{2})
\int\limits_t^{t_{2}}\phi_{j_{1}'}(t_{1})
d{\bf f}_{t_1}^{(i_1)}d{\bf f}_{t_2}^{(i_2)}
d{\bf f}_{t_3}^{(i_3)}\right\}\right\vert=
$$

\vspace{2mm}

$$
=\left\vert\int\limits_t^T \phi_{j_3}(s)\phi_{j_3'}(s)ds
\int\limits_t^{T}\phi_{j_2}(s)\phi_{j_2'}(s)ds
\int\limits_t^{T}\phi_{j_1}(s)\phi_{j_1'}(s)ds+\right.
$$

$$
+{\bf 1}_{\{i_1=i_2\}}\int\limits_t^T \phi_{j_3}(s)\phi_{j_3'}(s)ds
\int\limits_t^{T}\phi_{j_1}(s)\phi_{j_2'}(s)ds
\int\limits_t^{T}\phi_{j_2}(s)\phi_{j_1'}(s)ds+
$$

$$
+{\bf 1}_{\{i_2=i_3\}}\int\limits_t^T \phi_{j_1}(s)\phi_{j_1'}(s)ds
\int\limits_t^{T}\phi_{j_2}(s)\phi_{j_3'}(s)ds
\int\limits_t^{T}\phi_{j_3}(s)\phi_{j_2'}(s)ds+
$$

$$
+{\bf 1}_{\{i_1=i_3\}}\int\limits_t^T \phi_{j_1}(s)\phi_{j_3'}(s)ds
\int\limits_t^{T}\phi_{j_2}(s)\phi_{j_2'}(s)ds
\int\limits_t^{T}\phi_{j_3}(s)\phi_{j_1'}(s)ds+
$$

$$
+{\bf 1}_{\{i_1=i_2=i_3\}}\int\limits_t^T \phi_{j_2}(s)\phi_{j_3'}(s)ds
\int\limits_t^{T}\phi_{j_1}(s)\phi_{j_2'}(s)ds
\int\limits_t^{T}\phi_{j_3}(s)\phi_{j_1'}(s)ds+
$$

$$
\left.+{\bf 1}_{\{i_1=i_2=i_3\}}\int\limits_t^T \phi_{j_1}(s)\phi_{j_3'}(s)ds
\int\limits_t^{T}\phi_{j_3}(s)\phi_{j_2'}(s)ds
\int\limits_t^{T}\phi_{j_2}(s)\phi_{j_1'}(s)ds\right\vert=
$$

\vspace{2mm}

$$
=\biggl\vert{\bf 1}_{\{j_3=j_3'\}}{\bf 1}_{\{j_2=j_2'\}}{\bf 1}_{\{j_1=j_1'\}}+
{\bf 1}_{\{i_1=i_2\}} \cdot
{\bf 1}_{\{j_3=j_3'\}}{\bf 1}_{\{j_1=j_2'\}}{\bf 1}_{\{j_2=j_1'\}}+\biggr.
$$

$$
+{\bf 1}_{\{i_2=i_3\}} \cdot
{\bf 1}_{\{j_1=j_1'\}}{\bf 1}_{\{j_2=j_3'\}}{\bf 1}_{\{j_3=j_2'\}}+
{\bf 1}_{\{i_1=i_3\}} \cdot
{\bf 1}_{\{j_1=j_3'\}}{\bf 1}_{\{j_2=j_2'\}}{\bf 1}_{\{j_3=j_1'\}}+
$$

$$
\biggl.+{\bf 1}_{\{i_1=i_2=i_3\}}\cdot
{\bf 1}_{\{j_2=j_3'\}}{\bf 1}_{\{j_1=j_2'\}}{\bf 1}_{\{j_3=j_1'\}}+
{\bf 1}_{\{i_1=i_2=i_3\}}\cdot
{\bf 1}_{\{j_1=j_3'\}}{\bf 1}_{\{j_3=j_2'\}}{\bf 1}_{\{j_2=j_1'\}}
\biggr\vert\le
$$

\vspace{2mm}

$$
\le{\bf 1}_{\{j_3=j_3'\}}{\bf 1}_{\{j_2=j_2'\}}{\bf 1}_{\{j_1=j_1'\}}+
{\bf 1}_{\{j_3=j_3'\}}{\bf 1}_{\{j_1=j_2'\}}{\bf 1}_{\{j_2=j_1'\}}+
$$

$$
+
{\bf 1}_{\{j_1=j_1'\}}{\bf 1}_{\{j_2=j_3'\}}{\bf 1}_{\{j_3=j_2'\}}+
{\bf 1}_{\{j_1=j_3'\}}{\bf 1}_{\{j_2=j_2'\}}{\bf 1}_{\{j_3=j_1'\}}+
$$

$$
+
{\bf 1}_{\{j_2=j_3'\}}{\bf 1}_{\{j_1=j_2'\}}{\bf 1}_{\{j_3=j_1'\}}+
{\bf 1}_{\{j_1=j_3'\}}{\bf 1}_{\{j_3=j_2'\}}{\bf 1}_{\{j_2=j_1'\}}
=
$$

\vspace{2mm}
$$
=
\sum\limits_{(j_1',j_2',j_3')}
{\bf 1}_{\{j_1=j_1'\}}{\bf 1}_{\{j_2=j_2'\}}{\bf 1}_{\{j_3=j_3'\}},
$$

\vspace{3mm}
\noindent
where we used the relation

$$
\int\limits_t^T
\phi_{g}(s)\phi_{q}(s)ds={\bf 1}_{\{g=q\}},\ \ \ g,q=0, 1, 2\ldots
$$

\vspace{3mm}

Now consider the case of an arbitrary $k\in \mathbb{N}.$ We have

\vspace{1mm}
$$
{\sf M}\left\{\sum\limits_{(j_1,\ldots,j_k)}\sum\limits_{(j_1',\ldots,j_k')}
\int\limits_t^T \phi_{j_k}(t_k)
\ldots
\int\limits_t^{t_{2}}\phi_{j_{1}}(t_{1})
d{\bf f}_{t_1}^{(i_1)}\ldots
d{\bf f}_{t_k}^{(i_k)}\times\right.
$$

\vspace{3.5mm}
$$
\left.\times
\int\limits_t^T \phi_{j_k'}(t_k)
\ldots
\int\limits_t^{t_{2}}\phi_{j_1'}(t_{1})
d{\bf f}_{t_1}^{(i_1)}\ldots
d{\bf f}_{t_k}^{(i_k)}\right\}=
$$

\vspace{3.5mm}
$$
={\sf M}\left\{\sum\limits_{(j_1,\ldots,j_k)}\sum\limits_{(j_1',\ldots,j_k')}
\int\limits_t^T \phi_{j_k}(t_k)
\ldots
\int\limits_t^{t_{2}}\phi_{j_{1}}(t_{1})
d{\bf f}_{t_1}^{(i_1)}\ldots
d{\bf f}_{t_k}^{(i_k)}\times\right.
$$

\vspace{3.5mm}
$$
\left.\times
\int\limits_t^T \phi_{j_k'}(t_k)
\ldots
\int\limits_t^{t_{2}}\phi_{j_1'}(t_{1})
d{\bf f}_{t_1}^{(i_1')}\ldots
d{\bf f}_{t_k}^{(i_k')}\right\}=
$$

\vspace{3.5mm}
$$
=\sum\limits_{(j_1,\ldots,j_k)}\sum\limits_{(j_1',\ldots,j_k')}
{\bf 1}_{\{i_k=i_k'\}}\ldots {\bf 1}_{\{i_1=i_1'\}}\times
$$

\vspace{3.5mm}
$$
\times
\int\limits_t^T \phi_{j_k}(t_k)\phi_{j_k'}(t_k)
\ldots
\int\limits_t^{t_{2}}\phi_{j_1}(t_1)\phi_{j_1'}(t_1)
dt_1\ldots dt_k=
$$

\vspace{3.5mm}
$$
=\sum\limits_{(j_1',\ldots,j_k')}
{\bf 1}_{\{i_k=i_k'\}}\ldots {\bf 1}_{\{i_1=i_1'\}}\times
$$

\vspace{3.5mm}
$$
\times
\int\limits_t^T \phi_{j_k}(t_k)\phi_{j_k'}(t_k)dt_k
\ldots
\int\limits_t^T\phi_{j_1}(t_1)\phi_{j_1'}(t_1)dt_1=
$$

\vspace{3.5mm}
\begin{equation}
\label{wen100}
=\sum\limits_{(j_1',\ldots,j_k')}
{\bf 1}_{\{i_k=i_k'\}}\ldots {\bf 1}_{\{i_1=i_1'\}}
{\bf 1}_{\{j_k=j_k'\}}\ldots {\bf 1}_{\{j_1=j_1'\}},
\end{equation}

\vspace{6mm}
\noindent
where $(i_1',\ldots,i_k')=(i_1,\ldots,i_k).$
However, if 
$j_r$ swapped  with $j_q$ in the permutation $(j_1,\ldots,j_k)$,
then $i_r$ swapped  with $i_q$ in the permutation
$(i_1,\ldots,i_k)$ and 
if $j_r'$ swapped  with $j_q'$ in the permutation $(j_1',\ldots,j_k')$,
then $i_r'$ swapped  with $i_q'$ in the permutation
$(i_1',\ldots,i_k')$. From (\ref{wen100}) we obtain (\ref{hq11}).
The equality (\ref{tyty1}) is proved.

Note that the formula (\ref{tyty1}) (in the light of the results of 
\cite{20a} (Sect.~1.10, 1.11) can be 
interpreted as a consequence of the orthogonality of two random 
variables that are Hermite polynomials of vector random arguments.

From (\ref{tyty}) and (\ref{tyty1}) we obtain

\vspace{2mm}
$$
{\sf M}\left\{J[\psi^{(k)}]_{T,t}^p\xi[\psi^{(k)}]_{T,t}^{p+1,n}
\right\}=0.
$$

\vspace{5mm}

Due to (\ref{yeee2}), (\ref{yyye1}), and (\ref{yyye}) we can write 

\vspace{3mm}
$$
\xi[\psi^{(k)}]_{T,t}^{p+1,n}=J[\psi^{(k)}]_{T,t}^n-J[\psi^{(k)}]_{T,t}^p,
$$

\vspace{3mm}
$$
\hbox{\vtop{\offinterlineskip\halign{
\hfil#\hfil\cr
{\rm l.i.m.}\cr
$\stackrel{}{{}_{n\to \infty}}$\cr
}} }
\xi[\psi^{(k)}]_{T,t}^{p+1,n}=J[\psi^{(k)}]_{T,t}-J[\psi^{(k)}]_{T,t}^p
\stackrel{\rm def}{=}\xi[\psi^{(k)}]_{T,t}^{p+1}.
$$

\vspace{6mm}

We have

\vspace{-2mm}
$$
0\le \left|{\sf M}\left\{
\xi[\psi^{(k)}]_{T,t}^{p+1}J[\psi^{(k)}]_{T,t}^p\right\}\right|=
$$

\vspace{5mm}
$$
=
\left|{\sf M}\left\{\left(
\xi[\psi^{(k)}]_{T,t}^{p+1}-
\xi[\psi^{(k)}]_{T,t}^{p+1,n}+\xi[\psi^{(k)}]_{T,t}^{p+1,n}\right)
J[\psi^{(k)}]_{T,t}^p\right\}\right|=
$$

\vspace{5mm}
$$
\le 
\left|{\sf M}\left\{\left(
\xi[\psi^{(k)}]_{T,t}^{p+1}-
\xi[\psi^{(k)}]_{T,t}^{p+1,n}\right)
J[\psi^{(k)}]_{T,t}^p\right\}\right|+
\left|{\sf M}\left\{
\xi[\psi^{(k)}]_{T,t}^{p+1,n}J[\psi^{(k)}]_{T,t}^p\right\}\right|=
$$

\vspace{5mm}
$$
=\left|{\sf M}\left\{\left(
J[\psi^{(k)}]_{T,t}-
J[\psi^{(k)}]_{T,t}^{n}\right)
J[\psi^{(k)}]_{T,t}^p\right\}\right|\le
$$

\vspace{5mm}
$$
\le \sqrt{{\sf M}\left\{\left(
J[\psi^{(k)}]_{T,t}-
J[\psi^{(k)}]_{T,t}^{n}\right)^2\right\}}
\sqrt{{\sf M}\left\{\left(
J[\psi^{(k)}]_{T,t}^p\right)^2\right\}}\le
$$

\vspace{5mm}
$$
\le \sqrt{{\sf M}\left\{\left(
J[\psi^{(k)}]_{T,t}-
J[\psi^{(k)}]_{T,t}^{n}\right)^2\right\}}\times
$$

\vspace{4mm}
$$
\times
\left(\sqrt{{\sf M}\left\{\left(
J[\psi^{(k)}]_{T,t}^p - J[\psi^{(k)}]_{T,t}\right)^2\right\}}
+ \sqrt{{\sf M}\left\{\left(
J[\psi^{(k)}]_{T,t}\right)^2\right\}}\right)\le
$$

\vspace{4mm}
\begin{equation}
\label{rrre}
\le K \sqrt{{\sf M}\left\{\left(
J[\psi^{(k)}]_{T,t}-
J[\psi^{(k)}]_{T,t}^{n}\right)^2\right\}} \to 0\ \ \ {\rm if}\ \ \ n\to\infty,
\end{equation}

\vspace{7mm}
\noindent
where $K$ is a constant.

From (\ref{rrre}) it follows that

$$
{\sf M}\left\{
\xi[\psi^{(k)}]_{T,t}^{p+1}J[\psi^{(k)}]_{T,t}^p\right\}=0
$$

\noindent
or

$$
{\sf M}\left\{
\left(J[\psi^{(k)}]_{T,t}-J[\psi^{(k)}]_{T,t}^p\right)
J[\psi^{(k)}]_{T,t}^p\right\}=0.
$$

\vspace{5mm}

The last equality means that

\vspace{2mm}
\begin{equation}
\label{yyyw}
{\sf M}\left\{
J[\psi^{(k)}]_{T,t}J[\psi^{(k)}]_{T,t}^p\right\}=
{\sf M}\left\{
\left(J[\psi^{(k)}]_{T,t}^p\right)^2\right\}.
\end{equation}

\vspace{5mm}

Taking into account (\ref{yyyw}), we obtain

\vspace{2mm}
$$
{\sf M}\left\{\left(J[\psi^{(k)}]_{T,t}-
J[\psi^{(k)}]_{T,t}^p\right)^2\right\}=
{\sf M}\left\{\left(J[\psi^{(k)}]_{T,t}\right)^2\right\}+
$$

\vspace{4mm}
$$
+
{\sf M}\left\{\left(J[\psi^{(k)}]_{T,t}^p\right)^2\right\}
-2{\sf M}\left\{J[\psi^{(k)}]_{T,t}J[\psi^{(k)}]_{T,t}^p\right\}=
{\sf M}\left\{\left(J[\psi^{(k)}]_{T,t}\right)^2\right\}-
$$

\vspace{5mm}
$$
-
{\sf M}\left\{J[\psi^{(k)}]_{T,t}J[\psi^{(k)}]_{T,t}^p\right\}=
$$

\vspace{3mm}
\begin{equation}
\label{tttr}
=\int\limits_{[t,T]^k} K^2(t_1,\ldots,t_k)
dt_1\ldots dt_k - 
{\sf M}\left\{J[\psi^{(k)}]_{T,t}J[\psi^{(k)}]_{T,t}^p\right\}.
\end{equation}

\vspace{6mm}

Consider the value

\vspace{1mm}
$$
{\sf M}\left\{J[\psi^{(k)}]_{T,t}J[\psi^{(k)}]_{T,t}^p\right\}.
$$

\vspace{5mm}

From (\ref{yeee2}) and (\ref{ttt211}) we obtain

\vspace{2mm}
\begin{equation}
\label{z9}
J[\psi^{(k)}]_{T,t}^p=
\sum_{j_1=0}^{p}\ldots\sum_{j_k=0}^{p}
C_{j_k\ldots j_1}
\sum\limits_{(j_1,\ldots,j_k)}
\int\limits_t^T \phi_{j_k}(t_k)
\ldots
\int\limits_t^{t_{2}}\phi_{j_{1}}(t_{1})
d{\bf f}_{t_1}^{(i_1)}\ldots
d{\bf f}_{t_k}^{(i_k)}.
\end{equation}

\vspace{6mm}

After substituting (\ref{z9}) into (\ref{tttr}) we obtain (\ref{tttr11}).
Theorem 10 is proved.  

\vspace{2mm}

Let $J[\psi^{(k)}]_{T,t}^{p_1,\ldots,p_k}$ be the 
expression before passing to the limit
$\hbox{\vtop{\offinterlineskip\halign{
\hfil#\hfil\cr
{\rm l.i.m.}\cr
$\stackrel{}{{}_{p_1,\ldots,p_k\to \infty}}$\cr
}} }$ in (\ref{tyyy}). Denote

\vspace{3mm}
$$
{\sf M}\left\{\left(J[\psi^{(k)}]_{T,t}-
J[\psi^{(k)}]_{T,t}^{p_1,\ldots,p_k}\right)^2\right\}\stackrel{{\rm def}}
{=}E_k^{p_1,\ldots,p_k},
$$

\vspace{4mm}
$$
E_k^{p_1,\ldots,p_k}\stackrel{{\rm def}}{=}E_k^p\ \ \ \hbox{if}\ \ \
p_1=\ldots=p_k=p,
$$

\vspace{4mm}
$$
\left\Vert K\right\Vert_{L_2([t,T]^k)}^2=\int\limits_{[t,T]^k}
K^2(t_1,\ldots,t_k)dt_1\ldots dt_k\stackrel{{\rm def}}{=}I_k.
$$

\vspace{7mm}

In \cite{19}, \cite{20}-\cite{20a-new-x}, \cite{26} it was shown that 

\vspace{2mm}
\begin{equation}
\label{zzz0}
E_k^{p_1,\ldots,p_k}\le k!\left(I_k-\sum_{j_1=0}^{p_1}\ldots
\sum_{j_k=0}^{p_k}C^2_{j_k\ldots j_1}\right)
\end{equation}

\vspace{2mm}
\noindent
if 
$$
i_1,\ldots,i_k=1,\ldots,m\ \ \ (0<T-t<\infty)
$$
 
\vspace{2mm}
\noindent
or 
$$
i_1,\ldots,i_k=0, 1,\ldots,m\ \ \ (0<T-t<1).
$$

\vspace{7mm}

Moreover, 
in \cite{20a} (Sect.~1.1.9, 1.11, 1.12), \cite{26a} (Sect.~6, 15, 16)
it was shown that

\vspace{2mm}
$$
{\sf M}\left\{\left(J[\psi^{(k)}]_{T,t}-
J[\psi^{(k)}]_{T,t}^{p_1,\ldots,p_k}\right)^{2n}\right\}\le 
$$

\vspace{2mm}
$$
\le C_{n,k}
\left(I_k-\sum_{j_1=0}^{p_1}\ldots
\sum_{j_k=0}^{p_k}C^2_{j_k\ldots j_1}\right)^n,
$$

\vspace{4mm}
\noindent
where 
$$
C_{n,k}=(k!)^{n} (2n-1)^{nk}\ \ \ (n\in \mathbb{N}).
$$

\vspace{6mm}

{\bf Remark 1.}\ {\it Note that

\vspace{2mm}
$$
{\sf M}\left\{J[\psi^{(k)}]_{T,t}
\int\limits_t^T \phi_{j_k}(t_k)
\ldots
\int\limits_t^{t_{2}}\phi_{j_{1}}(t_{1})
d{\bf f}_{t_1}^{(i_1)}\ldots
d{\bf f}_{t_k}^{(i_k)}\right\}=
$$

\vspace{2mm}
$$
=\int\limits_t^T\psi_k(t_k) \phi_{j_k}(t_k)\ldots \int\limits_t^{t_{2}}
\psi_1(t_1)\phi_{j_1}(t_1) dt_1\ldots dt_k=
C_{j_k\ldots j_1}.
$$

\vspace{6mm}

Then from Theorem {\rm 10} for pairwise different $i_1,\ldots,i_k$ 
and for $i_1=\ldots=i_k$
we get

\vspace{2mm}
\begin{equation}
\label{zzz4}
E_k^p= I_k- \sum_{j_1,\ldots,j_k=0}^{p}
C_{j_k\ldots j_1}^2,
\end{equation}

\vspace{2mm}
$$
E_k^p= I_k - \sum_{j_1,\ldots,j_k=0}^{p}
C_{j_k\ldots j_1}\left(\sum\limits_{(j_1,\ldots,j_k)}
C_{j_k\ldots j_1}\right).
$$
}

\vspace{7mm}

Consider some examples of application of Theorem 10
$(i_1, \ldots ,i_5=1,\ldots,m)$

\vspace{2mm}
$$
E_2^p
=I_2
-\sum_{j_1,j_2=0}^p
C_{j_2j_1}^2-
\sum_{j_1,j_2=0}^p
C_{j_2j_1}C_{j_1j_2}\ \ \ (i_1=i_2),
$$

\vspace{2mm}
\begin{equation}
\label{zzz5}
E_3^p=I_3
-\sum_{j_3,j_2,j_1=0}^p C_{j_3j_2j_1}^2-
\sum_{j_3,j_2,j_1=0}^p C_{j_3j_1j_2}C_{j_3j_2j_1}\ \ \ (i_1=i_2\ne i_3),
\end{equation}

\vspace{2mm}
\begin{equation}
\label{zzz6}
E_3^p=I_3-
\sum_{j_3,j_2,j_1=0}^p C_{j_3j_2j_1}^2-
\sum_{j_3,j_2,j_1=0}^p C_{j_2j_3j_1}C_{j_3j_2j_1}\ \ \ (i_1\ne i_2=i_3),
\end{equation}

\vspace{2mm}
\begin{equation}
\label{zzz7}
E_3^p=I_3
-\sum_{j_3,j_2,j_1=0}^p C_{j_3j_2j_1}^2-
\sum_{j_3,j_2,j_1=0}^p C_{j_3j_2j_1}C_{j_1j_2j_3}\ \ \ (i_1=i_3\ne i_2),
\end{equation}

\vspace{2mm}
$$
E^p_4 = I_4 - \sum_{j_1,\ldots,j_4=0}^{p}
C_{j_4\ldots j_1}\Biggl(\sum\limits_{(j_1,j_2)}
C_{j_4\ldots j_1}\Biggr)\ \ \ (i_1=i_2\ne i_3, i_4;\ i_3\ne i_4),
$$

\vspace{2mm}
$$
E^p_4 = I_4 - \sum_{j_1,\ldots,j_4=0}^{p}
C_{j_4\ldots j_1}\Biggl(\sum\limits_{(j_1,j_4)}
C_{j_4\ldots j_1}\Biggr)\ \ \ (i_1=i_4\ne i_2, i_3;\ i_2\ne i_3),
$$

\vspace{3mm}
$$
E_4^p = I_4 -
\sum_{j_1,\ldots,j_4=0}^{p}
C_{j_4\ldots j_1}\Biggl(\sum\limits_{(j_1,j_2,j_3)}
C_{j_4\ldots j_1}\Biggr)\ \ \ (i_1=i_2=i_3\ne i_4),
$$

\vspace{3mm}
$$
E_4^p = I_4 -
 \sum_{j_1,\ldots,j_4=0}^{p}
C_{j_4\ldots j_1}\Biggl(\sum\limits_{(j_2,j_3,j_4)}
C_{j_4\ldots j_1}\Biggr)\ \ \ (i_2=i_3=i_4\ne i_1),
$$

\vspace{3mm}
$$
E^p_4 = I_4 - \sum_{j_1,\ldots,j_4=0}^{p}
C_{j_4\ldots j_1}\Biggl(\sum\limits_{(j_1,j_2)}\Biggl(
\sum\limits_{(j_3,j_4)}
C_{j_4\ldots j_1}\Biggr)\Biggr)\ \ \ (i_1=i_2\ne i_3=i_4),
$$

\vspace{3mm}
$$
E^p_4 = I_4 - \sum_{j_1,\ldots,j_4=0}^{p}
C_{j_4\ldots j_1}\Biggl(\sum\limits_{(j_1,j_4)}\Biggl(
\sum\limits_{(j_2,j_3)}
C_{j_4\ldots j_1}\Biggr)\Biggr)\ \ \ (i_1=i_4\ne i_2=i_3),
$$

\vspace{3mm}
$$
E_5^p = I_5 - \sum_{j_1,\ldots,j_5=0}^{p}
C_{j_5\ldots j_1}\Biggl(\sum\limits_{(j_2,j_4)}\Biggl(
\sum\limits_{(j_3,j_5)}
C_{j_5\ldots j_1}\Biggr)\Biggr)\ \ \ (i_1\ne i_2=i_4\ne i_3=i_5\ne i_1),
$$

\vspace{3mm}
$$
E^p_5 = I_5 - \sum_{j_1,\ldots,j_5=0}^{p}
C_{j_5\ldots j_1}\Biggl(\sum\limits_{(j_4,j_5)}\Biggl(
\sum\limits_{(j_1,j_2,j_3)}
C_{j_5\ldots j_1}\Biggr)\Biggr)\ \ \ (i_1=i_2=i_3\ne i_4=i_5),
$$

\vspace{3mm}
$$
E^p_5 = I_5 - \sum_{j_1,\ldots,j_5=0}^{p}
C_{j_5\ldots j_1}\Biggl(\sum\limits_{(j_1,j_3,j_4,j_5)}
C_{j_5\ldots j_1}\Biggr)\ \ \ (i_1=i_3=i_4=i_5\ne i_2).
$$

\vspace{5mm}

Consider a generalization of Theorem 10 to the case of an arbitrary 
complete ortho\-nor\-mal system of functions in the space $L_2([t, T])$
and $\psi_1(\tau),$ $\ldots,\psi_k(\tau)\in L_2([t, T])$

\vspace{2mm}              
               
{\bf Theorem 11} \cite{20a} (Sect.~1.12), \cite{26} (Sect.~6).
{\it Suppose that $\{\phi_j(x)\}_{j=0}^{\infty}$ 
is an arbitrary complete orthonormal system  
of functions in the space $L_2([t,T])$ and
$\psi_1(\tau),\ldots,\psi_k(\tau)\in L_2([t, T]),$  $i_1,\ldots, i_k=1,\ldots,m$.
Then

$$
{\sf M}\left\{\left(J[\psi^{(k)}]_{T,t}-
J[\psi^{(k)}]_{T,t}^{p}\right)^2\right\}=
\int\limits_{[t,T]^k}
K^2(t_1,\ldots,t_k)dt_1\ldots dt_k-
$$

\begin{equation}
\label{tttr11}
- \sum_{j_1,\ldots, j_k=0}^{p}
C_{j_k\ldots j_1}
{\sf M}\left\{J[\psi^{(k)}]_{T,t}
\sum\limits_{(j_1,\ldots,j_k)}
\int\limits_t^T \phi_{j_k}(t_k)
\ldots
\int\limits_t^{t_{2}}\phi_{j_{1}}(t_{1})
d{\bf f}_{t_1}^{(i_1)}\ldots
d{\bf f}_{t_k}^{(i_k)}\right\},
\end{equation}

\vspace{5mm}
\noindent
where $i_1,\ldots,i_k = 1,\ldots,m;$
the expression 

\vspace{-1mm}
$$
\sum\limits_{(j_1,\ldots,j_k)}
$$ 

\vspace{3mm}
\noindent
means the sum with respect to all
possible permutations 
$(j_1,\ldots,j_k)$. At the same time if 
$j_r$ swapped with $j_q$ in the permutation $(j_1,\ldots,j_k),$
then $i_r$ swapped with $i_q$ in the permutation
$(i_1,\ldots,i_k),$

$$
J[\psi^{(k)}]_{T,t}^p=
\sum\limits_{j_1,\ldots,j_k=0}^{p}
C_{j_k\ldots j_1}\Biggl(
\prod_{l=1}^k\zeta_{j_l}^{(i_l)}+\sum\limits_{r=1}^{[k/2]}
(-1)^r \times
\Biggr.
$$

\vspace{3mm}
$$
\times
\sum_{\stackrel{(\{\{g_1, g_2\}, \ldots, 
\{g_{2r-1}, g_{2r}\}\}, \{q_1, \ldots, q_{k-2r}\})}
{{}_{\{g_1, g_2, \ldots, 
g_{2r-1}, g_{2r}, q_1, \ldots, q_{k-2r}\}=\{1, 2, \ldots, k\}}}}
\prod\limits_{s=1}^r
{\bf 1}_{\{i_{g_{{}_{2s-1}}}=~i_{g_{{}_{2s}}}\ne 0\}}
\Biggl.{\bf 1}_{\{j_{g_{{}_{2s-1}}}=~j_{g_{{}_{2s}}}\}}
\prod_{l=1}^{k-2r}\zeta_{j_{q_l}}^{(i_{q_l})}\Biggr);
$$

\vspace{5mm}
\noindent
another notations are the same as in Theorems {\rm 2--4.}
}

\vspace{5mm}

\section{Approximation of Specific Iterated Ito
and Stratonovich Stochastic Integrals}

\vspace{5mm}

In this section we provide considerable practical material (based on 
Theorems 2--9 and Legendre
polynomials) concerning expansions and approximations of 
iterated Ito and 
Strato\-no\-vich stochastic integrals.
The question about what kind of functions (polynomial
or trigonometric) is more convenient for the mean-square approximation
of iterated stochastic integrals is also considered.

Let us consider the following iterated Ito and Stratonovich stochastic 
integrals

\vspace{-1mm}
\begin{equation}
\label{k1000}
I_{(l_1\ldots l_k)T,t}^{(i_1\ldots i_k)}
=\int\limits_t^T(t-t_k)^{l_k} \ldots \int\limits_t^{t_{2}}
(t-t_1)^{l_1} d{\bf f}_{t_1}^{(i_1)}\ldots
d{\bf f}_{t_k}^{(i_k)},
\end{equation}

\begin{equation}
\label{k1001}
I_{(l_1\ldots l_k)T,t}^{*(i_1\ldots i_k)}
={\int\limits_t^{*}}^T (t-t_k)^{l_k} \ldots {\int\limits_t^{*}}^{t_2}
(t-t_1)^{l_1} d{\bf f}_{t_1}^{(i_1)}\ldots
d{\bf f}_{t_k}^{(i_k)},
\end{equation}

\vspace{4mm}
\noindent
where $i_1,\ldots, i_k=1,\dots,m,$\ \ $l_1,\ldots,l_k=0, 1,\ldots$

The complete orthonormal system of Legendre polynomials in the 
space $L_2([t,T])$ looks as follows

\begin{equation}
\label{4009}
\phi_j(x)=\sqrt{\frac{2j+1}{T-t}}P_j\left(\left(
x-\frac{T+t}{2}\right)\frac{2}{T-t}\right),\ \ \ j=0, 1, 2,\ldots,
\end{equation}

\vspace{3mm}
\noindent
where $P_j(x)$ is the Legendre polynomial. 

Using the
system of 
functions (\ref{4009}) as well as Theorems 2--9,
we obtain the following expansions of iterated 
Ito and Stratonovich stochastic integrals (\ref{k1000}),
(\ref{k1001}) \cite{20a}

\vspace{2mm}
$$
I_{(0)T,t}^{(i_1)}=\sqrt{T-t}\zeta_0^{(i_1)},
$$

\vspace{2mm}
\begin{equation}
\label{4002}
I_{(1)T,t}^{(i_1)}=-\frac{(T-t)^{3/2}}{2}\left(\zeta_0^{(i_1)}+
\frac{1}{\sqrt{3}}\zeta_1^{(i_1)}\right),
\end{equation}

\vspace{3mm}
\begin{equation}
\label{4003}
I_{(2)T,t}^{(i_1)}=\frac{(T-t)^{5/2}}{3}\left(\zeta_0^{(i_1)}+
\frac{\sqrt{3}}{2}\zeta_1^{(i_1)}+
\frac{1}{2\sqrt{5}}\zeta_2^{(i_1)}\right),
\end{equation}

\vspace{3mm}
\begin{equation}
I_{(00)T,t}^{*(i_1 i_2)}=
\frac{T-t}{2}\left(\zeta_0^{(i_1)}\zeta_0^{(i_2)}+\sum_{i=1}^{\infty}
\frac{1}{\sqrt{4i^2-1}}\left(
\zeta_{i-1}^{(i_1)}\zeta_{i}^{(i_2)}-
\zeta_i^{(i_1)}\zeta_{i-1}^{(i_2)}\right)\right),
\label{4004}
\end{equation}

\vspace{6mm}

$$
I_{(00)T,t}^{(i_1 i_2)}=
\frac{T-t}{2}\left(\zeta_0^{(i_1)}\zeta_0^{(i_2)}+\sum_{i=1}^{\infty}
\frac{1}{\sqrt{4i^2-1}}\left(
\zeta_{i-1}^{(i_1)}\zeta_{i}^{(i_2)}-
\zeta_i^{(i_1)}\zeta_{i-1}^{(i_2)}\right) - {\bf 1}_{\{i_1=i_2\}}\right),
$$

\vspace{7mm}

$$
I_{(01)T,t}^{*(i_1 i_2)}=-\frac{T-t}{2}I_{(00)T,t}^{*(i_1 i_2)}
-\frac{(T-t)^2}{4}\Biggl(
\frac{1}{\sqrt{3}}\zeta_0^{(i_1)}\zeta_1^{(i_2)}+\Biggr.
$$

\vspace{2mm}
$$
+\Biggl.\sum_{i=0}^{\infty}\Biggl(
\frac{(i+2)\zeta_i^{(i_1)}\zeta_{i+2}^{(i_2)}
-(i+1)\zeta_{i+2}^{(i_1)}\zeta_{i}^{(i_2)}}
{\sqrt{(2i+1)(2i+5)}(2i+3)}-
\frac{\zeta_i^{(i_1)}\zeta_{i}^{(i_2)}}{(2i-1)(2i+3)}\Biggr)\Biggr),
$$

\vspace{7mm}

$$
I_{(10)T,t}^{*(i_1 i_2)}=-\frac{T-t}{2}I_{(00)T,t}^{*(i_1 i_2)}
-\frac{(T-t)^2}{4}\Biggl(
\frac{1}{\sqrt{3}}\zeta_0^{(i_2)}\zeta_1^{(i_1)}+\Biggr.
$$

\vspace{2mm}
\begin{equation}
\label{4006}
+\Biggl.\sum_{i=0}^{\infty}\Biggl(
\frac{(i+1)\zeta_{i+2}^{(i_2)}\zeta_{i}^{(i_1)}
-(i+2)\zeta_{i}^{(i_2)}\zeta_{i+2}^{(i_1)}}
{\sqrt{(2i+1)(2i+5)}(2i+3)}+
\frac{\zeta_i^{(i_1)}\zeta_{i}^{(i_2)}}{(2i-1)(2i+3)}\Biggr)\Biggr),
\end{equation}

\vspace{6mm}
\noindent 
or
$$
I_{(01)T,t}^{*(i_1i_2)}
=\hbox{\vtop{\offinterlineskip\halign{
\hfil#\hfil\cr
{\rm l.i.m.}\cr
$\stackrel{}{{}_{p\to \infty}}$\cr
}} }
\sum_{j_1,j_2=0}^{p}
C_{j_2j_1}^{01}
\zeta_{j_1}^{(i_1)}\zeta_{j_2}^{(i_2)},
$$

\vspace{2mm}
$$
I_{(10)T,t}^{*(i_1i_2)}
=\hbox{\vtop{\offinterlineskip\halign{
\hfil#\hfil\cr
{\rm l.i.m.}\cr
$\stackrel{}{{}_{p\to \infty}}$\cr
}} }
\sum_{j_1,j_2=0}^{p}
C_{j_2j_1}^{10}
\zeta_{j_1}^{(i_1)}\zeta_{j_2}^{(i_2)},
$$

\vspace{3mm}
\noindent
where

\vspace{-2mm}
$$
C_{j_2j_1}^{01}
=\frac{\sqrt{(2j_1+1)(2j_2+1)}}{8}(T-t)^{2}\bar
C_{j_2j_1}^{01},
$$

\vspace{1mm}
$$
C_{j_2j_1}^{10}
=\frac{\sqrt{(2j_1+1)(2j_2+1)}}{8}(T-t)^{2}\bar
C_{j_2j_1}^{10},
$$

\vspace{1mm}
$$
\bar C_{j_2j_1}^{01}=-\int\limits_{-1}^{1}(1+y)P_{j_2}(y)
\int\limits_{-1}^{y}
P_{j_1}(x)dx dy,
$$

\vspace{1mm}
$$
\bar C_{j_2j_1}^{10}=-\int\limits_{-1}^{1}P_{j_2}(y)
\int\limits_{-1}^{y}
(1+x)P_{j_1}(x)dx dy;
$$

\vspace{4mm}

$$
I_{(10)T,t}^{(i_1 i_2)}=
I_{(10)T,t}^{*(i_1 i_2)}+
\frac{1}{4}{\bf 1}_{\{i_1=i_2\}}(T-t)^2,\ \ \
I_{(01)T,t}^{(i_1 i_2)}=
I_{(01)T,t}^{*(i_1 i_2)}+
\frac{1}{4}{\bf 1}_{\{i_1=i_2\}}(T-t)^2\ \ \ {\rm w.\ p.\ 1},
$$

\vspace{5mm}

$$
I_{(01)T,t}^{(i_1 i_2)}=
-\frac{T-t}{2}
I_{(00)T,t}^{(i_1 i_2)}
-\frac{(T-t)^2}{4}\Biggl(
\frac{1}{\sqrt{3}}\zeta_0^{(i_1)}\zeta_1^{(i_2)}+\Biggr.
$$

\vspace{2mm}
$$
+\Biggl.\sum_{i=0}^{\infty}\Biggl(
\frac{(i+2)\zeta_i^{(i_1)}\zeta_{i+2}^{(i_2)}
-(i+1)\zeta_{i+2}^{(i_1)}\zeta_{i}^{(i_2)}}
{\sqrt{(2i+1)(2i+5)}(2i+3)}-
\frac{\zeta_i^{(i_1)}\zeta_{i}^{(i_2)}}{(2i-1)(2i+3)}\Biggr)\Biggr),
$$

\vspace{7mm}

$$
I_{(10)T,t}^{(i_1 i_2)}=
-\frac{T-t}{2}I_{(00)T,t}^{(i_1 i_2)}
-\frac{(T-t)^2}{4}\Biggl(
\frac{1}{\sqrt{3}}\zeta_0^{(i_2)}\zeta_1^{(i_1)}+\Biggr.
$$

\vspace{2mm}
$$
+\Biggl.\sum_{i=0}^{\infty}\Biggl(
\frac{(i+1)\zeta_{i+2}^{(i_2)}\zeta_{i}^{(i_1)}
-(i+2)\zeta_{i}^{(i_2)}\zeta_{i+2}^{(i_1)}}
{\sqrt{(2i+1)(2i+5)}(2i+3)}+
\frac{\zeta_i^{(i_1)}\zeta_{i}^{(i_2)}}{(2i-1)(2i+3)}\Biggr)\Biggr),
$$

\vspace{8mm}
\noindent
or
$$
I_{(01)T,t}^{(i_1 i_2)}=
\hbox{\vtop{\offinterlineskip\halign{
\hfil#\hfil\cr
{\rm l.i.m.}\cr
$\stackrel{}{{}_{p\to \infty}}$\cr
}} }
\sum_{j_1,j_2=0}^{p}
C_{j_2j_1}^{01}\Biggl(\zeta_{j_1}^{(i_1)}\zeta_{j_2}^{(i_2)}
-{\bf 1}_{\{i_1=i_2\}}
{\bf 1}_{\{j_1=j_2\}}\Biggr),
$$

\vspace{5mm}
$$
I_{(10)T,t}^{(i_1 i_2)}=
\hbox{\vtop{\offinterlineskip\halign{
\hfil#\hfil\cr
{\rm l.i.m.}\cr
$\stackrel{}{{}_{p\to \infty}}$\cr
}} }
\sum_{j_1,j_2=0}^{p}
C_{j_2j_1}^{10}\Biggl(\zeta_{j_1}^{(i_1)}\zeta_{j_2}^{(i_2)}
-{\bf 1}_{\{i_1=i_2\}}
{\bf 1}_{\{j_1=j_2\}}\Biggr),
$$

\vspace{7mm}

$$
I_{(000)T,t}^{*(i_1 i_2 i_3)}=
\hbox{\vtop{\offinterlineskip\halign{
\hfil#\hfil\cr
{\rm l.i.m.}\cr
$\stackrel{}{{}_{p\to \infty}}$\cr
}} }
\sum\limits_{j_1, j_2, j_3=0}^{p}
C_{j_3 j_2 j_1}\zeta_{j_1}^{(i_1)}\zeta_{j_2}^{(i_2)}\zeta_{j_3}^{(i_3)},
$$

\vspace{7mm}

$$
I_{(000)T,t}^{(i_1i_2i_3)}
=\hbox{\vtop{\offinterlineskip\halign{
\hfil#\hfil\cr
{\rm l.i.m.}\cr
$\stackrel{}{{}_{p\to \infty}}$\cr
}} }
\sum_{j_1,j_2,j_3=0}^{p}
C_{j_3j_2j_1}
\Biggl(
\zeta_{j_1}^{(i_1)}\zeta_{j_2}^{(i_2)}\zeta_{j_3}^{(i_3)}
-{\bf 1}_{\{i_1=i_2\ne 0\}}
{\bf 1}_{\{j_1=j_2\}}
\zeta_{j_3}^{(i_3)}-
\Biggr.
$$

\vspace{1mm}
\begin{equation}
\label{zzz1}
\Biggl.
-{\bf 1}_{\{i_2=i_3\ne 0\}}
{\bf 1}_{\{j_2=j_3\}}
\zeta_{j_1}^{(i_1)}-
{\bf 1}_{\{i_1=i_3\ne 0\}}
{\bf 1}_{\{j_1=j_3\}}
\zeta_{j_2}^{(i_2)}\Biggr),
\end{equation}

\vspace{7mm}

$$
I_{(000)T,t}^{(i_1 i_1 i_1)}
=\frac{1}{6}(T-t)^{3/2}\biggl(
\left(\zeta_0^{(i_1)}\right)^3-3
\zeta_0^{(i_1)}\biggr)\ \ \ \hbox{w.\ p.\ 1},
$$

\vspace{4mm}

$$
I_{(000)T,t}^{*(i_1 i_1 i_1)}
=\frac{1}{6}(T-t)^{3/2}
\left(\zeta_0^{(i_1)}\right)^3\ \ \ \hbox{w.\ p.\ 1},
$$

\vspace{5mm}
\noindent
where

\begin{equation}
\label{zzz2}
C_{j_3j_2j_1}
=\frac{\sqrt{(2j_1+1)(2j_2+1)(2j_3+1)}}{8}(T-t)^{3/2}\bar
C_{j_3j_2j_1},
\end{equation}

\begin{equation}
\label{zzz3}
\bar C_{j_3j_2j_1}=\int\limits_{-1}^{1}P_{j_3}(z)
\int\limits_{-1}^{z}P_{j_2}(y)
\int\limits_{-1}^{y}
P_{j_1}(x)dx dy dz;
\end{equation}

\vspace{6mm}

$$
I_{(000)T,t}^{(i_1 i_2 i_3)}=I_{(000)T,t}^{*(i_1 i_2 i_3)}
+{\bf 1}_{\{i_1=i_2\ne 0\}}\frac{1}{2}I_{(1)T,t}^{(i_3)}-
{\bf 1}_{\{i_2=i_3\ne 0\}}\frac{1}{2}\left((T-t)
I_{(0)T,t}^{(i_1)}+I_{(1)T,t}^{(i_1)}\right)\ \ \ \hbox{w.\ p.\ 1},
$$

\vspace{10mm}

$$
I_{(02)T,t}^{*(i_1 i_2)}=-\frac{(T-t)^2}{4}I_{(00)T,t}^{*(i_1 i_2)}
-(T-t)I_{(01)T,t}^{*(i_1 i_2)}+
\frac{(T-t)^3}{8}\Biggl[
\frac{2}{3\sqrt{5}}\zeta_2^{(i_2)}\zeta_0^{(i_1)}+\Biggr.
$$

\vspace{4mm}
$$
+\frac{1}{3}\zeta_0^{(i_1)}\zeta_0^{(i_2)}+
\sum_{i=0}^{\infty}\Biggl(
\frac{(i+2)(i+3)\zeta_{i+3}^{(i_2)}\zeta_{i}^{(i_1)}
-(i+1)(i+2)\zeta_{i}^{(i_2)}\zeta_{i+3}^{(i_1)}}
{\sqrt{(2i+1)(2i+7)}(2i+3)(2i+5)}+
\Biggr.
$$

\vspace{4mm}
\begin{equation}
\label{leto1}
\Biggl.\Biggl.+
\frac{(i^2+i-3)\zeta_{i+1}^{(i_2)}\zeta_{i}^{(i_1)}
-(i^2+3i-1)\zeta_{i}^{(i_2)}\zeta_{i+1}^{(i_1)}}
{\sqrt{(2i+1)(2i+3)}(2i-1)(2i+5)}\Biggr)\Biggr],
\end{equation}

\vspace{11mm}

$$
I_{(20)T,t}^{*(i_1 i_2)}=-\frac{(T-t)^2}{4}I_{(00)T,t}^{*(i_1 i_2)}
-(T-t)I_{(10)T,t}^{*(i_1 i_2)}+
\frac{(T-t)^3}{8}\Biggl[
\frac{2}{3\sqrt{5}}\zeta_0^{(i_2)}\zeta_2^{(i_1)}+\Biggr.
$$

\vspace{4mm}
$$
+\frac{1}{3}\zeta_0^{(i_1)}\zeta_0^{(i_2)}+
\sum_{i=0}^{\infty}\Biggl(
\frac{(i+1)(i+2)\zeta_{i+3}^{(i_2)}\zeta_{i}^{(i_1)}
-(i+2)(i+3)\zeta_{i}^{(i_2)}\zeta_{i+3}^{(i_1)}}
{\sqrt{(2i+1)(2i+7)}(2i+3)(2i+5)}+
\Biggr.
$$

\vspace{4mm}
\begin{equation}
\label{leto2}
\Biggl.\Biggl.+
\frac{(i^2+3i-1)\zeta_{i+1}^{(i_2)}\zeta_{i}^{(i_1)}
-(i^2+i-3)\zeta_{i}^{(i_2)}\zeta_{i+1}^{(i_1)}}
{\sqrt{(2i+1)(2i+3)}(2i-1)(2i+5)}\Biggr)\Biggr],
\end{equation}

\vspace{11mm}

$$
I_{(11)T,t}^{*(i_1 i_2)}=-\frac{(T-t)^2}{4}I_{(00)T,t}^{*(i_1 i_2)}-
\frac{(T-t)}{2}\left(I_{(10)T,t}^{*(i_1 i_2)}+
I_{(01)T,t}^{*(i_1 i_2)}\right)+
$$

\vspace{4mm}
$$
+\frac{(T-t)^3}{8}\Biggl[
\frac{1}{3}\zeta_1^{(i_1)}\zeta_1^{(i_2)}+
\sum_{i=0}^{\infty}\Biggl(
\frac{(i+1)(i+3)\left(\zeta_{i+3}^{(i_2)}\zeta_{i}^{(i_1)}
-\zeta_{i}^{(i_2)}\zeta_{i+3}^{(i_1)}\right)}
{\sqrt{(2i+1)(2i+7)}(2i+3)(2i+5)}+
\Biggr.
$$

\vspace{4mm}
\begin{equation}
\label{leto3}
\Biggl.\Biggl.+
\frac{(i+1)^2\left(\zeta_{i+1}^{(i_2)}\zeta_{i}^{(i_1)}
-\zeta_{i}^{(i_2)}\zeta_{i+1}^{(i_1)}\right)}
{\sqrt{(2i+1)(2i+3)}(2i-1)(2i+5)}\Biggr)\Biggr],
\end{equation}

\vspace{8mm}
or
$$
I_{(02)T,t}^{*(i_1 i_2)}=
\hbox{\vtop{\offinterlineskip\halign{
\hfil#\hfil\cr
{\rm l.i.m.}\cr
$\stackrel{}{{}_{p\to \infty}}$\cr
}} }
\sum_{j_1,j_2=0}^{p}
C_{j_2j_1}^{02}\zeta_{j_1}^{(i_1)}\zeta_{j_2}^{(i_2)},
$$

\vspace{3mm}

$$
I_{(20)T,t}^{*(i_1 i_2)}=
\hbox{\vtop{\offinterlineskip\halign{
\hfil#\hfil\cr
{\rm l.i.m.}\cr
$\stackrel{}{{}_{p\to \infty}}$\cr
}} }
\sum_{j_1,j_2=0}^{p}
C_{j_2j_1}^{20}\zeta_{j_1}^{(i_1)}\zeta_{j_2}^{(i_2)},
$$

\vspace{3mm}

$$
I_{(11)T,t}^{*(i_1 i_2)}=
\hbox{\vtop{\offinterlineskip\halign{
\hfil#\hfil\cr
{\rm l.i.m.}\cr
$\stackrel{}{{}_{p\to \infty}}$\cr
}} }
\sum_{j_1,j_2=0}^{p}
C_{j_2j_1}^{11}\zeta_{j_1}^{(i_1)}\zeta_{j_2}^{(i_2)},
$$

\vspace{5mm}
\noindent
where

\vspace{-2mm}
$$
C_{j_2j_1}^{02}=
\frac{\sqrt{(2j_1+1)(2j_2+1)}}{16}(T-t)^{3}\bar
C_{j_2j_1}^{02},
$$

\vspace{3mm}
$$
C_{j_2j_1}^{20}=
\frac{\sqrt{(2j_1+1)(2j_2+1)}}{16}(T-t)^{3}\bar
C_{j_2j_1}^{20},
$$

\vspace{3mm}

$$
C_{j_2j_1}^{11}=
\frac{\sqrt{(2j_1+1)(2j_2+1)}}{16}(T-t)^{3}\bar
C_{j_2j_1}^{11}, 
$$

\vspace{3mm}
$$
\bar C_{j_2j_1}^{02}=
\int\limits_{-1}^{1}P_{j_2}(y)(y+1)^2
\int\limits_{-1}^{y}
P_{j_1}(x)dx dy,
$$

\vspace{3mm}
$$
\bar C_{j_2j_1}^{20}=
\int\limits_{-1}^{1}P_{j_2}(y)
\int\limits_{-1}^{y}
P_{j_1}(x)(x+1)^2 dx dy,
$$

\vspace{3mm}
$$
\bar C_{j_2j_1}^{11}=
\int\limits_{-1}^{1}P_{j_2}(y)(y+1)
\int\limits_{-1}^{y}
P_{j_1}(x)(x+1)dx dy;
$$

\vspace{8mm}

$$
I_{(11)T,t}^{*(i_1 i_1)}=\frac{1}{2}\left(I_{(1)T,t}^{(i_1)}
\right)^2\ \ \ \hbox{w.\ p.\ 1,}
$$

\vspace{5mm}

$$
I_{(02)T,t}^{(i_1 i_2)}=
I_{(02)T,t}^{*(i_1 i_2)}-
\frac{1}{6}{\bf 1}_{\{i_1=i_2\}}(T-t)^3,\ \ \
I_{(20)T,t}^{(i_1 i_2)}=
I_{(20)T,t}^{*(i_1 i_2)}-
\frac{1}{6}{\bf 1}_{\{i_1=i_2\}}(T-t)^3\ \ \ \hbox{w.\ p.\ 1},
$$

\vspace{5mm}

$$
I_{(11)T,t}^{(i_1 i_2)}=
I_{(11)T,t}^{*(i_1 i_2)}-
\frac{1}{6}{\bf 1}_{\{i_1=i_2\}}(T-t)^3\ \ \ \hbox{w.\ p.\ 1},
$$

\vspace{11mm}

$$
I_{(02)T,t}^{(i_1 i_2)}
=-\frac{(T-t)^2}{4}I_{(00)T,t}^{(i_1 i_2)}
-(T-t) I_{01_{T,t}}^{(i_1 i_2)}+
\frac{(T-t)^3}{8}\Biggl[
\frac{2}{3\sqrt{5}}\zeta_2^{(i_2)}\zeta_0^{(i_1)}+\Biggr.
$$

\vspace{4mm}
$$
+\frac{1}{3}\zeta_0^{(i_1)}\zeta_0^{(i_2)}+
\sum_{i=0}^{\infty}\Biggl(
\frac{(i+2)(i+3)\zeta_{i+3}^{(i_2)}\zeta_{i}^{(i_1)}
-(i+1)(i+2)\zeta_{i}^{(i_2)}\zeta_{i+3}^{(i_1)}}
{\sqrt{(2i+1)(2i+7)}(2i+3)(2i+5)}+
\Biggr.
$$

\vspace{4mm}
$$
\Biggl.\Biggl.+
\frac{(i^2+i-3)\zeta_{i+1}^{(i_2)}\zeta_{i}^{(i_1)}
-(i^2+3i-1)\zeta_{i}^{(i_2)}\zeta_{i+1}^{(i_1)}}
{\sqrt{(2i+1)(2i+3)}(2i-1)(2i+5)}\Biggr)\Biggr] - 
\frac{1}{24}{\bf 1}_{\{i_1=i_2\}}{(T-t)^3},
$$

\vspace{11mm}

$$
I_{(20)T,t}^{(i_1 i_2)}=-\frac{(T-t)^2}{4}
I_{(00)T,t}^{(i_1 i_2)}
-(T-t) I_{(10)T,t}^{(i_1 i_2)}+
\frac{(T-t)^3}{8}\Biggl[
\frac{2}{3\sqrt{5}}\zeta_0^{(i_2)}\zeta_2^{(i_1)}+\Biggr.
$$

\vspace{4mm}
$$
+\frac{1}{3}\zeta_0^{(i_1)}\zeta_0^{(i_2)}+
\sum_{i=0}^{\infty}\Biggl(
\frac{(i+1)(i+2)\zeta_{i+3}^{(i_2)}\zeta_{i}^{(i_1)}
-(i+2)(i+3)\zeta_{i}^{(i_2)}\zeta_{i+3}^{(i_1)}}
{\sqrt{(2i+1)(2i+7)}(2i+3)(2i+5)}+
\Biggr.
$$

\vspace{4mm}
$$
\Biggl.\Biggl.+
\frac{(i^2+3i-1)\zeta_{i+1}^{(i_2)}\zeta_{i}^{(i_1)}
-(i^2+i-3)\zeta_{i}^{(i_2)}\zeta_{i+1}^{(i_1)}}
{\sqrt{(2i+1)(2i+3)}(2i-1)(2i+5)}\Biggr)\Biggr] - 
\frac{1}{24}{\bf 1}_{\{i_1=i_2\}}{(T-t)^3},
$$

\vspace{11mm}

$$
I_{(11)T,t}^{(i_1 i_2)}
=-\frac{(T-t)^2}{4}I_{(00)T,t}^{(i_1 i_2)}
-\frac{T-t}{2}\left(
I_{(10)T,t}^{(i_1 i_2)}+
I_{(01)T,t}^{(i_1 i_2)}\right)+
\frac{(T-t)^3}{8}\Biggl[
\frac{1}{3}\zeta_1^{(i_1)}\zeta_1^{(i_2)}+\Biggr.
$$

\vspace{4mm}
$$
+
\sum_{i=0}^{\infty}\Biggl(
\frac{(i+1)(i+3)\left(\zeta_{i+3}^{(i_2)}\zeta_{i}^{(i_1)}
-\zeta_{i}^{(i_2)}\zeta_{i+3}^{(i_1)}\right)}
{\sqrt{(2i+1)(2i+7)}(2i+3)(2i+5)}+
\Biggr.
$$

\vspace{4mm}
$$
\Biggl.\Biggl.
+\frac{(i+1)^2\left(\zeta_{i+1}^{(i_2)}\zeta_{i}^{(i_1)}
-\zeta_{i}^{(i_2)}\zeta_{i+1}^{(i_1)}\right)}
{\sqrt{(2i+1)(2i+3)}(2i-1)(2i+5)}\Biggr)\Biggr] - 
\frac{1}{24}{\bf 1}_{\{i_1=i_2\}}{(T-t)^3},
$$

\vspace{8mm}
or
$$
I_{(02)T,t}^{(i_1 i_2)}=
\hbox{\vtop{\offinterlineskip\halign{
\hfil#\hfil\cr
{\rm l.i.m.}\cr
$\stackrel{}{{}_{p\to \infty}}$\cr
}} }
\sum_{j_1,j_2=0}^p
C_{j_2j_1}^{02}\Biggl(\zeta_{j_1}^{(i_1)}\zeta_{j_2}^{(i_2)}
-{\bf 1}_{\{i_1=i_2\}}
{\bf 1}_{\{j_1=j_2\}}\Biggr),
$$

\vspace{4mm}
$$
I_{(20)T,t}^{(i_1 i_2)}=
\hbox{\vtop{\offinterlineskip\halign{
\hfil#\hfil\cr
{\rm l.i.m.}\cr
$\stackrel{}{{}_{p\to \infty}}$\cr
}} }
\sum_{j_1,j_2=0}^{p}
C_{j_2j_1}^{20}\Biggl(\zeta_{j_1}^{(i_1)}\zeta_{j_2}^{(i_2)}
-{\bf 1}_{\{i_1=i_2\}}
{\bf 1}_{\{j_1=j_2\}}\Biggr),
$$

\vspace{4mm}
$$
I_{(11)T,t}^{(i_1 i_2)}=
\hbox{\vtop{\offinterlineskip\halign{
\hfil#\hfil\cr
{\rm l.i.m.}\cr
$\stackrel{}{{}_{p\to \infty}}$\cr
}} }
\sum_{j_1,j_2=0}^{p}
C_{j_2j_1}^{11}\Biggl(\zeta_{j_1}^{(i_1)}\zeta_{j_2}^{(i_2)}
-{\bf 1}_{\{i_1=i_2\}}
{\bf 1}_{\{j_1=j_2\}}\Biggr),
$$

\vspace{9mm}

$$
I_{(3)T,t}^{(i_1)}=-\frac{(T-t)^{7/2}}{4}\left(\zeta_0^{(i_1)}+
\frac{3\sqrt{3}}{5}\zeta_1^{(i_1)}+
\frac{1}{\sqrt{5}}\zeta_2^{(i_1)}+
\frac{1}{5\sqrt{7}}\zeta_3^{(i_1)}\right),
$$

\vspace{7mm}

$$
I_{(0000)T,t}^{*(i_1 i_2 i_3 i_4)}=
\hbox{\vtop{\offinterlineskip\halign{
\hfil#\hfil\cr
{\rm l.i.m.}\cr
$\stackrel{}{{}_{p\to \infty}}$\cr
}} }
\sum\limits_{j_1, j_2, j_3, j_4=0}^{p}
C_{j_4 j_3 j_2 j_1}\zeta_{j_1}^{(i_1)}\zeta_{j_2}^{(i_2)}\zeta_{j_3}^{(i_3)}
\zeta_{j_4}^{(i_4)},
$$

\vspace{8mm}

$$
I_{(0000)T,t}^{(i_1 i_2 i_3 i_4)}
=
\hbox{\vtop{\offinterlineskip\halign{
\hfil#\hfil\cr
{\rm l.i.m.}\cr
$\stackrel{}{{}_{p\to \infty}}$\cr
}} }
\sum_{j_1,j_2,j_3,j_4=0}^{p}
C_{j_4 j_3 j_2 j_1}\Biggl(
\prod_{l=1}^4\zeta_{j_l}^{(i_l)}
-\Biggr.
$$
$$
-
{\bf 1}_{\{i_1=i_2\ne 0\}}
{\bf 1}_{\{j_1=j_2\}}
\zeta_{j_3}^{(i_3)}
\zeta_{j_4}^{(i_4)}
-
{\bf 1}_{\{i_1=i_3\ne 0\}}
{\bf 1}_{\{j_1=j_3\}}
\zeta_{j_2}^{(i_2)}
\zeta_{j_4}^{(i_4)}-
$$
$$
-
{\bf 1}_{\{i_1=i_4\ne 0\}}
{\bf 1}_{\{j_1=j_4\}}
\zeta_{j_2}^{(i_2)}
\zeta_{j_3}^{(i_3)}
-
{\bf 1}_{\{i_2=i_3\ne 0\}}
{\bf 1}_{\{j_2=j_3\}}
\zeta_{j_1}^{(i_1)}
\zeta_{j_4}^{(i_4)}-
$$
$$
-
{\bf 1}_{\{i_2=i_4\ne 0\}}
{\bf 1}_{\{j_2=j_4\}}
\zeta_{j_1}^{(i_1)}
\zeta_{j_3}^{(i_3)}
-
{\bf 1}_{\{i_3=i_4\ne 0\}}
{\bf 1}_{\{j_3=j_4\}}
\zeta_{j_1}^{(i_1)}
\zeta_{j_2}^{(i_2)}+
$$
$$
+
{\bf 1}_{\{i_1=i_2\ne 0\}}
{\bf 1}_{\{j_1=j_2\}}
{\bf 1}_{\{i_3=i_4\ne 0\}}
{\bf 1}_{\{j_3=j_4\}}+
$$
$$
+
{\bf 1}_{\{i_1=i_3\ne 0\}}
{\bf 1}_{\{j_1=j_3\}}
{\bf 1}_{\{i_2=i_4\ne 0\}}
{\bf 1}_{\{j_2=j_4\}}+
$$
\begin{equation}
\label{zzz10}
+\Biggl.
{\bf 1}_{\{i_1=i_4\ne 0\}}
{\bf 1}_{\{j_1=j_4\}}
{\bf 1}_{\{i_2=i_3\ne 0\}}
{\bf 1}_{\{j_2=j_3\}}\Biggr),
\end{equation}

\vspace{7mm}

$$
I_{(0000)T,t}^{(i_1i_1i_1i_1)}=
\frac{1}{24}(T-t)^2
\left(\left(\zeta_0^{(i_1)}\right)^4-
6\left(\zeta_0^{(i_1)}\right)^2+3\right)\ \ \ \hbox{w.\ p.\ 1},
$$

\vspace{4mm}

$$
I_{(0000)T,t}^{*(i_1i_1i_1i_1)}=
\frac{1}{24}(T-t)^2
\left(\zeta_0^{(i_1)}\right)^4\ \ \ \hbox{w.\ p.\ 1},
$$

\vspace{4mm}
\noindent
where

\begin{equation}
\label{zzz11}
C_{j_4j_3j_2j_1}=
\frac{\sqrt{(2j_1+1)(2j_2+1)(2j_3+1)(2j_4+1)}}{16}(T-t)^{2}\bar
C_{j_4j_3j_2j_1},
\end{equation}

\vspace{2mm}

\begin{equation}
\label{zzz12}
\bar C_{j_4j_3j_2j_1}=\int\limits_{-1}^{1}P_{j_4}(u)
\int\limits_{-1}^{u}P_{j_3}(z)
\int\limits_{-1}^{z}P_{j_2}(y)
\int\limits_{-1}^{y}
P_{j_1}(x)dx dy dz du;
\end{equation}

\vspace{6mm}

$$
I_{(001)T,t}^{*(i_1i_2i_3)}
=\hbox{\vtop{\offinterlineskip\halign{
\hfil#\hfil\cr
{\rm l.i.m.}\cr
$\stackrel{}{{}_{p\to \infty}}$\cr
}} }
\sum_{j_1,j_2,j_3=0}^{p}
C_{j_3 j_2 j_1}^{001}
\zeta_{j_1}^{(i_1)}\zeta_{j_2}^{(i_2)}\zeta_{j_3}^{(i_3)},
$$

\vspace{4mm}

$$
I_{(010)T,t}^{*(i_1i_2i_3)}
=\hbox{\vtop{\offinterlineskip\halign{
\hfil#\hfil\cr
{\rm l.i.m.}\cr
$\stackrel{}{{}_{p\to \infty}}$\cr
}} }
\sum_{j_1,j_2,j_3=0}^{p}
C_{j_3 j_2 j_1}^{010}
\zeta_{j_1}^{(i_1)}\zeta_{j_2}^{(i_2)}\zeta_{j_3}^{(i_3)},
$$

\vspace{4mm}

$$
I_{(100)T,t}^{*(i_1i_2i_3)}
=\hbox{\vtop{\offinterlineskip\halign{
\hfil#\hfil\cr
{\rm l.i.m.}\cr
$\stackrel{}{{}_{p\to \infty}}$\cr
}} }
\sum_{j_1,j_2,j_3=0}^{p}
C_{j_3 j_2 j_1}^{100}
\zeta_{j_1}^{(i_1)}\zeta_{j_2}^{(i_2)}\zeta_{j_3}^{(i_3)},
$$

\vspace{9mm}

$$
I_{(001)T,t}^{(i_1i_2i_3)}
=\hbox{\vtop{\offinterlineskip\halign{
\hfil#\hfil\cr
{\rm l.i.m.}\cr
$\stackrel{}{{}_{p\to \infty}}$\cr
}} }
\sum_{j_1,j_2,j_3=0}^{p}
C_{j_3j_2j_1}^{001}\Biggl(
\zeta_{j_1}^{(i_1)}\zeta_{j_2}^{(i_2)}\zeta_{j_3}^{(i_3)}
-{\bf 1}_{\{i_1=i_2\}}
{\bf 1}_{\{j_1=j_2\}}
\zeta_{j_3}^{(i_3)}-
\Biggr.
$$

\vspace{1mm}
\begin{equation}
\label{sss1}
\Biggl.
-{\bf 1}_{\{i_2=i_3\}}
{\bf 1}_{\{j_2=j_3\}}
\zeta_{j_1}^{(i_1)}-
{\bf 1}_{\{i_1=i_3\}}
{\bf 1}_{\{j_1=j_3\}}
\zeta_{j_2}^{(i_2)}\Biggr),
\end{equation}

\vspace{5mm}

$$
I_{(010)T,t}^{(i_1i_2i_3)}
=\hbox{\vtop{\offinterlineskip\halign{
\hfil#\hfil\cr
{\rm l.i.m.}\cr
$\stackrel{}{{}_{p\to \infty}}$\cr
}} }
\sum_{j_1,j_2,j_3=0}^{p}
C_{j_3j_2j_1}^{010}\Biggl(
\zeta_{j_1}^{(i_1)}\zeta_{j_2}^{(i_2)}\zeta_{j_3}^{(i_3)}
-{\bf 1}_{\{i_1=i_2\}}
{\bf 1}_{\{j_1=j_2\}}
\zeta_{j_3}^{(i_3)}-
\Biggr.
$$

\vspace{1mm}
\begin{equation}
\label{sss2}
\Biggl.
-{\bf 1}_{\{i_2=i_3\}}
{\bf 1}_{\{j_2=j_3\}}
\zeta_{j_1}^{(i_1)}-
{\bf 1}_{\{i_1=i_3\}}
{\bf 1}_{\{j_1=j_3\}}
\zeta_{j_2}^{(i_2)}\Biggr),
\end{equation}

\vspace{5mm}

$$
I_{(100)T,t}^{(i_1i_2i_3)}
=\hbox{\vtop{\offinterlineskip\halign{
\hfil#\hfil\cr
{\rm l.i.m.}\cr
$\stackrel{}{{}_{p\to \infty}}$\cr
}} }
\sum_{j_1,j_2,j_3=0}^{p}
C_{j_3j_2j_1}^{100}\Biggl(
\zeta_{j_1}^{(i_1)}\zeta_{j_2}^{(i_2)}\zeta_{j_3}^{(i_3)}
-{\bf 1}_{\{i_1=i_2\}}
{\bf 1}_{\{j_1=j_2\}}
\zeta_{j_3}^{(i_3)}-
\Biggr.
$$

\vspace{1mm}
\begin{equation}
\label{sss3}
\Biggl.
-{\bf 1}_{\{i_2=i_3\}}
{\bf 1}_{\{j_2=j_3\}}
\zeta_{j_1}^{(i_1)}-
{\bf 1}_{\{i_1=i_3\}}
{\bf 1}_{\{j_1=j_3\}}
\zeta_{j_2}^{(i_2)}\Biggr),
\end{equation}

\vspace{5mm}
\noindent
where

\vspace{-3mm}
$$
C_{j_3j_2j_1}^{001}
=\frac{\sqrt{(2j_1+1)(2j_2+1)(2j_3+1)}}{16}(T-t)^{5/2}\bar
C_{j_3j_2j_1}^{001},
$$

\vspace{4mm}

$$
C_{j_3j_2j_1}^{010}
=\frac{\sqrt{(2j_1+1)(2j_2+1)(2j_3+1)}}{16}(T-t)^{5/2}\bar
C_{j_3j_2j_1}^{010},
$$

\vspace{4mm}

$$
C_{j_3j_2j_1}^{100}
=\frac{\sqrt{(2j_1+1)(2j_2+1)(2j_3+1)}}{16}(T-t)^{5/2}\bar
C_{j_3j_2j_1}^{100},
$$

\vspace{6mm}

$$
\bar C_{j_3j_2j_1}^{100}=-
\int\limits_{-1}^{1}P_{j_3}(z)
\int\limits_{-1}^{z}P_{j_2}(y)
\int\limits_{-1}^{y}
P_{j_1}(x)(x+1)dx dy dz,
$$

\vspace{4mm}

$$
\bar C_{j_3j_2j_1}^{010}=-
\int\limits_{-1}^{1}P_{j_3}(z)
\int\limits_{-1}^{z}P_{j_2}(y)(y+1)
\int\limits_{-1}^{y}
P_{j_1}(x)dx dy dz,
$$

\vspace{4mm}

$$
\bar C_{j_3j_2j_1}^{001}=-
\int\limits_{-1}^{1}P_{j_3}(z)(z+1)
\int\limits_{-1}^{z}P_{j_2}(y)
\int\limits_{-1}^{y}
P_{j_1}(x)dx dy dz;
$$

\vspace{9mm}

$$
I_{(lll)T,t}^{(i_1i_1i_1)}=
\frac{1}{6}\left(\left(I_{(l)T,t}^{(i_1)}\right)^3-
3I_{(l)T,t}^{(i_1)}\Delta_{l(T,t)}\right)\ \ \ \hbox{w.\ p.\ 1},
$$

\vspace{4mm}
$$
I_{(lll)T,t}^{*(iii)}=
\frac{1}{6}\left(I_{(l)T,t}^{(i_1)}\right)^3\ \ \ \hbox{w.\ p.\ 1},
$$

\vspace{4mm}

$$
I_{(llll)T,t}^{(i_1i_1i_1i_1)}=
\frac{1}{24}\left(\left(I_{(l)T,t}^{(i_1)}\right)^4-
6\left(I_{(l)T,t}^{(i_1)}\right)^2\Delta_{(l)T,t}+3
\left(\Delta_{(l)T,t}\right)^2\right)\ \ \ \hbox{w.\ p.\ 1},
$$

\vspace{4mm}

$$
I_{(llll)T,t}^{*(i_1i_1i_1i_1)}=
\frac{1}{24}\left(I_{(l)T,t}^{(i_1)}\right)^4\ \ \ \hbox{w.\ p.\ 1},
$$

\vspace{4mm}
\noindent
where

\vspace{-4mm}
$$
I_{(l)T,t}^{(i_1)}=\sum_{j=0}^l C_j^l \zeta_j^{(i_1)}\ \ \ \hbox{w.\ p.\ 1},
$$       

\vspace{4mm}

$$
\Delta_{l(T,t)}=\int\limits_t^T(t-s)^{2l}ds,\ \ \
C_j^l=\int\limits_t^T(t-s)^l\phi_j(s)ds;
$$

\vspace{5mm}

$$
I_{(00000)T,t}^{*(i_1 i_2 i_3 i_4 i_5)}=
\hbox{\vtop{\offinterlineskip\halign{
\hfil#\hfil\cr
{\rm l.i.m.}\cr
$\stackrel{}{{}_{p\to \infty}}$\cr
}} }
\sum\limits_{j_1, j_2, j_3, j_4, j_5=0}^{p}
C_{j_5j_4 j_3 j_2 j_1}
\zeta_{j_1}^{(i_1)}\zeta_{j_2}^{(i_2)}\zeta_{j_3}^{(i_3)}
\zeta_{j_4}^{(i_4)}\zeta_{j_5}^{(i_5)},
$$

\vspace{7mm}

$$
I_{(00000)T,t}^{(i_1 i_2 i_3 i_4 i_5)}
=
\hbox{\vtop{\offinterlineskip\halign{
\hfil#\hfil\cr
{\rm l.i.m.}\cr
$\stackrel{}{{}_{p\to \infty}}$\cr
}} }
\sum_{j_1,j_2,j_3,j_4,j_5=0}^p
C_{j_5 j_4 j_3 j_2 j_1}\Biggl(
\prod_{l=1}^5\zeta_{j_l}^{(i_l)}
-\Biggr.
$$
$$
-
{\bf 1}_{\{i_1=i_2\}}
{\bf 1}_{\{j_1=j_2\}}
\zeta_{j_3}^{(i_3)}
\zeta_{j_4}^{(i_4)}
\zeta_{j_5}^{(i_5)}-
{\bf 1}_{\{i_1=i_3\}}
{\bf 1}_{\{j_1=j_3\}}
\zeta_{j_2}^{(i_2)}
\zeta_{j_4}^{(i_4)}
\zeta_{j_5}^{(i_5)}-
$$
$$
-
{\bf 1}_{\{i_1=i_4\}}
{\bf 1}_{\{j_1=j_4\}}
\zeta_{j_2}^{(i_2)}
\zeta_{j_3}^{(i_3)}
\zeta_{j_5}^{(i_5)}-
{\bf 1}_{\{i_1=i_5\}}
{\bf 1}_{\{j_1=j_5\}}
\zeta_{j_2}^{(i_2)}
\zeta_{j_3}^{(i_3)}
\zeta_{j_4}^{(i_4)}-
$$
$$
-
{\bf 1}_{\{i_2=i_3\}}
{\bf 1}_{\{j_2=j_3\}}
\zeta_{j_1}^{(i_1)}
\zeta_{j_4}^{(i_4)}
\zeta_{j_5}^{(i_5)}-
{\bf 1}_{\{i_2=i_4\}}
{\bf 1}_{\{j_2=j_4\}}
\zeta_{j_1}^{(i_1)}
\zeta_{j_3}^{(i_3)}
\zeta_{j_5}^{(i_5)}-
$$
$$
-
{\bf 1}_{\{i_2=i_5\}}
{\bf 1}_{\{j_2=j_5\}}
\zeta_{j_1}^{(i_1)}
\zeta_{j_3}^{(i_3)}
\zeta_{j_4}^{(i_4)}
-{\bf 1}_{\{i_3=i_4\}}
{\bf 1}_{\{j_3=j_4\}}
\zeta_{j_1}^{(i_1)}
\zeta_{j_2}^{(i_2)}
\zeta_{j_5}^{(i_5)}-
$$
$$
-
{\bf 1}_{\{i_3=i_5\}}
{\bf 1}_{\{j_3=j_5\}}
\zeta_{j_1}^{(i_1)}
\zeta_{j_2}^{(i_2)}
\zeta_{j_4}^{(i_4)}
-{\bf 1}_{\{i_4=i_5\}}
{\bf 1}_{\{j_4=j_5\}}
\zeta_{j_1}^{(i_1)}
\zeta_{j_2}^{(i_2)}
\zeta_{j_3}^{(i_3)}+
$$
$$
+
{\bf 1}_{\{i_1=i_2\}}
{\bf 1}_{\{j_1=j_2\}}
{\bf 1}_{\{i_3=i_4\}}
{\bf 1}_{\{j_3=j_4\}}\zeta_{j_5}^{(i_5)}+
{\bf 1}_{\{i_1=i_2\}}
{\bf 1}_{\{j_1=j_2\}}
{\bf 1}_{\{i_3=i_5\}}
{\bf 1}_{\{j_3=j_5\}}\zeta_{j_4}^{(i_4)}+
$$
$$
+
{\bf 1}_{\{i_1=i_2\}}
{\bf 1}_{\{j_1=j_2\}}
{\bf 1}_{\{i_4=i_5\}}
{\bf 1}_{\{j_4=j_5\}}\zeta_{j_3}^{(i_3)}+
{\bf 1}_{\{i_1=i_3\}}
{\bf 1}_{\{j_1=j_3\}}
{\bf 1}_{\{i_2=i_4\}}
{\bf 1}_{\{j_2=j_4\}}\zeta_{j_5}^{(i_5)}+
$$
$$
+
{\bf 1}_{\{i_1=i_3\}}
{\bf 1}_{\{j_1=j_3\}}
{\bf 1}_{\{i_2=i_5\}}
{\bf 1}_{\{j_2=j_5\}}\zeta_{j_4}^{(i_4)}+
{\bf 1}_{\{i_1=i_3\}}
{\bf 1}_{\{j_1=j_3\}}
{\bf 1}_{\{i_4=i_5\}}
{\bf 1}_{\{j_4=j_5\}}\zeta_{j_2}^{(i_2)}+
$$
$$
+
{\bf 1}_{\{i_1=i_4\}}
{\bf 1}_{\{j_1=j_4\}}
{\bf 1}_{\{i_2=i_3\}}
{\bf 1}_{\{j_2=j_3\}}\zeta_{j_5}^{(i_5)}+
{\bf 1}_{\{i_1=i_4\}}
{\bf 1}_{\{j_1=j_4\}}
{\bf 1}_{\{i_2=i_5\}}
{\bf 1}_{\{j_2=j_5\}}\zeta_{j_3}^{(i_3)}+
$$
$$
+
{\bf 1}_{\{i_1=i_4\}}
{\bf 1}_{\{j_1=j_4\}}
{\bf 1}_{\{i_3=i_5\}}
{\bf 1}_{\{j_3=j_5\}}\zeta_{j_2}^{(i_2)}+
{\bf 1}_{\{i_1=i_5\}}
{\bf 1}_{\{j_1=j_5\}}
{\bf 1}_{\{i_2=i_3\}}
{\bf 1}_{\{j_2=j_3\}}\zeta_{j_4}^{(i_4)}+
$$
$$
+
{\bf 1}_{\{i_1=i_5\}}
{\bf 1}_{\{j_1=j_5\}}
{\bf 1}_{\{i_2=i_4\}}
{\bf 1}_{\{j_2=j_4\}}\zeta_{j_3}^{(i_3)}+
{\bf 1}_{\{i_1=i_5\}}
{\bf 1}_{\{j_1=j_5\}}
{\bf 1}_{\{i_3=i_4\}}
{\bf 1}_{\{j_3=j_4\}}\zeta_{j_2}^{(i_2)}+
$$
$$
+
{\bf 1}_{\{i_2=i_3\}}
{\bf 1}_{\{j_2=j_3\}}
{\bf 1}_{\{i_4=i_5\}}
{\bf 1}_{\{j_4=j_5\}}\zeta_{j_1}^{(i_1)}+
{\bf 1}_{\{i_2=i_4\}}
{\bf 1}_{\{j_2=j_4\}}
{\bf 1}_{\{i_3=i_5\}}
{\bf 1}_{\{j_3=j_5\}}\zeta_{j_1}^{(i_1)}+
$$
\begin{equation}
\label{sss4}
+\Biggl.
{\bf 1}_{\{i_2=i_5\}}
{\bf 1}_{\{j_2=j_5\}}
{\bf 1}_{\{i_3=i_4\}}
{\bf 1}_{\{j_3=j_4\}}\zeta_{j_1}^{(i_1)}\Biggr),
\end{equation}

\vspace{7mm}

$$         
I_{(00000)T,t}^{(i_1i_1i_1i_1i_1)}=
\frac{1}{120}(T-t)^{5/2}
\left(\left(\zeta_0^{(i_1)}\right)^5-
10\left(\zeta_0^{(i_1)}\right)^3+15\zeta_0^{(i_1)}\right)\ \ \ 
\hbox{w.\ p.\ 1},
$$

\vspace{2mm}

$$
I_{(00000)T,t}^{*(i_1i_1i_1i_1i_1)}=
\frac{1}{120}(T-t)^{5/2}\left(\zeta_0^{(i_1)}\right)^5\ \ \ \hbox{w.\ p.\ 1},
$$

\vspace{5mm}
\noindent
where

\vspace{2mm}
$$
C_{j_5j_4 j_3 j_2 j_1}=
\frac{\sqrt{(2j_1+1)(2j_2+1)(2j_3+1)(2j_4+1)(2j_5+1)}}{32}(T-t)^{5/2}\bar
C_{j_5j_4 j_3 j_2 j_1},
$$

\vspace{3mm}

$$
\bar C_{j_5j_4 j_3 j_2 j_1}=
\int\limits_{-1}^{1}P_{j_5}(v)
\int\limits_{-1}^{v}P_{j_4}(u)
\int\limits_{-1}^{u}P_{j_3}(z)
\int\limits_{-1}^{z}P_{j_2}(y)
\int\limits_{-1}^{y}
P_{j_1}(x)dx dy dz du dv;
$$

\vspace{4mm}

$$
I_{(0001)T,t}^{*(i_1i_2i_3)}
=\hbox{\vtop{\offinterlineskip\halign{
\hfil#\hfil\cr
{\rm l.i.m.}\cr
$\stackrel{}{{}_{p\to \infty}}$\cr
}} }
\sum_{j_1,j_2,j_3,j_4=0}^{p}
C_{j_4j_3 j_2 j_1}^{0001}
\zeta_{j_1}^{(i_1)}\zeta_{j_2}^{(i_2)}\zeta_{j_3}^{(i_3)}\zeta_{j_4}^{(i_4)},
$$

\vspace{2mm}
$$
I_{(0010)T,t}^{*(i_1i_2i_3)}
=\hbox{\vtop{\offinterlineskip\halign{
\hfil#\hfil\cr
{\rm l.i.m.}\cr
$\stackrel{}{{}_{p\to \infty}}$\cr
}} }
\sum_{j_1,j_2,j_3,j_4=0}^{p}
C_{j_4j_3 j_2 j_1}^{0010}
\zeta_{j_1}^{(i_1)}\zeta_{j_2}^{(i_2)}\zeta_{j_3}^{(i_3)}\zeta_{j_4}^{(i_4)},
$$

\vspace{2mm}

$$
I_{(0100)T,t}^{*(i_1i_2i_3)}
=\hbox{\vtop{\offinterlineskip\halign{
\hfil#\hfil\cr
{\rm l.i.m.}\cr
$\stackrel{}{{}_{p\to \infty}}$\cr
}} }
\sum_{j_1,j_2,j_3,j_4=0}^{p}
C_{j_4j_3 j_2 j_1}^{0100}
\zeta_{j_1}^{(i_1)}\zeta_{j_2}^{(i_2)}\zeta_{j_3}^{(i_3)}\zeta_{j_4}^{(i_4)},
$$

\vspace{2mm}

$$
I_{(1000)T,t}^{*(i_1i_2i_3)}
=\hbox{\vtop{\offinterlineskip\halign{
\hfil#\hfil\cr
{\rm l.i.m.}\cr
$\stackrel{}{{}_{p\to \infty}}$\cr
}} }
\sum_{j_1,j_2,j_3,j_4=0}^{p}
C_{j_4j_3 j_2 j_1}^{1000}
\zeta_{j_1}^{(i_1)}\zeta_{j_2}^{(i_2)}\zeta_{j_3}^{(i_3)}\zeta_{j_4}^{(i_4)},
$$

\vspace{6mm}

$$
I_{(0001)T,t}^{(i_1 i_2 i_3 i_4)}
=\hbox{\vtop{\offinterlineskip\halign{
\hfil#\hfil\cr
{\rm l.i.m.}\cr
$\stackrel{}{{}_{p\to \infty}}$\cr
}} }
\sum_{j_1,j_2,j_3,j_4=0}^{p}
C_{j_4 j_3 j_2 j_1}^{0001}\Biggl(
\zeta_{j_1}^{(i_1)}\zeta_{j_2}^{(i_2)}\zeta_{j_3}^{(i_3)}\zeta_{j_4}^{(i_4)}
-\Biggr.
$$
$$
-
{\bf 1}_{\{i_1=i_2\}}
{\bf 1}_{\{j_1=j_2\}}
\zeta_{j_3}^{(i_3)}
\zeta_{j_4}^{(i_4)}
-
{\bf 1}_{\{i_1=i_3\}}
{\bf 1}_{\{j_1=j_3\}}
\zeta_{j_2}^{(i_2)}
\zeta_{j_4}^{(i_4)}-
$$
$$
-
{\bf 1}_{\{i_1=i_4\}}
{\bf 1}_{\{j_1=j_4\}}
\zeta_{j_2}^{(i_2)}
\zeta_{j_3}^{(i_3)}
-
{\bf 1}_{\{i_2=i_3\}}
{\bf 1}_{\{j_2=j_3\}}
\zeta_{j_1}^{(i_1)}
\zeta_{j_4}^{(i_4)}-
$$
$$
-
{\bf 1}_{\{i_2=i_4\}}
{\bf 1}_{\{j_2=j_4\}}
\zeta_{j_1}^{(i_1)}
\zeta_{j_3}^{(i_3)}
-
{\bf 1}_{\{i_3=i_4\}}
{\bf 1}_{\{j_3=j_4\}}
\zeta_{j_1}^{(i_1)}
\zeta_{j_2}^{(i_2)}+
$$
$$
+
{\bf 1}_{\{i_1=i_2\}}
{\bf 1}_{\{j_1=j_2\}}
{\bf 1}_{\{i_3=i_4\}}
{\bf 1}_{\{j_3=j_4\}}+
{\bf 1}_{\{i_1=i_3\}}
{\bf 1}_{\{j_1=j_3\}}
{\bf 1}_{\{i_2=i_4\}}
{\bf 1}_{\{j_2=j_4\}}+
$$
$$
+\Biggl.
{\bf 1}_{\{i_1=i_4\}}
{\bf 1}_{\{j_1=j_4\}}
{\bf 1}_{\{i_2=i_3\}}
{\bf 1}_{\{j_2=j_3\}}\Biggr),
$$

\vspace{6mm}

$$
I_{(0010)T,t}^{(i_1 i_2 i_3 i_4)}
=\hbox{\vtop{\offinterlineskip\halign{
\hfil#\hfil\cr
{\rm l.i.m.}\cr
$\stackrel{}{{}_{p\to \infty}}$\cr
}} }
\sum_{j_1,j_2,j_3,j_4=0}^{p}
C_{j_4 j_3 j_2 j_1}^{0010}\Biggl(
\zeta_{j_1}^{(i_1)}\zeta_{j_2}^{(i_2)}\zeta_{j_3}^{(i_3)}\zeta_{j_4}^{(i_4)}
-\Biggr.
$$
$$
-
{\bf 1}_{\{i_1=i_2\}}
{\bf 1}_{\{j_1=j_2\}}
\zeta_{j_3}^{(i_3)}
\zeta_{j_4}^{(i_4)}
-
{\bf 1}_{\{i_1=i_3\}}
{\bf 1}_{\{j_1=j_3\}}
\zeta_{j_2}^{(i_2)}
\zeta_{j_4}^{(i_4)}-
$$
$$
-
{\bf 1}_{\{i_1=i_4\}}
{\bf 1}_{\{j_1=j_4\}}
\zeta_{j_2}^{(i_2)}
\zeta_{j_3}^{(i_3)}
-
{\bf 1}_{\{i_2=i_3\}}
{\bf 1}_{\{j_2=j_3\}}
\zeta_{j_1}^{(i_1)}
\zeta_{j_4}^{(i_4)}-
$$
$$
-
{\bf 1}_{\{i_2=i_4\}}
{\bf 1}_{\{j_2=j_4\}}
\zeta_{j_1}^{(i_1)}
\zeta_{j_3}^{(i_3)}
-
{\bf 1}_{\{i_3=i_4\}}
{\bf 1}_{\{j_3=j_4\}}
\zeta_{j_1}^{(i_1)}
\zeta_{j_2}^{(i_2)}+
$$
$$
+
{\bf 1}_{\{i_1=i_2\}}
{\bf 1}_{\{j_1=j_2\}}
{\bf 1}_{\{i_3=i_4\}}
{\bf 1}_{\{j_3=j_4\}}+
{\bf 1}_{\{i_1=i_3\}}
{\bf 1}_{\{j_1=j_3\}}
{\bf 1}_{\{i_2=i_4\}}
{\bf 1}_{\{j_2=j_4\}}+
$$
$$
+\Biggl.
{\bf 1}_{\{i_1=i_4\}}
{\bf 1}_{\{j_1=j_4\}}
{\bf 1}_{\{i_2=i_3\}}
{\bf 1}_{\{j_2=j_3\}}\Biggr),
$$

\vspace{6mm}

$$
I_{(0100)T,t}^{(i_1 i_2 i_3 i_4)}
=\hbox{\vtop{\offinterlineskip\halign{
\hfil#\hfil\cr
{\rm l.i.m.}\cr
$\stackrel{}{{}_{p\to \infty}}$\cr
}} }
\sum_{j_1,j_2,j_3,j_4=0}^{p}
C_{j_4 j_3 j_2 j_1}^{0100}\Biggl(
\zeta_{j_1}^{(i_1)}\zeta_{j_2}^{(i_2)}\zeta_{j_3}^{(i_3)}\zeta_{j_4}^{(i_4)}
-\Biggr.
$$
$$
-
{\bf 1}_{\{i_1=i_2\}}
{\bf 1}_{\{j_1=j_2\}}
\zeta_{j_3}^{(i_3)}
\zeta_{j_4}^{(i_4)}
-
{\bf 1}_{\{i_1=i_3\}}
{\bf 1}_{\{j_1=j_3\}}
\zeta_{j_2}^{(i_2)}
\zeta_{j_4}^{(i_4)}-
$$
$$
-
{\bf 1}_{\{i_1=i_4\}}
{\bf 1}_{\{j_1=j_4\}}
\zeta_{j_2}^{(i_2)}
\zeta_{j_3}^{(i_3)}
-
{\bf 1}_{\{i_2=i_3\}}
{\bf 1}_{\{j_2=j_3\}}
\zeta_{j_1}^{(i_1)}
\zeta_{j_4}^{(i_4)}-
$$
$$
-
{\bf 1}_{\{i_2=i_4\}}
{\bf 1}_{\{j_2=j_4\}}
\zeta_{j_1}^{(i_1)}
\zeta_{j_3}^{(i_3)}
-
{\bf 1}_{\{i_3=i_4\}}
{\bf 1}_{\{j_3=j_4\}}
\zeta_{j_1}^{(i_1)}
\zeta_{j_2}^{(i_2)}+
$$
$$
+
{\bf 1}_{\{i_1=i_2\}}
{\bf 1}_{\{j_1=j_2\}}
{\bf 1}_{\{i_3=i_4\}}
{\bf 1}_{\{j_3=j_4\}}+
{\bf 1}_{\{i_1=i_3\}}
{\bf 1}_{\{j_1=j_3\}}
{\bf 1}_{\{i_2=i_4\}}
{\bf 1}_{\{j_2=j_4\}}+
$$
$$
+\Biggl.
{\bf 1}_{\{i_1=i_4\}}
{\bf 1}_{\{j_1=j_4\}}
{\bf 1}_{\{i_2=i_3\}}
{\bf 1}_{\{j_2=j_3\}}\Biggr),
$$

\vspace{6mm}

$$
I_{(1000)T,t}^{(i_1 i_2 i_3 i_4)}
=\hbox{\vtop{\offinterlineskip\halign{
\hfil#\hfil\cr
{\rm l.i.m.}\cr
$\stackrel{}{{}_{p\to \infty}}$\cr
}} }
\sum_{j_1,j_2,j_3,j_4=0}^{p}
C_{j_4 j_3 j_2 j_1}^{1000}\Biggl(
\zeta_{j_1}^{(i_1)}\zeta_{j_2}^{(i_2)}\zeta_{j_3}^{(i_3)}\zeta_{j_4}^{(i_4)}
-\Biggr.
$$
$$
-
{\bf 1}_{\{i_1=i_2\}}
{\bf 1}_{\{j_1=j_2\}}
\zeta_{j_3}^{(i_3)}
\zeta_{j_4}^{(i_4)}
-
{\bf 1}_{\{i_1=i_3\}}
{\bf 1}_{\{j_1=j_3\}}
\zeta_{j_2}^{(i_2)}
\zeta_{j_4}^{(i_4)}-
$$
$$
-
{\bf 1}_{\{i_1=i_4\}}
{\bf 1}_{\{j_1=j_4\}}
\zeta_{j_2}^{(i_2)}
\zeta_{j_3}^{(i_3)}
-
{\bf 1}_{\{i_2=i_3\}}
{\bf 1}_{\{j_2=j_3\}}
\zeta_{j_1}^{(i_1)}
\zeta_{j_4}^{(i_4)}-
$$
$$
-
{\bf 1}_{\{i_2=i_4\}}
{\bf 1}_{\{j_2=j_4\}}
\zeta_{j_1}^{(i_1)}
\zeta_{j_3}^{(i_3)}
-
{\bf 1}_{\{i_3=i_4\}}
{\bf 1}_{\{j_3=j_4\}}
\zeta_{j_1}^{(i_1)}
\zeta_{j_2}^{(i_2)}+
$$
$$
+
{\bf 1}_{\{i_1=i_2\}}
{\bf 1}_{\{j_1=j_2\}}
{\bf 1}_{\{i_3=i_4\}}
{\bf 1}_{\{j_3=j_4\}}+
{\bf 1}_{\{i_1=i_3\}}
{\bf 1}_{\{j_1=j_3\}}
{\bf 1}_{\{i_2=i_4\}}
{\bf 1}_{\{j_2=j_4\}}+
$$
$$
+\Biggl.
{\bf 1}_{\{i_1=i_4\}}
{\bf 1}_{\{j_1=j_4\}}
{\bf 1}_{\{i_2=i_3\}}
{\bf 1}_{\{j_2=j_3\}}\Biggr),
$$

\vspace{3mm}
\noindent
where

$$
C_{j_4j_3j_2j_1}^{0001}
=\frac{\sqrt{(2j_1+1)(2j_2+1)(2j_3+1)(2j_4+1)}}{32}(T-t)^{3}\bar
C_{j_4j_3j_2j_1}^{0001},
$$

\vspace{4mm}
$$
C_{j_3j_2j_1}^{0010}
=\frac{\sqrt{(2j_1+1)(2j_2+1)(2j_3+1)(2j_4+1)}}{32}(T-t)^{3}\bar
C_{j_4j_3j_2j_1}^{0010},
$$

\vspace{4mm}
$$
C_{j_4j_3j_2j_1}^{0100}=
\frac{\sqrt{(2j_1+1)(2j_2+1)(2j_3+1)(2j_4+1)}}{32}(T-t)^{3}\bar
C_{j_3j_2j_1}^{0100},
$$

\vspace{4mm}
$$
C_{j_4j_3j_2j_1}^{1000}
=\frac{\sqrt{(2j_1+1)(2j_2+1)(2j_3+1)(2j_4+1)}}{32}(T-t)^{3}\bar
C_{j_4j_3j_2j_1}^{1000},
$$

\vspace{4mm}
$$
\bar C_{j_4j_3j_2j_1}^{1000}=-
\int\limits_{-1}^{1}P_{j_4}(u)
\int\limits_{-1}^{u}P_{j_3}(z)
\int\limits_{-1}^{z}P_{j_2}(y)
\int\limits_{-1}^{y}
P_{j_1}(x)(x+1)dx dy dz du,
$$

\vspace{3mm}
$$
\bar C_{j_4j_3j_2j_1}^{0100}=-
\int\limits_{-1}^{1}P_{j_4}(u)
\int\limits_{-1}^{u}P_{j_3}(z)
\int\limits_{-1}^{z}P_{j_2}(y)(y+1)
\int\limits_{-1}^{y}
P_{j_1}(x)dx dy dz du,
$$

\vspace{3mm}
$$
\bar C_{j_4j_3j_2j_1}^{0010}=-
\int\limits_{-1}^{1}P_{j_4}(u)
\int\limits_{-1}^{u}P_{j_3}(z)(z+1)
\int\limits_{-1}^{z}P_{j_2}(y)
\int\limits_{-1}^{y}
P_{j_1}(x)dx dy dz du,
$$

\vspace{3mm}
$$
\bar C_{j_4j_3j_2j_1}^{0001}=-
\int\limits_{-1}^{1}P_{j_4}(u)(u+1)
\int\limits_{-1}^{u}P_{j_3}(z)
\int\limits_{-1}^{z}P_{j_2}(y)
\int\limits_{-1}^{y}
P_{j_1}(x)dx dy dz du;
$$

\vspace{6mm}

$$
I_{(000000)T,t}^{*(i_1 i_2 i_3 i_4 i_5 i_6)}=
\hbox{\vtop{\offinterlineskip\halign{
\hfil#\hfil\cr
{\rm l.i.m.}\cr
$\stackrel{}{{}_{p\to \infty}}$\cr
}} }
\sum\limits_{j_1, j_2, j_3, j_4, j_5, j_6=0}^{p}
C_{j_6j_5j_4 j_3 j_2 j_1}
\zeta_{j_1}^{(i_1)}\zeta_{j_2}^{(i_2)}\zeta_{j_3}^{(i_3)}
\zeta_{j_4}^{(i_4)}\zeta_{j_5}^{(i_5)}\zeta_{j_6}^{(i_6)},
$$

\vspace{7mm}

$$
I_{(000000)T,t}^{(i_1 i_2 i_3 i_4 i_5 i_6)}
=\hbox{\vtop{\offinterlineskip\halign{
\hfil#\hfil\cr
{\rm l.i.m.}\cr
$\stackrel{}{{}_{p\to \infty}}$\cr
}} }\sum_{j_1,j_2,j_3,j_4,j_5,j_6=0}^{p}
C_{j_6 j_5 j_4 j_3 j_2 j_1}\Biggl(
\prod_{l=1}^6
\zeta_{j_l}^{(i_l)}
-\Biggr.
$$
$$
-
{\bf 1}_{\{j_1=j_6\}}
{\bf 1}_{\{i_1=i_6\}}
\zeta_{j_2}^{(i_2)}
\zeta_{j_3}^{(i_3)}
\zeta_{j_4}^{(i_4)}
\zeta_{j_5}^{(i_5)}-
{\bf 1}_{\{j_2=j_6\}}
{\bf 1}_{\{i_2=i_6\}}
\zeta_{j_1}^{(i_1)}
\zeta_{j_3}^{(i_3)}
\zeta_{j_4}^{(i_4)}
\zeta_{j_5}^{(i_5)}-
$$
$$
-
{\bf 1}_{\{j_3=j_6\}}
{\bf 1}_{\{i_3=i_6\}}
\zeta_{j_1}^{(i_1)}
\zeta_{j_2}^{(i_2)}
\zeta_{j_4}^{(i_4)}
\zeta_{j_5}^{(i_5)}-
{\bf 1}_{\{j_4=j_6\}}
{\bf 1}_{\{i_4=i_6\}}
\zeta_{j_1}^{(i_1)}
\zeta_{j_2}^{(i_2)}
\zeta_{j_3}^{(i_3)}
\zeta_{j_5}^{(i_5)}-
$$
$$
-
{\bf 1}_{\{j_5=j_6\}}
{\bf 1}_{\{i_5=i_6\}}
\zeta_{j_1}^{(i_1)}
\zeta_{j_2}^{(i_2)}
\zeta_{j_3}^{(i_3)}
\zeta_{j_4}^{(i_4)}-
{\bf 1}_{\{j_1=j_2\}}
{\bf 1}_{\{i_1=i_2\}}
\zeta_{j_3}^{(i_3)}
\zeta_{j_4}^{(i_4)}
\zeta_{j_5}^{(i_5)}
\zeta_{j_6}^{(i_6)}-
$$
$$
-
{\bf 1}_{\{j_1=j_3\}}
{\bf 1}_{\{i_1=i_3\}}
\zeta_{j_2}^{(i_2)}
\zeta_{j_4}^{(i_4)}
\zeta_{j_5}^{(i_5)}
\zeta_{j_6}^{(i_6)}-
{\bf 1}_{\{j_1=j_4\}}
{\bf 1}_{\{i_1=i_4\}}
\zeta_{j_2}^{(i_2)}
\zeta_{j_3}^{(i_3)}
\zeta_{j_5}^{(i_5)}
\zeta_{j_6}^{(i_6)}-
$$
$$
-
{\bf 1}_{\{j_1=j_5\}}
{\bf 1}_{\{i_1=i_5\}}
\zeta_{j_2}^{(i_2)}
\zeta_{j_3}^{(i_3)}
\zeta_{j_4}^{(i_4)}
\zeta_{j_6}^{(i_6)}-
{\bf 1}_{\{j_2=j_3\}}
{\bf 1}_{\{i_2=i_3\}}
\zeta_{j_1}^{(i_1)}
\zeta_{j_4}^{(i_4)}
\zeta_{j_5}^{(i_5)}
\zeta_{j_6}^{(i_6)}-
$$
$$
-
{\bf 1}_{\{j_2=j_4\}}
{\bf 1}_{\{i_2=i_4\}}
\zeta_{j_1}^{(i_1)}
\zeta_{j_3}^{(i_3)}
\zeta_{j_5}^{(i_5)}
\zeta_{j_6}^{(i_6)}-
{\bf 1}_{\{j_2=j_5\}}
{\bf 1}_{\{i_2=i_5\}}
\zeta_{j_1}^{(i_1)}
\zeta_{j_3}^{(i_3)}
\zeta_{j_4}^{(i_4)}
\zeta_{j_6}^{(i_6)}-
$$
$$
-
{\bf 1}_{\{j_3=j_4\}}
{\bf 1}_{\{i_3=i_4\}}
\zeta_{j_1}^{(i_1)}
\zeta_{j_2}^{(i_2)}
\zeta_{j_5}^{(i_5)}
\zeta_{j_6}^{(i_6)}-
{\bf 1}_{\{j_3=j_5\}}
{\bf 1}_{\{i_3=i_5\}}
\zeta_{j_1}^{(i_1)}
\zeta_{j_2}^{(i_2)}
\zeta_{j_4}^{(i_4)}
\zeta_{j_6}^{(i_6)}-
$$
$$
-
{\bf 1}_{\{j_4=j_5\}}
{\bf 1}_{\{i_4=i_5\}}
\zeta_{j_1}^{(i_1)}
\zeta_{j_2}^{(i_2)}
\zeta_{j_3}^{(i_3)}
\zeta_{j_6}^{(i_6)}+
$$
$$
+
{\bf 1}_{\{j_1=j_2\}}
{\bf 1}_{\{i_1=i_2\}}
{\bf 1}_{\{j_3=j_4\}}
{\bf 1}_{\{i_3=i_4\}}
\zeta_{j_5}^{(i_5)}
\zeta_{j_6}^{(i_6)}
+
{\bf 1}_{\{j_1=j_2\}}
{\bf 1}_{\{i_1=i_2\}}
{\bf 1}_{\{j_3=j_5\}}
{\bf 1}_{\{i_3=i_5\}}
\zeta_{j_4}^{(i_4)}
\zeta_{j_6}^{(i_6)}+
$$
$$
+
{\bf 1}_{\{j_1=j_2\}}
{\bf 1}_{\{i_1=i_2\}}
{\bf 1}_{\{j_4=j_5\}}
{\bf 1}_{\{i_4=i_5\}}
\zeta_{j_3}^{(i_3)}
\zeta_{j_6}^{(i_6)}
+
{\bf 1}_{\{j_1=j_3\}}
{\bf 1}_{\{i_1=i_3\}}
{\bf 1}_{\{j_2=j_4\}}
{\bf 1}_{\{i_2=i_4\}}
\zeta_{j_5}^{(i_5)}
\zeta_{j_6}^{(i_6)}+
$$
$$
+
{\bf 1}_{\{j_1=j_3\}}
{\bf 1}_{\{i_1=i_3\}}
{\bf 1}_{\{j_2=j_5\}}
{\bf 1}_{\{i_2=i_5\}}
\zeta_{j_4}^{(i_4)}
\zeta_{j_6}^{(i_6)}
+
{\bf 1}_{\{j_1=j_3\}}
{\bf 1}_{\{i_1=i_3\}}
{\bf 1}_{\{j_4=j_5\}}
{\bf 1}_{\{i_4=i_5\}}
\zeta_{j_2}^{(i_2)}
\zeta_{j_6}^{(i_6)}+
$$
$$
+
{\bf 1}_{\{j_1=j_4\}}
{\bf 1}_{\{i_1=i_4\}}
{\bf 1}_{\{j_2=j_3\}}
{\bf 1}_{\{i_2=i_3\}}
\zeta_{j_5}^{(i_5)}
\zeta_{j_6}^{(i_6)}
+
{\bf 1}_{\{j_1=j_4\}}
{\bf 1}_{\{i_1=i_4\}}
{\bf 1}_{\{j_2=j_5\}}
{\bf 1}_{\{i_2=i_5\}}
\zeta_{j_3}^{(i_3)}
\zeta_{j_6}^{(i_6)}+
$$
$$
+
{\bf 1}_{\{j_1=j_4\}}
{\bf 1}_{\{i_1=i_4\}}
{\bf 1}_{\{j_3=j_5\}}
{\bf 1}_{\{i_3=i_5\}}
\zeta_{j_2}^{(i_2)}
\zeta_{j_6}^{(i_6)}
+
{\bf 1}_{\{j_1=j_5\}}
{\bf 1}_{\{i_1=i_5\}}
{\bf 1}_{\{j_2=j_3\}}
{\bf 1}_{\{i_2=i_3\}}
\zeta_{j_4}^{(i_4)}
\zeta_{j_6}^{(i_6)}+
$$
$$
+
{\bf 1}_{\{j_1=j_5\}}
{\bf 1}_{\{i_1=i_5\}}
{\bf 1}_{\{j_2=j_4\}}
{\bf 1}_{\{i_2=i_4\}}
\zeta_{j_3}^{(i_3)}
\zeta_{j_6}^{(i_6)}
+
{\bf 1}_{\{j_1=j_5\}}
{\bf 1}_{\{i_1=i_5\}}
{\bf 1}_{\{j_3=j_4\}}
{\bf 1}_{\{i_3=i_4\}}
\zeta_{j_2}^{(i_2)}
\zeta_{j_6}^{(i_6)}+
$$
$$
+
{\bf 1}_{\{j_2=j_3\}}
{\bf 1}_{\{i_2=i_3\}}
{\bf 1}_{\{j_4=j_5\}}
{\bf 1}_{\{i_4=i_5\}}
\zeta_{j_1}^{(i_1)}
\zeta_{j_6}^{(i_6)}
+
{\bf 1}_{\{j_2=j_4\}}
{\bf 1}_{\{i_2=i_4\}}
{\bf 1}_{\{j_3=j_5\}}
{\bf 1}_{\{i_3=i_5\}}
\zeta_{j_1}^{(i_1)}
\zeta_{j_6}^{(i_6)}+
$$
$$
+
{\bf 1}_{\{j_2=j_5\}}
{\bf 1}_{\{i_2=i_5\}}
{\bf 1}_{\{j_3=j_4\}}
{\bf 1}_{\{i_3=i_4\}}
\zeta_{j_1}^{(i_1)}
\zeta_{j_6}^{(i_6)}
+
{\bf 1}_{\{j_6=j_1\}}
{\bf 1}_{\{i_6=i_1\}}
{\bf 1}_{\{j_3=j_4\}}
{\bf 1}_{\{i_3=i_4\}}
\zeta_{j_2}^{(i_2)}
\zeta_{j_5}^{(i_5)}+
$$
$$
+
{\bf 1}_{\{j_6=j_1\}}
{\bf 1}_{\{i_6=i_1\}}
{\bf 1}_{\{j_3=j_5\}}
{\bf 1}_{\{i_3=i_5\}}
\zeta_{j_2}^{(i_2)}
\zeta_{j_4}^{(i_4)}
+
{\bf 1}_{\{j_6=j_1\}}
{\bf 1}_{\{i_6=i_1\}}
{\bf 1}_{\{j_2=j_5\}}
{\bf 1}_{\{i_2=i_5\}}
\zeta_{j_3}^{(i_3)}
\zeta_{j_4}^{(i_4)}+
$$
$$
+
{\bf 1}_{\{j_6=j_1\}}
{\bf 1}_{\{i_6=i_1\}}
{\bf 1}_{\{j_2=j_4\}}
{\bf 1}_{\{i_2=i_4\}}
\zeta_{j_3}^{(i_3)}
\zeta_{j_5}^{(i_5)}
+
{\bf 1}_{\{j_6=j_1\}}
{\bf 1}_{\{i_6=i_1\}}
{\bf 1}_{\{j_4=j_5\}}
{\bf 1}_{\{i_4=i_5\}}
\zeta_{j_2}^{(i_2)}
\zeta_{j_3}^{(i_3)}+
$$
$$
+
{\bf 1}_{\{j_6=j_1\}}
{\bf 1}_{\{i_6=i_1\}}
{\bf 1}_{\{j_2=j_3\}}
{\bf 1}_{\{i_2=i_3\}}
\zeta_{j_4}^{(i_4)}
\zeta_{j_5}^{(i_5)}
+
{\bf 1}_{\{j_6=j_2\}}
{\bf 1}_{\{i_6=i_2\}}
{\bf 1}_{\{j_3=j_5\}}
{\bf 1}_{\{i_3=i_5\}}
\zeta_{j_1}^{(i_1)}
\zeta_{j_4}^{(i_4)}+
$$
$$
+
{\bf 1}_{\{j_6=j_2\}}
{\bf 1}_{\{i_6=i_2\}}
{\bf 1}_{\{j_4=j_5\}}
{\bf 1}_{\{i_4=i_5\}}
\zeta_{j_1}^{(i_1)}
\zeta_{j_3}^{(i_3)}
+
{\bf 1}_{\{j_6=j_2\}}
{\bf 1}_{\{i_6=i_2\}}
{\bf 1}_{\{j_3=j_4\}}
{\bf 1}_{\{i_3=i_4\}}
\zeta_{j_1}^{(i_1)}
\zeta_{j_5}^{(i_5)}+
$$
$$
+
{\bf 1}_{\{j_6=j_2\}}
{\bf 1}_{\{i_6=i_2\}}
{\bf 1}_{\{j_1=j_5\}}
{\bf 1}_{\{i_1=i_5\}}
\zeta_{j_3}^{(i_3)}
\zeta_{j_4}^{(i_4)}
+
{\bf 1}_{\{j_6=j_2\}}
{\bf 1}_{\{i_6=i_2\}}
{\bf 1}_{\{j_1=j_4\}}
{\bf 1}_{\{i_1=i_4\}}
\zeta_{j_3}^{(i_3)}
\zeta_{j_5}^{(i_5)}+
$$
$$
+
{\bf 1}_{\{j_6=j_2\}}
{\bf 1}_{\{i_6=i_2\}}
{\bf 1}_{\{j_1=j_3\}}
{\bf 1}_{\{i_1=i_3\}}
\zeta_{j_4}^{(i_4)}
\zeta_{j_5}^{(i_5)}
+
{\bf 1}_{\{j_6=j_3\}}
{\bf 1}_{\{i_6=i_3\}}
{\bf 1}_{\{j_2=j_5\}}
{\bf 1}_{\{i_2=i_5\}}
\zeta_{j_1}^{(i_1)}
\zeta_{j_4}^{(i_4)}+
$$
$$
+
{\bf 1}_{\{j_6=j_3\}}
{\bf 1}_{\{i_6=i_3\}}
{\bf 1}_{\{j_4=j_5\}}
{\bf 1}_{\{i_4=i_5\}}
\zeta_{j_1}^{(i_1)}
\zeta_{j_2}^{(i_2)}
+
{\bf 1}_{\{j_6=j_3\}}
{\bf 1}_{\{i_6=i_3\}}
{\bf 1}_{\{j_2=j_4\}}
{\bf 1}_{\{i_2=i_4\}}
\zeta_{j_1}^{(i_1)}
\zeta_{j_5}^{(i_5)}+
$$
$$
+
{\bf 1}_{\{j_6=j_3\}}
{\bf 1}_{\{i_6=i_3\}}
{\bf 1}_{\{j_1=j_5\}}
{\bf 1}_{\{i_1=i_5\}}
\zeta_{j_2}^{(i_2)}
\zeta_{j_4}^{(i_4)}
+
{\bf 1}_{\{j_6=j_3\}}
{\bf 1}_{\{i_6=i_3\}}
{\bf 1}_{\{j_1=j_4\}}
{\bf 1}_{\{i_1=i_4\}}
\zeta_{j_2}^{(i_2)}
\zeta_{j_5}^{(i_5)}+
$$
$$
+
{\bf 1}_{\{j_6=j_3\}}
{\bf 1}_{\{i_6=i_3\}}
{\bf 1}_{\{j_1=j_2\}}
{\bf 1}_{\{i_1=i_2\}}
\zeta_{j_4}^{(i_4)}
\zeta_{j_5}^{(i_5)}
+
{\bf 1}_{\{j_6=j_4\}}
{\bf 1}_{\{i_6=i_4\}}
{\bf 1}_{\{j_3=j_5\}}
{\bf 1}_{\{i_3=i_5\}}
\zeta_{j_1}^{(i_1)}
\zeta_{j_2}^{(i_2)}+
$$
$$
+
{\bf 1}_{\{j_6=j_4\}}
{\bf 1}_{\{i_6=i_4\}}
{\bf 1}_{\{j_2=j_5\}}
{\bf 1}_{\{i_2=i_5\}}
\zeta_{j_1}^{(i_1)}
\zeta_{j_3}^{(i_3)}
+
{\bf 1}_{\{j_6=j_4\}}
{\bf 1}_{\{i_6=i_4\}}
{\bf 1}_{\{j_2=j_3\}}
{\bf 1}_{\{i_2=i_3\}}
\zeta_{j_1}^{(i_1)}
\zeta_{j_5}^{(i_5)}+
$$
$$
+
{\bf 1}_{\{j_6=j_4\}}
{\bf 1}_{\{i_6=i_4\}}
{\bf 1}_{\{j_1=j_5\}}
{\bf 1}_{\{i_1=i_5\}}
\zeta_{j_2}^{(i_2)}
\zeta_{j_3}^{(i_3)}
+
{\bf 1}_{\{j_6=j_4\}}
{\bf 1}_{\{i_6=i_4\}}
{\bf 1}_{\{j_1=j_3\}}
{\bf 1}_{\{i_1=i_3\}}
\zeta_{j_2}^{(i_2)}
\zeta_{j_5}^{(i_5)}+
$$
$$
+
{\bf 1}_{\{j_6=j_4\}}
{\bf 1}_{\{i_6=i_4\}}
{\bf 1}_{\{j_1=j_2\}}
{\bf 1}_{\{i_1=i_2\}}
\zeta_{j_3}^{(i_3)}
\zeta_{j_5}^{(i_5)}
+
{\bf 1}_{\{j_6=j_5\}}
{\bf 1}_{\{i_6=i_5\}}
{\bf 1}_{\{j_3=j_4\}}
{\bf 1}_{\{i_3=i_4\}}
\zeta_{j_1}^{(i_1)}
\zeta_{j_2}^{(i_2)}+
$$
$$
+
{\bf 1}_{\{j_6=j_5\}}
{\bf 1}_{\{i_6=i_5\}}
{\bf 1}_{\{j_2=j_4\}}
{\bf 1}_{\{i_2=i_4\}}
\zeta_{j_1}^{(i_1)}
\zeta_{j_3}^{(i_3)}
+
{\bf 1}_{\{j_6=j_5\}}
{\bf 1}_{\{i_6=i_5\}}
{\bf 1}_{\{j_2=j_3\}}
{\bf 1}_{\{i_2=i_3\}}
\zeta_{j_1}^{(i_1)}
\zeta_{j_4}^{(i_4)}+
$$
$$
+
{\bf 1}_{\{j_6=j_5\}}
{\bf 1}_{\{i_6=i_5\}}
{\bf 1}_{\{j_1=j_4\}}
{\bf 1}_{\{i_1=i_4\}}
\zeta_{j_2}^{(i_2)}
\zeta_{j_3}^{(i_3)}
+
{\bf 1}_{\{j_6=j_5\}}
{\bf 1}_{\{i_6=i_5\}}
{\bf 1}_{\{j_1=j_3\}}
{\bf 1}_{\{i_1=i_3\}}
\zeta_{j_2}^{(i_2)}
\zeta_{j_4}^{(i_4)}+
$$
$$
+
{\bf 1}_{\{j_6=j_5\}}
{\bf 1}_{\{i_6=i_5\}}
{\bf 1}_{\{j_1=j_2\}}
{\bf 1}_{\{i_1=i_2\}}
\zeta_{j_3}^{(i_3)}
\zeta_{j_4}^{(i_4)}-
$$
$$
-
{\bf 1}_{\{j_6=j_1\}}
{\bf 1}_{\{i_6=i_1\}}
{\bf 1}_{\{j_2=j_5\}}
{\bf 1}_{\{i_2=i_5\}}
{\bf 1}_{\{j_3=j_4\}}
{\bf 1}_{\{i_3=i_4\}}-
$$
$$
-
{\bf 1}_{\{j_6=j_1\}}
{\bf 1}_{\{i_6=i_1\}}
{\bf 1}_{\{j_2=j_4\}}
{\bf 1}_{\{i_2=i_4\}}
{\bf 1}_{\{j_3=j_5\}}
{\bf 1}_{\{i_3=i_5\}}-
$$
$$
-
{\bf 1}_{\{j_6=j_1\}}
{\bf 1}_{\{i_6=i_1\}}
{\bf 1}_{\{j_2=j_3\}}
{\bf 1}_{\{i_2=i_3\}}
{\bf 1}_{\{j_4=j_5\}}
{\bf 1}_{\{i_4=i_5\}}-
$$
$$
-               
{\bf 1}_{\{j_6=j_2\}}
{\bf 1}_{\{i_6=i_2\}}
{\bf 1}_{\{j_1=j_5\}}
{\bf 1}_{\{i_1=i_5\}}
{\bf 1}_{\{j_3=j_4\}}
{\bf 1}_{\{i_3=i_4\}}-
$$
$$
-
{\bf 1}_{\{j_6=j_2\}}
{\bf 1}_{\{i_6=i_2\}}
{\bf 1}_{\{j_1=j_4\}}
{\bf 1}_{\{i_1=i_4\}}
{\bf 1}_{\{j_3=j_5\}}
{\bf 1}_{\{i_3=i_5\}}-
$$
$$
-
{\bf 1}_{\{j_6=j_2\}}
{\bf 1}_{\{i_6=i_2\}}
{\bf 1}_{\{j_1=j_3\}}
{\bf 1}_{\{i_1=i_3\}}
{\bf 1}_{\{j_4=j_5\}}
{\bf 1}_{\{i_4=i_5\}}-
$$
$$
-
{\bf 1}_{\{j_6=j_3\}}
{\bf 1}_{\{i_6=i_3\}}
{\bf 1}_{\{j_1=j_5\}}
{\bf 1}_{\{i_1=i_5\}}
{\bf 1}_{\{j_2=j_4\}}
{\bf 1}_{\{i_2=i_4\}}-
$$
$$
-
{\bf 1}_{\{j_6=j_3\}}
{\bf 1}_{\{i_6=i_3\}}
{\bf 1}_{\{j_1=j_4\}}
{\bf 1}_{\{i_1=i_4\}}
{\bf 1}_{\{j_2=j_5\}}
{\bf 1}_{\{i_2=i_5\}}-
$$
$$
-
{\bf 1}_{\{j_3=j_6\}}
{\bf 1}_{\{i_3=i_6\}}
{\bf 1}_{\{j_1=j_2\}}
{\bf 1}_{\{i_1=i_2\}}
{\bf 1}_{\{j_4=j_5\}}
{\bf 1}_{\{i_4=i_5\}}-
$$
$$
-
{\bf 1}_{\{j_6=j_4\}}
{\bf 1}_{\{i_6=i_4\}}
{\bf 1}_{\{j_1=j_5\}}
{\bf 1}_{\{i_1=i_5\}}
{\bf 1}_{\{j_2=j_3\}}
{\bf 1}_{\{i_2=i_3\}}-
$$
$$
-
{\bf 1}_{\{j_6=j_4\}}
{\bf 1}_{\{i_6=i_4\}}
{\bf 1}_{\{j_1=j_3\}}
{\bf 1}_{\{i_1=i_3\}}
{\bf 1}_{\{j_2=j_5\}}
{\bf 1}_{\{i_2=i_5\}}-
$$
$$
-
{\bf 1}_{\{j_6=j_4\}}
{\bf 1}_{\{i_6=i_4\}}
{\bf 1}_{\{j_1=j_2\}}
{\bf 1}_{\{i_1=i_2\}}
{\bf 1}_{\{j_3=j_5\}}
{\bf 1}_{\{i_3=i_5\}}-
$$
$$
-
{\bf 1}_{\{j_6=j_5\}}
{\bf 1}_{\{i_6=i_5\}}
{\bf 1}_{\{j_1=j_4\}}
{\bf 1}_{\{i_1=i_4\}}
{\bf 1}_{\{j_2=j_3\}}
{\bf 1}_{\{i_2=i_3\}}-
$$
$$
-
{\bf 1}_{\{j_6=j_5\}}
{\bf 1}_{\{i_6=i_5\}}
{\bf 1}_{\{j_1=j_2\}}
{\bf 1}_{\{i_1=i_2\}}
{\bf 1}_{\{j_3=j_4\}}
{\bf 1}_{\{i_3=i_4\}}-
$$
$$
\Biggl.-
{\bf 1}_{\{j_6=j_5\}}
{\bf 1}_{\{i_6=i_5\}}
{\bf 1}_{\{j_1=j_3\}}
{\bf 1}_{\{i_1=i_3\}}
{\bf 1}_{\{j_2=j_4\}}
{\bf 1}_{\{i_2=i_4\}}\Biggr),
$$

\vspace{6mm}

$$         
I_{(000000)T,t}^{(i_1i_1i_1i_1i_1i_1)}=
\frac{1}{720}(T-t)^{3}
\left(\left(\zeta_0^{(i_1)}\right)^6-
15\left(\zeta_0^{(i_1)}\right)^4+45\left(\zeta_0^{(i_1)}\right)^2-
15\right)\ \ \ 
\hbox{w.\ p.\ 1},
$$

\vspace{3mm}
$$
I_{(000000)T,t}^{*(i_1i_1i_1i_1i_1i_1)}=
\frac{1}{720}(T-t)^{3}\left(\zeta_0^{(i_1)}\right)^6\ \ \ \hbox{w.\ p.\ 1},
$$

\vspace{3mm}
\noindent
where

$$
C_{j_6j_5j_4 j_3 j_2 j_1}
=\frac{\sqrt{(2j_1+1)(2j_2+1)(2j_3+1)
(2j_4+1)(2j_5+1)(2j_6+1)}}{64}(T-t)^{3}\bar
C_{j_6j_5j_4 j_3 j_2 j_1},
$$

\vspace{3mm}

$$
\bar C_{j_6j_5j_4 j_3 j_2 j_1}=
\int\limits_{-1}^{1}P_{j_6}(w)
\int\limits_{-1}^{w}P_{j_5}(v)
\int\limits_{-1}^{v}P_{j_4}(u)
\int\limits_{-1}^{u}P_{j_3}(z)
\int\limits_{-1}^{z}P_{j_2}(y)
\int\limits_{-1}^{y}
P_{j_1}(x)dx dy dz du dv dw.
$$

\vspace{5mm}

Consider the approximation $I_{(00)T,t}^{*(i_1 i_2)q}$
of the iterated stochastic 
integral $I_{(00)T,t}^{*(i_1 i_2)}$ obtained from 
(\ref{4004}) by replacing $\infty$ on $q$.

It is easy to prove that \cite{3} (1997)

\vspace{3mm}
$$
{\sf M}\biggl\{\left(I_{(00)T,t}^{*(i_1 i_2)}-
I_{(00)T,t}^{*(i_1 i_2)q}
\right)^2\biggr\}
=\frac{(T-t)^2}{2}\left(\frac{1}{2}-\sum_{i=1}^q
\frac{1}{4i^2-1}\right)\ \ \ (i_1\ne i_2).
$$

\vspace{6mm}

Further, using Theorems 10, 11,  we obtain for $i_1\ne i_2$

\vspace{3mm}
$$
{\sf M}\biggl\{\left(I_{(10)T,t}^{*(i_1 i_2)}-I_{(10)T,t}^{*(i_1 i_2)q}
\right)^2\biggr\}=
{\sf M}\biggl\{\left(I_{(01)T,t}^{*(i_1 i_2)}-
I_{(01)T,t}^{*(i_1 i_2)q}\right)^2\biggr\}=
$$

\vspace{3mm}
$$
=\frac{(T-t)^4}{16}\left(\frac{5}{9}-
2\sum_{i=2}^q\frac{1}{4i^2-1}-
\sum_{i=1}^q
\frac{1}{(2i-1)^2(2i+3)^2}
-\sum_{i=0}^q\frac{(i+2)^2+(i+1)^2}{(2i+1)(2i+5)(2i+3)^2}
\right).
$$

\vspace{8mm}

For the case $i_1=i_2$
using Theorems 10, 11,  we have

\vspace{3mm}
$$
{\sf M}\biggl\{\left(I_{(10)T,t}^{(i_1 i_1)}-
I_{(10)T,t}^{(i_1 i_1)q}
\right)^2\biggr\}=
{\sf M}\biggl\{\left(I_{(01)T,t}^{(i_1 i_1)}-
I_{(01)T,t}^{(i_1 i_1)q}\right)^2\biggr\}=
$$

\vspace{3mm}
\begin{equation}
\label{2007ura1}
=\frac{(T-t)^4}{16}\left(\frac{1}{9}-
\sum_{i=0}^{q}
\frac{1}{(2i+1)(2i+5)(2i+3)^2}
-2\sum_{i=1}^{q}
\frac{1}{(2i-1)^2(2i+3)^2}\right).
\end{equation}

\vspace{8mm}

On the basis of 
the presented 
expansions of 
iterated stochastic integrals we 
can see that increasing of multiplicities of these integrals 
or degree indexes of their weight functions 
leads
to a noticeable complication of formulas 
for the mentioned expansions. 

However, increasing of the mentioned parameters leads to increasing 
of orders of smallness with respect to $T-t$ in the mean-square sense 
for iterated stochastic integrals. As a result, this feature
leads to a sharp decrease  
of member 
quantities
in expansions of iterated stochastic 
integrals, which are required for achieving the acceptable accuracy
of approximation. In the context of it, let us consider the approach 
to approximation of iterated stochastic integrals, which 
provides a possibility to obtain the mean-square approximations of 
the required accuracy without the using of
complex expansions.

Let us analyze the following approximation of triple stochastic integral 
on the base of (\ref{zzz1})

\vspace{3mm}
$$
I_{(000)T,t}^{(i_1i_2i_3)q_1}
=\sum_{j_1,j_2,j_3=0}^{q_1}
C_{j_3j_2j_1}\Biggl(
\zeta_{j_1}^{(i_1)}\zeta_{j_2}^{(i_2)}\zeta_{j_3}^{(i_3)}
-{\bf 1}_{\{i_1=i_2\}}
{\bf 1}_{\{j_1=j_2\}}
\zeta_{j_3}^{(i_3)}-
\Biggr.
$$

\vspace{1mm}
\begin{equation}
\label{sad001}
\Biggl.
-{\bf 1}_{\{i_2=i_3\}}
{\bf 1}_{\{j_2=j_3\}}
\zeta_{j_1}^{(i_1)}-
{\bf 1}_{\{i_1=i_3\}}
{\bf 1}_{\{j_1=j_3\}}
\zeta_{j_2}^{(i_2)}\Biggr),
\end{equation}

\vspace{5mm}
\noindent
where $C_{j_3j_2j_1}$ is defined by (\ref{zzz2}), (\ref{zzz3}).

In particular, from
(\ref{sad001}) for $i_1\ne i_2$, 
$i_2\ne i_3$, $i_1\ne i_3$
we obtain

\vspace{2mm}
\begin{equation}
\label{38}
I_{(000)T,t}^{(i_1i_2i_3)q_1}=
\sum_{j_1,j_2,j_3=0}^{q_1}
C_{j_3j_2j_1}
\zeta_{j_1}^{(i_1)}\zeta_{j_2}^{(i_2)}\zeta_{j_3}^{(i_3)}.
\end{equation}

\vspace{4mm}

Using (\ref{zzz0}), (\ref{zzz4}), (\ref{zzz5})--(\ref{zzz7}), we get

\vspace{2mm}
$$
{\sf M}\left\{\left(
I_{(000)T,t}^{(i_1i_2 i_3)}-
I_{(000)T,t}^{(i_1i_2 i_3)q_1}\right)^2\right\}=
$$

\begin{equation}
\label{39}
=
\frac{(T-t)^{3}}{6}
-\sum_{j_1,j_2,j_3=0}^{q_1}
C_{j_3j_2j_1}^2\ \ \ (i_1\ne i_2, i_1\ne i_3, i_2\ne i_3),
\end{equation}

\vspace{7mm}

$$
{\sf M}\left\{\left(
I_{(000)T,t}^{(i_1i_2 i_3)}-
I_{(000)T,t}^{(i_1i_2 i_3)q_1}\right)^2\right\}=
$$

\begin{equation}
\label{39a}
=
\frac{(T-t)^{3}}{6}-\sum_{j_1,j_2,j_3=0}^{q_1}
C_{j_3j_2j_1}^2
-\sum_{j_1,j_2,j_3=0}^{q_1}
C_{j_2j_3j_1}C_{j_3j_2j_1}\ \ \ (i_1\ne i_2=i_3),
\end{equation}

\vspace{7mm}

$$
{\sf M}\left\{\left(
I_{(000)T,t}^{(i_1i_2 i_3)}-
I_{(000)T,t}^{(i_1i_2 i_3)q_1}\right)^2\right\}=
$$

\begin{equation}
\label{39b}
=
\frac{(T-t)^{3}}{6}-\sum_{j_1,j_2,j_3=0}^{q_1}
C_{j_3j_2j_1}^2
-\sum_{j_1,j_2,j_3=0}^{q_1}
C_{j_3j_2j_1}C_{j_1j_2j_3}\ \ \ (i_1=i_3\ne i_2),
\end{equation}

\vspace{7mm}

$$
{\sf M}\left\{\left(
I_{(000)T,t}^{(i_1i_2 i_3)}-
I_{(000)T,t}^{(i_1i_2 i_3)q_1}\right)^2\right\}=
$$

\begin{equation}
\label{39c}
=
\frac{(T-t)^{3}}{6}-\sum_{j_1,j_2,j_3=0}^{q_1}
C_{j_3j_2j_1}^2
-\sum_{j_1,j_2,j_3=0}^{q_1}
C_{j_3j_1j_2}C_{j_3j_2j_1}\ \ \ (i_1=i_2\ne i_3),
\end{equation}

\vspace{7mm}

\begin{equation}
\label{leto1041}
{\sf M}\left\{\left(
I_{(000)T,t}^{(i_1i_2 i_3)}-
I_{(000)T,t}^{(i_1i_2 i_3)q_1}\right)^2\right\}\le
6\left(\frac{(T-t)^{3}}{6}-\sum_{j_1,j_2,j_3=0}^{q_1}
C_{j_3j_2j_1}^2\right)\ \ \ (i_1, i_2, i_3=1,\ldots,m).
\end{equation}

\vspace{6mm}

We can act similarly with more complicated 
iterated stochastic integrals. For example, for the 
approximation of stochastic integral
$I_{(0000)T,t}^{(i_1 i_2 i_3 i_4)}$ 
we can write (see (\ref{zzz10}))

\vspace{2mm}
$$
I_{(0000)T,t}^{(i_1 i_2 i_3 i_4)q_2}=
\sum_{j_1,j_2,j_3,j_4=0}^{q_2}
C_{j_4 j_3 j_2 j_1}\Biggl(
\prod_{l=1}^4\zeta_{j_l}^{(i_l)}
-\Biggr.
$$
$$
-
{\bf 1}_{\{i_1=i_2\}}
{\bf 1}_{\{j_1=j_2\}}
\zeta_{j_3}^{(i_3)}
\zeta_{j_4}^{(i_4)}
-
{\bf 1}_{\{i_1=i_3\}}
{\bf 1}_{\{j_1=j_3\}}
\zeta_{j_2}^{(i_2)}
\zeta_{j_4}^{(i_4)}-
$$
$$
-
{\bf 1}_{\{i_1=i_4\}}
{\bf 1}_{\{j_1=j_4\}}
\zeta_{j_2}^{(i_2)}
\zeta_{j_3}^{(i_3)}
-
{\bf 1}_{\{i_2=i_3\}}
{\bf 1}_{\{j_2=j_3\}}
\zeta_{j_1}^{(i_1)}
\zeta_{j_4}^{(i_4)}-
$$
$$
-
{\bf 1}_{\{i_2=i_4\}}
{\bf 1}_{\{j_2=j_4\}}
\zeta_{j_1}^{(i_1)}
\zeta_{j_3}^{(i_3)}
-
{\bf 1}_{\{i_3=i_4\}}
{\bf 1}_{\{j_3=j_4\}}
\zeta_{j_1}^{(i_1)}
\zeta_{j_2}^{(i_2)}+
$$
$$
+
{\bf 1}_{\{i_1=i_2\}}
{\bf 1}_{\{j_1=j_2\}}
{\bf 1}_{\{i_3=i_4\}}
{\bf 1}_{\{j_3=j_4\}}
+
{\bf 1}_{\{i_1=i_3\}}
{\bf 1}_{\{j_1=j_3\}}
{\bf 1}_{\{i_2=i_4\}}
{\bf 1}_{\{j_2=j_4\}}+
$$
$$
+\Biggl.
{\bf 1}_{\{i_1=i_4\}}
{\bf 1}_{\{j_1=j_4\}}
{\bf 1}_{\{i_2=i_3\}}
{\bf 1}_{\{j_2=j_3\}}\Biggr),
$$

\vspace{7mm}
\noindent
where $C_{j_4 j_3 j_2 j_1}$ is defined by (\ref{zzz11}), (\ref{zzz12}).
Moreover, according to (\ref{zzz0})

\vspace{2mm}
$$
{\sf M}\left\{\left(
I_{(0000)T,t}^{(i_1i_2 i_3 i_4)}-
I_{(0000)T,t}^{(i_1i_2 i_3 i_4)q_2}\right)^2\right\}\le
24\left(\frac{(T-t)^{4}}{24}-\sum_{j_1,j_2,j_3,j_4=0}^{q_2}
C_{j_4j_3j_2j_1}^2\right)\ \ \ (i_1, i_2, i_3, i_4=1,\ldots,m).
$$

\vspace{5mm}

For pairwise different $i_1, i_2, i_3, i_4=1,\ldots,m$ from 
(\ref{zzz4}) we obtain

\vspace{2mm}
\begin{equation}
\label{r7}
{\sf M}\left\{\left(
I_{(0000)T,t}^{(i_1i_2 i_3i_4)}-
I_{(0000)T,t}^{(i_1i_2 i_3i_4)q_2}\right)^2\right\}=
\frac{(T-t)^{4}}{24}
-\sum_{j_1,j_2,j_3,j_4=0}^{q_2}
C_{j_4j_3j_2j_1}^2.
\end{equation}

\vspace{5mm}

\noindent
\begin{figure}
\begin{center}
\centerline{Table 1.\ Coefficients $\bar C_{3j_2j_1}.$}
\vspace{4mm}
\begin{tabular}{|c|c|c|c|c|c|c|c|c|}
\hline
${}_{j_2} {}^{j_1}$&0&1&2&3&4&5&6\\
\hline
0&$0$&$\frac{2}{105}$&$0$&$-\frac{4}{315}$&$0$&$\frac{2}{693}$&0\\
\hline
1&$\frac{4}{105}$&0&$-\frac{2}{315}$&0&$-\frac{8}{3465}$&0&$\frac{10}{9009}$\\
\hline
2&$\frac{2}{35}$&$-\frac{2}{105}$&$0$&$\frac{4}{3465}$&
$0$&$-\frac{74}{45045}$&0\\
\hline
3&$\frac{2}{315}$&$0$&$-\frac{2}{3465}$&0&
$\frac{16}{45045}$&0&$-\frac{10}{9009}$\\
\hline
4&$-\frac{2}{63}$&$\frac{46}{3465}$&0&$-\frac{32}{45045}$&
0&$\frac{2}{9009}$&0\\
\hline
5&$-\frac{10}{693}$&0&$\frac{38}{9009}$&0&
$-\frac{4}{9009}$&0&$\frac{122}{765765}$\\
\hline
6&$0$&$-\frac{10}{3003}$&$0$&$\frac{20}{9009}$&$0$&$-\frac{226}{765765}$&$0$\\
\hline
\end{tabular}
\end{center}
\vspace{7mm}
\begin{center}
\centerline{Table 2.\ Coefficients $\bar C_{21j_2j_1}.$}
\vspace{4mm}
\begin{tabular}{|c|c|c|c|}
\hline
${}_{j_2} {}^{j_1}$&0&1&2\\
\hline
0&$\frac{2}{21}$&$-\frac{2}{45}$&$\frac{2}{315}$\\
\hline
1&$\frac{2}{315}$&$\frac{2}{315}$&$-\frac{2}{225}$\\
\hline
2&$-\frac{2}{105}$&$\frac{2}{225}$&$\frac{2}{1155}$\\
\hline
\end{tabular}
\end{center}
\vspace{7mm}
\begin{center}
\centerline{Table 3.\ Coefficients $\bar C_{101j_2j_1}.$}
\vspace{4mm}
\begin{tabular}{|c|c|c|}
\hline
${}_{j_2} {}^{j_1}$&0&1\\
\hline
0&$\frac{4}{315}$&$0$\\
\hline
1&$\frac{4}{315}$&$-\frac{8}{945}$\\
\hline
\end{tabular}
\end{center}
\vspace{5mm}
\end{figure}

Using Theorems 10, 11, we can calculate exactly the left-hand
side of (\ref{r7})
for any possible combinations
of $i_1, i_2, i_3, i_4$. These relations were obtained in 
\cite{20}-\cite{20a-new-x}, \cite{26}.

In Tables 1--3, we have some examples of exact 
values of the Fourier--Legendre coefficients 
(here and further in this article the Fourier--Legendre coefficients 
have been calculated exactly using Derive (computer algebra
system)). Note that in \cite{Kuz-Kuz}, \cite{Mikh-1} we used
the database with 270,000 exactly cal\-cu\-la\-ted
Fourier--Legendre coefficients. These coefficients \cite{Kuz-Kuz}, \cite{Mikh-1}
were calculated using the Python programming
language.

Assume that $q_1=6$.
Calculating the value of expression 
(\ref{39}) for $q_1=6,$ 
$i_1\ne i_2,$ $i_1\ne i_3,$ $i_3\ne i_2$,
we obtain 

\vspace{2mm}
$$
{\sf M}\left\{\left(
I_{(000)T,t}^{(i_1i_2 i_3)}-
I_{(000)T,t}^{(i_1i_2 i_3)q_1}\right)^2\right\}\approx
0.01956(T-t)^3.
$$

\vspace{5mm}

Let us choose, for example, $q_2=2.$ 
In the case of pairwise different
$i_1, i_2, i_3, i_4$ we have from (\ref{r7}) the following  
approximate equality

\vspace{2mm}
\begin{equation}
\label{46000}
{\sf M}\left\{\left(
I_{(0000)T,t}^{(i_1i_2i_3 i_4)}-
I_{(0000)T,t}^{(i_1i_2i_3 i_4)q_2}\right)^2\right\}\approx
0.0236084(T-t)^4.
\end{equation}

\vspace{5mm}

Consider the approximations 

\vspace{2mm}
$$
I_{(001)T,t}^{(i_1i_2i_3)q_3},\ \ \
I_{(010)T,t}^{(i_1i_2i_3)q_3},\ \ \
I_{(100)T,t}^{(i_1i_2i_3)q_3},\ \ \
I_{(00000)T,t}^{(i_1i_2i_3i_4 i_5)q_4}
$$

\vspace{4mm}
\noindent
based on the expansions (\ref{sss1})--(\ref{sss4}).

Assume that 
$q_3=2,$ $q_4=1.$  
In the case of pairwise different 
$i_1, \ldots, i_5$ we obtain

\vspace{2mm}
$$
{\sf M}\left\{\left(
I_{(100)T,t}^{(i_1i_2 i_3)}-
I_{(100)T,t}^{(i_1i_2 i_3)q_3}\right)^2\right\}=
\frac{(T-t)^{5}}{60}-\sum_{j_1,j_2,j_3=0}^{2}
\left(C_{j_3j_2j_1}^{100}\right)^2\approx 0.00815429(T-t)^5,
$$

\vspace{4mm}
$$
{\sf M}\left\{\left(
I_{(010)T,t}^{(i_1i_2 i_3)}-
I_{(010)T,t}^{(i_1i_2 i_3)q_3}\right)^2\right\}=
\frac{(T-t)^{5}}{20}-\sum_{j_1,j_2,j_3=0}^{2}
\left(C_{j_3j_2j_1}^{010}\right)^2\approx 0.0173903(T-t)^5,
$$

\vspace{4mm}
$$
{\sf M}\left\{\left(
I_{(001)T,t}^{(i_1i_2 i_3)}-
I_{(001)T,t}^{(i_1i_2 i_3)q_3}\right)^2\right\}=
\frac{(T-t)^5}{10}-\sum_{j_1,j_2,j_3=0}^{2}
\left(C_{j_3j_2j_1}^{001}\right)^2
\approx 0.0252801(T-t)^5,
$$

\vspace{6mm}
$$
{\sf M}\left\{\left(
I_{(00000)T,t}^{(i_1i_2i_3i_4 i_5)}-
I_{(00000)T,t}^{(i_1i_2i_3i_4 i_5)q_4}\right)^2\right\}=
$$

\vspace{2mm}
$$
=
\frac{(T-t)^5}{120}-\sum_{j_1,j_2,j_3,j_4,j_5=0}^{1}
C_{j_5j_4j_3j_2j_1}^2\approx 0.00759105(T-t)^5.
$$

\vspace{7mm}

Note that from (\ref{zzz0}) we have

\vspace{2mm}
$$
{\sf M}\left\{\left(
I_{(00000)T,t}^{(i_1i_2 i_3 i_4 i_5)}-
I_{(0000)T,t}^{(i_1i_2 i_3 i_4 i_5)q_4}\right)^2\right\}\le
120\left(\frac{(T-t)^{5}}{120}-\sum_{j_1,j_2,j_3,j_4,j_5=0}^{q_4}
C_{j_5j_4j_3j_2j_1}^2\right),
$$

\vspace{5mm}
\noindent
where $i_1, \ldots, i_5=1,\ldots,m$.

Let us consider the expansions of 
Ito stochastic integrals $I_{(1)T,t}^{(i_1)},$
$I_{(2)T,t}^{(i_1)}$ based on the Milstein approach from
\cite{Mi2}, which was mentioned in Sect.~1 
(also see \cite{KlPl2}, \cite{KPW})

\vspace{2mm}
\begin{equation}
\label{ww1}
I_{(1)T,t}^{(i_1)q}=-\frac{{(T-t)}^{3/2}}{2}
\left(\zeta_0^{(i_1)}-\frac{\sqrt{2}}{\pi}\left(\sum_{r=1}^{q}
\frac{1}{r}
\zeta_{2r-1}^{(i_1)}+\sqrt{\alpha_q}\xi_q^{(i_1)}\right)
\right),
\end{equation}

\vspace{3mm}
\begin{equation}
\label{ww2}
I_{(2)T,t}^{(i_1)q}=
(T-t)^{5/2}\left(
\frac{1}{3}\zeta_0^{(i_1)}+\frac{1}{\sqrt{2}\pi^2}
\left(\sum_{r=1}^{q}\frac{1}{r^2}\zeta_{2r}^{(i_1)}+
\sqrt{\beta_q}
\mu_q^{(i_1)}\right)-\right.
$$

\vspace{1mm}
$$
\left.-\frac{1}{\sqrt{2}\pi}\left(\sum_{r=1}^q
\frac{1}{r}\zeta_{2r-1}^{(i_1)}+
\sqrt{\alpha_q}\xi_q^{(i_1)}\right)\right),
\end{equation}

\vspace{5mm}
\noindent
where $\zeta_j^{(i)}$ is defined by the formula (\ref{rr23}),
$\phi_j(\tau)$ is a complete orthonormal system of trigonometric
functions in the space
$L_2([t, T]),$ and
$\zeta_0^{(i)},$ $\zeta_{2r}^{(i)},$
$\zeta_{2r-1}^{(i)},$ $\xi_q^{(i)},$ $\mu_q^{(i)}$ $(r=1,\ldots,q,$\
$i=1,\ldots,m)$ are independent 
standard Gaussian random variables, $i_1=1,\ldots,m,$ 

\vspace{1mm}
$$
\xi_q^{(i)}=\frac{1}{\sqrt{\alpha_q}}\sum_{r=q+1}^{\infty}
\frac{1}{r}\zeta_{2r-1}^{(i)},\ \ \
\alpha_q=\frac{\pi^2}{6}-\sum_{r=1}^q\frac{1}{r^2},
$$

\vspace{1mm}
$$
\mu_q^{(i)}=\frac{1}{\sqrt{\beta_q}}\sum_{r=q+1}^{\infty}
\frac{1}{r^2}~\zeta_{2r}^{(i)},\ \ \ 
\beta_q=\frac{\pi^4}{90}-\sum_{r=1}^q\frac{1}{r^4}.
$$

\vspace{4mm}

It is obvious that (\ref{ww1}), (\ref{ww2})
significantly more complicated compared to
(\ref{4002}), (\ref{4003}).

Another example of obvious advantage of the Legendre polynomials 
over the trigonometric functions (in the framework of the considered 
problem) is the truncated expansion 
of iterated Stratonovich stochastic integral
$I_{(10)T, t}^{*(i_1 i_2)}$ obtained by Theorem 5, in which
instead of the double Fourier--Legendre series (see (\ref{4004}),
(\ref{4006})) is taken 
the double trigonometric Fourier series

\vspace{1mm}
$$
I_{(10)T,t}^{*(i_1 i_2)q}=-(T-t)^{2}\Biggl(\frac{1}{6}
\zeta_{0}^{(i_1)}\zeta_{0}^{(i_2)}-\frac{1}{2\sqrt{2}\pi}
\sqrt{\alpha_q}\xi_q^{(i_2)}\zeta_0^{(i_1)}+\Biggr.
$$

\vspace{1mm}
$$
+\frac{1}{2\sqrt{2}\pi^2}\sqrt{\beta_q}\Biggl(
\mu_q^{(i_2)}\zeta_0^{(i_1)}-2\mu_q^{(i_1)}\zeta_0^{(i_2)}\Biggr)+
$$

\vspace{1mm}
$$
+\frac{1}{2\sqrt{2}}\sum_{r=1}^{q}
\Biggl(-\frac{1}{\pi r}
\zeta_{2r-1}^{(i_2)}
\zeta_{0}^{(i_1)}+
\frac{1}{\pi^2 r^2}\left(
\zeta_{2r}^{(i_2)}
\zeta_{0}^{(i_1)}-
2\zeta_{2r}^{(i_1)}
\zeta_{0}^{(i_2)}\right)\Biggr)-
$$

\vspace{1mm}
$$
-
\frac{1}{2\pi^2}\sum_{\stackrel{r,l=1}{{}_{r\ne l}}}^{q}
\frac{1}{r^2-l^2}\Biggl(
\zeta_{2r}^{(i_1)}
\zeta_{2l}^{(i_2)}+
\frac{l}{r}
\zeta_{2r-1}^{(i_1)}
\zeta_{2l-1}^{(i_2)}
\Biggr)+
$$

\vspace{1mm}
$$
+
\sum_{r=1}^{q}
\Biggl(\frac{1}{4\pi r}\left(
\zeta_{2r}^{(i_1)}
\zeta_{2r-1}^{(i_2)}-
\zeta_{2r-1}^{(i_1)}
\zeta_{2r}^{(i_2)}\right)+
$$

\vspace{1mm}
\begin{equation}
\label{944}
+
\Biggl.\Biggl.
\frac{1}{8\pi^2 r^2}\left(
3\zeta_{2r-1}^{(i_1)}
\zeta_{2r-1}^{(i_2)}+
\zeta_{2r}^{(i_2)}
\zeta_{2r}^{(i_1)}\right)\Biggr)\Biggr),
\end{equation}

\vspace{6mm}
\noindent
where the meaning of the notations included 
in (\ref{ww1}), (\ref{ww2}) is saved.

\vspace{5mm}

\section{Theorems 2--9 from Point
of View of the Wong--Zakai Approximation}

\vspace{5mm}

The iterated Ito stochastic integrals and solutions
of Ito SDEs are complex and important func\-ti\-o\-nals
from the independent components ${\bf f}_{s}^{(i)},$
$i=1,\ldots,m$ of the multidimensional
Wiener process ${\bf f}_{s},$ $s\in[0, T].$
Let ${\bf f}_{s}^{(i)p},$ $p\in\mathbb{N}$ 
be some approximation of
${\bf f}_{s}^{(i)},$
$i=1,\ldots,m$.
Suppose that 
${\bf f}_{s}^{(i)p}$
converges to
${\bf f}_{s}^{(i)},$
$i=1,\ldots,m$ if $p\to\infty$ in some sense and has
differentiable sample trajectories.

A natural question arises: if we replace 
${\bf f}_{s}^{(i)}$
by ${\bf f}_{s}^{(i)p},$
$i=1,\ldots,m$ in the functionals
mentioned above, will the resulting
functionals converge to the original
functionals from the components 
${\bf f}_{s}^{(i)},$
$i=1,\ldots,m$ of the multidimentional
Wiener process ${\bf f}_{s}$?
The answere to this question is negative 
in the general case. However, 
in the pioneering works of Wong E. and Zakai M. \cite{W-Z-1},
\cite{W-Z-2},
it was shown that under the special conditions and 
for some types of approximations 
of the Wiener process the answere is affirmative
with one peculiarity: the convergence takes place 
to the iterated Stratonovich stochastic integrals
and solutions of Stratonovich SDEs and not to iterated 
Ito stochastic integrals and solutions
of Ito SDEs.
The piecewise 
linear approximation 
as well as the regularization by convolution 
\cite{W-Z-1}-\cite{Watanabe} relate the 
mentioned types of approximations
of the Wiener process. The above approximation 
of stochastic integrals and solutions of SDEs 
is often called the Wong--Zakai approximation.

Let ${\bf w}_{\tau},$ $\tau\in[0, T]$ is a random vector with 
an $m+1$ components: ${\bf w}_{\tau}^{(i)}={\bf f}_{\tau}^{(i)}$ 
for $i=1,\ldots,m$ and 
${\bf w}_{\tau}^{(0)}=\tau,$\ 
${\bf f}_{\tau}^{(i)}$ $(i=1,\ldots,m)$
are independent standard Wiener processes.

It is well known that the following representation 
takes place \cite{Lipt}, \cite{7e}

\begin{equation}
\label{um1x}
{\bf w}_{\tau}^{(i)}-{\bf w}_{t}^{(i)}=
\sum_{j=0}^{\infty}\int\limits_t^{\tau}
\phi_j(s)ds\ \zeta_j^{(i)},\ \ \ \zeta_j^{(i)}=
\int\limits_t^T \phi_j(s)d{\bf w}_s^{(i)},
\end{equation}

\vspace{4mm}
\noindent
where $\tau\in[t, T],$ $t\ge 0,$
$\{\phi_j(x)\}_{j=0}^{\infty}$ is an arbitrary complete 
orthonormal system of functions in the space $L_2([t, T]),$ and
$\zeta_j^{(i)}$ are independent standard Gaussian 
random variables for various $i$ or $j.$
Moreover, the series (\ref{um1x}) converges for any $\tau\in [t, T]$
in the mean-square sense.

Let ${\bf w}_{\tau}^{(i)p}-{\bf w}_{t}^{(i)p}$ be 
the mean-square approximation of the process
${\bf w}_{\tau}^{(i)}-{\bf w}_{t}^{(i)},$
which has the following form

\vspace{-3mm}
\begin{equation}
\label{um1xx}
{\bf w}_{\tau}^{(i)p}-{\bf w}_{t}^{(i)p}=
\sum_{j=0}^{p}\int\limits_t^{\tau}
\phi_j(s)ds\ \zeta_j^{(i)}.
\end{equation}

\vspace{3mm}

From (\ref{um1xx}) we obtain

\vspace{-4mm}
\begin{equation}
\label{um1xxx}
d{\bf w}_{\tau}^{(i)p}=
\sum_{j=0}^{p}
\phi_j(\tau)\zeta_j^{(i)} d\tau.
\end{equation}

\vspace{4mm}

Consider the following iterated Riemann--Stieltjes
integral

\begin{equation}
\label{um1xxxx}
\int\limits_t^T
\psi_k(t_k)\ldots \int\limits_t^{t_2}\psi_1(t_1)
d{\bf w}_{t_1}^{(i_1)p_1}\ldots d{\bf w}_{t_k}^{(i_k)p_k},
\end{equation}

\vspace{4mm}
\noindent
where $p_1,\ldots,p_k\in \mathbb{N}$,\ \ $i_1,\ldots,i_k=0,1,\ldots,m,$ 

\begin{equation}
\label{um1xxx1}
d{\bf w}_{\tau}^{(i)p}=
\left\{\begin{matrix}
d{\bf f}_{\tau}^{(i)p}\ &\hbox{\rm for}\ \ \ i=1,\ldots,m\cr\cr\cr
d\tau^p\ &\hbox{\rm for}\ \ \ i=0
\end{matrix}
,\right.
\end{equation}

\vspace{4mm}
\noindent
and $d{\bf f}_{\tau}^{(i)p},$ $d\tau^p$ are defined by the relation (\ref{um1xxx}).

Let us substitute (\ref{um1xxx}) into (\ref{um1xxxx})

\begin{equation}
\label{um1xxxx1}
\int\limits_t^T
\psi_k(t_k)\ldots \int\limits_t^{t_2}\psi_1(t_1)
d{\bf w}_{t_1}^{(i_1)p_1}\ldots d{\bf w}_{t_k}^{(i_k)p_k}=
\sum\limits_{j_1=0}^{p_1}\ldots \sum\limits_{j_k=0}^{p_k}
C_{j_k \ldots j_1}\prod\limits_{l=1}^k \zeta_{j_l}^{(i_l)},
\end{equation}

\vspace{4mm}
\noindent
where 
$$
\zeta_j^{(i)}=\int\limits_t^T \phi_j(s)d{\bf w}_s^{(i)}
$$ 

\vspace{2mm}
\noindent
are independent standard Gaussian random variables for various 
$i$ or $j$ (in the case when $i\ne 0$),
${\bf w}_{s}^{(i)}={\bf f}_{s}^{(i)}$ for
$i=1,\ldots,m$ and 
${\bf w}_{s}^{(0)}=s,$

$$
C_{j_k \ldots j_1}=\int\limits_t^T\psi_k(t_k)\phi_{j_k}(t_k)\ldots
\int\limits_t^{t_2}
\psi_1(t_1)\phi_{j_1}(t_1)
dt_1\ldots dt_k
$$

\vspace{4mm}
\noindent
is the Fourier coefficient.

To best of our knowledge \cite{W-Z-1}-\cite{Watanabe}
the approximations of the Wiener process
in the Wong--Zakai approximation must satisfy fairly strong
restrictions
\cite{Watanabe}
(see Definition 7.1, pp.~480--481).
Moreover, approximations of the Wiener process that are
similar to (\ref{um1xx})
were not considered in \cite{W-Z-1}, \cite{W-Z-2}
(also see \cite{Watanabe}, Theorems 7.1, 7.2).
Therefore, the proof of analogs of Theorems 7.1 and 7.2 \cite{Watanabe}
for approximations of the Wiener 
process based on its series expansion (\ref{um1x})
should be carried out separately.
Thus, the mean-square convergence of the right-hand side
of (\ref{um1xxxx1}) to the iterated Stratonovich stochastic integral 
(\ref{str})
does not follow from the results of the papers
\cite{W-Z-1}, \cite{W-Z-2} (also see \cite{Watanabe},
Theorems 7.1, 7.2).

Nevetheless, the authors of the publications \cite{KlPl2}
(Sect.~5.8, pp.~202--204), \cite{KPS} (pp.~82-84),
\cite{KPW} (pp.~438-439),  
\cite{Zapad-9} (pp.~263-264) use (without rigorous proof)
the Wong--Zakai approximation 
\cite{W-Z-1}-\cite{Watanabe} based
on the series expansion of the Brownian bridge process \cite{Mi2}.

From the other hand, Theorems 2--9 from this 
paper can be considered as the proof of the
Wong--Zakai approximation for the iterated 
Stratonovich stochastic integrals (\ref{str}) of multiplicities 1 to 6
based on the approximation (\ref{um1xx}) of the Wiener process.
At that, the Riemann--Stieltjes integrals (\ref{um1xxxx}) converge
(according to Theorems 2--9)
to the appropriate Stratonovich 
stochastic integrals (\ref{str}). Recall that
$\{\phi_j(x)\}_{j=0}^{\infty}$ (see (\ref{um1x}), (\ref{um1xx}), and
Theorems 5--9)
is a complete 
orthonormal system of Legendre polynomials or 
trigonometric functions 
in the space $L_2([t, T])$.

To illustrate the above reasoning, 
consider two examples for the case $k=2,$
$\psi_1(s),$ $\psi_2(s)\equiv 1;$ $i_1, i_2=1,\ldots,m.$

The first example relates to the piecewise linear approximation
of the multidimensional Wiener process (these approximations 
were considered in \cite{W-Z-1}-\cite{Watanabe}).

Let ${\bf b}_{\Delta}^{(i)}(t),$ $t\in[0, T]$ be the piecewise
linear approximation of the $i$th component ${\bf f}_t^{(i)}$
of the multidimensional standard Wiener process ${\bf f}_t,$
$t\in [0, T]$ with independent components
${\bf f}_t^{(i)},$ $i=1,\ldots,m,$ i.e.

$$
{\bf b}_{\Delta}^{(i)}(t)={\bf f}_{k\Delta}^{(i)}+
\frac{t-k\Delta}{\Delta}\Delta{\bf f}_{k\Delta}^{(i)},
$$

\vspace{3mm}
\noindent
where 

\vspace{-2mm}
$$
\Delta{\bf f}_{k\Delta}^{(i)}={\bf f}_{(k+1)\Delta}^{(i)}-
{\bf f}_{k\Delta}^{(i)},\ \ \
t\in[k\Delta, (k+1)\Delta),\ \ \ k=0, 1,\ldots, N-1.
$$

\vspace{4mm}

Note that w.~p.~1

\vspace{-1mm}
\begin{equation}
\label{pridum}
\frac{d{\bf b}_{\Delta}^{(i)}}{dt}(t)=
\frac{\Delta{\bf f}_{k\Delta}^{(i)}}{\Delta},\ \ \
t\in[k\Delta, (k+1)\Delta),\ \ \ k=0, 1,\ldots, N-1.
\end{equation}

\vspace{4mm}

Consider the following iterated Riemann--Stieltjes
integral

\vspace{1mm}
$$
\int\limits_0^T
\int\limits_0^{s}
d{\bf b}_{\Delta}^{(i_1)}(\tau)d{\bf b}_{\Delta}^{(i_2)}(s),\ \ \ 
i_1,i_2=1,\ldots,m.
$$

\vspace{4mm}

Using (\ref{pridum}) and additive property of Riemann--Stieltjes integrals, 
we can write w.~p.~1

\vspace{2mm}
$$
\int\limits_0^T
\int\limits_0^{s}
d{\bf b}_{\Delta}^{(i_1)}(\tau)d{\bf b}_{\Delta}^{(i_2)}(s)=
\int\limits_0^T
\int\limits_0^{s}
\frac{d{\bf b}_{\Delta}^{(i_1)}}{d\tau}(\tau)d\tau
\frac{d {\bf b}_{\Delta}^{(i_2)}}{d s}(s)
ds =
$$

\vspace{4mm}
$$
=
\sum\limits_{l=0}^{N-1}\int\limits_{l\Delta}^{(l+1)\Delta}
\left(
\sum\limits_{q=0}^{l-1}\int\limits_{q\Delta}^{(q+1)\Delta}
\frac{\Delta{\bf f}_{q\Delta}^{(i_1)}}{\Delta}d\tau+
\int\limits_{l\Delta}^{s}
\frac{\Delta{\bf f}_{l\Delta}^{(i_1)}}{\Delta}d\tau\right)
\frac{\Delta{\bf f}_{l\Delta}^{(i_2)}}{\Delta}ds=
$$

\vspace{4mm}
$$
=\sum\limits_{l=0}^{N-1}\sum\limits_{q=0}^{l-1}
\Delta{\bf f}_{q\Delta}^{(i_1)}
\Delta{\bf f}_{l\Delta}^{(i_2)}+
\frac{1}{\Delta^2}\sum\limits_{l=0}^{N-1}
\Delta{\bf f}_{l\Delta}^{(i_1)}
\Delta{\bf f}_{l\Delta}^{(i_2)}
\int\limits_{l\Delta}^{(l+1)\Delta}
\int\limits_{l\Delta}^{s}d\tau ds=
$$

\vspace{4mm}
\begin{equation}
\label{oh-ty}
=\sum\limits_{l=0}^{N-1}\sum\limits_{q=0}^{l-1}
\Delta{\bf f}_{q\Delta}^{(i_1)}
\Delta{\bf f}_{l\Delta}^{(i_2)}+
\frac{1}{2}\sum\limits_{l=0}^{N-1}
\Delta{\bf f}_{l\Delta}^{(i_1)}
\Delta{\bf f}_{l\Delta}^{(i_2)}.
\end{equation}

\vspace{6mm}

Using (\ref{oh-ty}) and standard relation between Stratonovich
and Ito stochastic integrals, it 
is not difficult to show 
that

\vspace{1mm}
$$
\hbox{\vtop{\offinterlineskip\halign{
\hfil#\hfil\cr
{\rm l.i.m.}\cr
$\stackrel{}{{}_{N\to \infty}}$\cr
}} }
\int\limits_0^T
\int\limits_0^{s}
d{\bf b}_{\Delta}^{(i_1)}(\tau)d{\bf b}_{\Delta}^{(i_2)}(s)=
\int\limits_0^T
\int\limits_0^{s}
d{\bf f}_{\tau}^{(i_1)}d{\bf f}_{s}^{(i_2)}+
\frac{1}{2}{\bf 1}_{\{i_1=i_2\}}\int\limits_0^T ds=
$$

\vspace{4mm}
\begin{equation}
\label{uh-111}
=
{\int\limits_0^{*}}^T
{\int\limits_0^{*}}^s
d{\bf f}_{\tau}^{(i_1)}d{\bf f}_{s}^{(i_2)},
\end{equation}

\vspace{6mm}
\noindent
where $\Delta\to 0$ if $N\to\infty$ ($N\Delta=T$).

Obviously, (\ref{uh-111}) agrees with Theorem 7.1 (see \cite{Watanabe},
p.~486).

The next example relates to the approximation
of the Wiener process based on its series expansion
(\ref{um1x}) for $t=0$, where
$\{\phi_j(x)\}_{j=0}^{\infty}$ 
is a complete 
orthonormal system of Legendre polynomials or 
trigonometric functions 
in the space $L_2([0, T])$.

Consider the following iterated Riemann--Stieltjes
integral

\begin{equation}
\label{abcd1}
\int\limits_0^T
\int\limits_0^{s}
d{\bf f}_{\tau}^{(i_1)p}d{\bf f}_{s}^{(i_2)p},\ \ \ 
i_1,i_2=1,\ldots,m,
\end{equation}

\vspace{4mm}
\noindent
where $d{\bf f}_{\tau}^{(i)p}$ is defined by the
relation
(\ref{um1xxx}).

Let us substitute (\ref{um1xxx}) into (\ref{abcd1}) 

\vspace{-1mm}
\begin{equation}
\label{set18}
\int\limits_0^T
\int\limits_0^{s}
d{\bf f}_{\tau}^{(i_1)p}d{\bf f}_{s}^{(i_2)p}=
\sum\limits_{j_1,j_2=0}^p
C_{j_2 j_1} \zeta_{j_1}^{(i_1)}\zeta_{j_2}^{(i_2)},
\end{equation}

\vspace{3mm}
\noindent
where 
$$
C_{j_2 j_1}=
\int\limits_0^T \phi_{j_2}(s)\int\limits_0^s
\phi_{j_1}(\tau)d\tau ds
$$

\vspace{4mm}
\noindent
is the Fourier coefficient; another notations 
are the same as in (\ref{um1xxxx1}).

As we noted above, approximations of the Wiener process that are
similar to (\ref{um1xx})
were not considered in \cite{W-Z-1}, \cite{W-Z-2}
(also see Theorems 7.1, 7.2 in \cite{Watanabe}).
Furthermore, the extension of the results of Theorems 7.1 and 7.2
\cite{Watanabe} to the case under consideration is
not obvious.

On the other hand, we can apply the theory built in Chapters 1 and 2
of the monographs \cite{20a}-\cite{20a-new-x}. More precisely, 
using 
Theorem 5 for the case $k=2$,  
we obtain from (\ref{set18}) the desired result

\vspace{-1mm}
\begin{equation}
\label{umen-bl}
\hbox{\vtop{\offinterlineskip\halign{
\hfil#\hfil\cr
{\rm l.i.m.}\cr
$\stackrel{}{{}_{p\to \infty}}$\cr
}} }
\int\limits_0^T
\int\limits_0^{s}
d{\bf f}_{\tau}^{(i_1)p}d{\bf f}_{s}^{(i_2)p}=
\hbox{\vtop{\offinterlineskip\halign{
\hfil#\hfil\cr
{\rm l.i.m.}\cr
$\stackrel{}{{}_{p\to \infty}}$\cr
}} }
\sum\limits_{j_1,j_2=0}^p
C_{j_2 j_1} \zeta_{j_1}^{(i_1)}\zeta_{j_2}^{(i_2)}=
{\int\limits_0^{*}}^T
{\int\limits_0^{*}}^s
d{\bf f}_{\tau}^{(i_1)}d{\bf f}_{s}^{(i_2)}.
\end{equation}

\vspace{5mm}

From the other hand, by Theorems 2, 4
(see (\ref{a2})) for the case
$k=2$ we obtain from (\ref{set18}) the following relation

\vspace{-1mm}
$$
\hbox{\vtop{\offinterlineskip\halign{
\hfil#\hfil\cr
{\rm l.i.m.}\cr
$\stackrel{}{{}_{p\to \infty}}$\cr
}} }
\int\limits_0^T
\int\limits_0^{s}
d{\bf f}_{\tau}^{(i_1)p}d{\bf f}_{s}^{(i_2)p}=
\hbox{\vtop{\offinterlineskip\halign{
\hfil#\hfil\cr
{\rm l.i.m.}\cr
$\stackrel{}{{}_{p\to \infty}}$\cr
}} }
\sum\limits_{j_1,j_2=0}^p
C_{j_2 j_1} \zeta_{j_1}^{(i_1)}\zeta_{j_2}^{(i_2)}=
$$

\vspace{1mm}
$$
=
\hbox{\vtop{\offinterlineskip\halign{
\hfil#\hfil\cr
{\rm l.i.m.}\cr
$\stackrel{}{{}_{p\to \infty}}$\cr
}} }
\sum\limits_{j_1,j_2=0}^p
C_{j_2 j_1} \biggl(\zeta_{j_1}^{(i_1)}\zeta_{j_2}^{(i_2)}-
{\bf 1}_{\{i_1=i_2\}}{\bf 1}_{\{j_1=j_2\}}\biggr)+
{\bf 1}_{\{i_1=i_2\}}\sum\limits_{j_1=0}^{\infty}
C_{j_1 j_1}=
$$

\begin{equation}
\label{umen-blx}
=
\int\limits_0^T
\int\limits_0^{s}
d{\bf f}_{\tau}^{(i_1)}d{\bf f}_{s}^{(i_2)}+
{\bf 1}_{\{i_1=i_2\}}\sum\limits_{j_1=0}^{\infty}
C_{j_1 j_1}.
\end{equation}

\vspace{5mm}

Since
$$
\sum\limits_{j_1=0}^{\infty}
C_{j_1 j_1}=\frac{1}{2}\sum\limits_{j_1=0}^{\infty}
\left(\int\limits_0^T \phi_j(\tau)d\tau\right)^2
=\frac{1}{2}
\left(\int\limits_0^T \phi_0(\tau)d\tau\right)^2=\frac{1}{2}
\int\limits_0^T ds,
$$

\vspace{5mm}
\noindent
then from (\ref{umen-blx}) 
and standard relation between Stratonovich
and Ito stochastic integrals we obtain (\ref{umen-bl}).

\end{document}